\documentclass[12pt]{article}
\usepackage[cmtip,arrow]{xy}
\usepackage{pb-diagram, pb-xy}
\usepackage{amssymb,epsfig,amsfonts}
\newtheorem{defin}{Definition}
\newtheorem{prop}{Proposition}
\newtheorem{nt}{Remark}
\newtheorem{Th}{Theorem}
\newtheorem{lemma}{Lemma}
\newtheorem{cons}{Corollary}

\newtheorem{example}{Example}

\newfont{\ssdbl}{msbm8}
\newfont{\sdbl}{msbm9}
\newfont{\dbl}{msbm10 at 12pt}
\newcommand{\eqdef}{\stackrel{\rm def}{=}}
\newcommand{\proof}{{\bf Proof\ }}
\newcommand{\Ob}{\mathop {\rm Ob}}

\newcommand{\oo}{{\cal O}}
\newcommand{\ff}{{\cal F}}
\newcommand{\kk}{{\cal K}}

\newcommand{\Hom}{\mathop {\rm Hom}}
\newcommand{\Aut}{\mathop {\rm Aut}}
\newcommand{\Mor}{\mathop {\rm Mor}\nolimits}

\newcommand{\Lim}{\mathop {\rm lim}}

\newcommand{\Spec}{\mathop {\rm Spec}}

\newcommand{\Id}{\mathop {\rm Id}}
\newcommand{\rank}{\mathop {\rm rank}}

\newcommand{\da}{\mathbb{A}}
\newcommand{\dz}{\mathbb{Z}}
\newcommand{\dc}{\mathbb{C}}
\newcommand{\dq}{\mathbb{Q}}
\newcommand{\dr}{\mathbb{R}}
\newcommand{\dt}{\mathbb{T}}

\newcommand{\sdz}{\dz}
\newcommand{\sdr}{\dr}
\newcommand{\sdt}{\dt}
\newcommand{\sdc}{\mathbb{C}}
\newcommand{\ssdc}{\mathbb{C}}

\newcommand{\veps}{\varepsilon}

\newcommand{\Ker}{{\rm Ker}\:}

\newcommand{\f}{{\cal F}}
\newcommand{\lto}{\longrightarrow}

\newcommand{\dn}{\mathbb{N}}
\newcommand{\sdn}{\dn}
\newcommand{\df}{\mathbb{F}}

\newcommand{\cc}{\dc}
\newcommand{\F}{{\bf F}}
\newcommand{\I}{{\bf I}}
\newcommand{\D}{{\cal D}}

\newcommand{\s}{{\cal S}}

\newcommand{\oc}{\otimes_{\cc}}
\newcommand{\znakneponyaten}{{\bf 1}}

\def\Q{{\mathbb Q}}

\newcommand{\de}{{\bf 1}}
\newcommand{\ic}{\in \Ob(C_2^{ar})}
\newcommand{\nsubset}{\not{\subset}}

\begin{document}

\author{D.V. Osipov, A.N. Parshin
\footnote{Both authors are supported by RFBR (grant no.
08-01-00095-a), by a program of President of RF for supporting of
Leading Scientific Schools (grant no. NSh-4713.2010.1), by INTAS (grant
no. 05-1000008-8118); besides the first author is supported  by a
program of President of RF for supporting of young russian
scientists (grant no. MK-864.2008.1).} }

\title{Harmonic analysis on local fields and adelic
spaces II}
\date{}

\maketitle

\abstract{We develop harmonic analysis in  some categories of
filtered abelian groups and vector spaces. These categories contain
as objects local fields and adelic spaces arising from arithmetical
surfaces. Some structure theorems are proven for quotients of the
adelic groups of algebraic and arithmetical surfaces.}

\section{Introduction}

The paper is a sequel of the work~\cite{OsipPar}. The goal of this
part is two-fold.

First, we introduce and investigate a category $C_n^{fin}$ of
filtered abelian groups which is inductively constructed from the
category $C_0^n$ of finite abelian groups instead of the category
$C_n$ of filtered vector spaces over a field $k$ which was
inductively constructed  from the category $C_0$ of
finite-dimensional vector spaces over $k$ and was considered in
\cite{Osip}. Next, we consider some categories of filtered vector
spaces over $\dr$ or $\dc$ and define the mixed categories
$C_0^{ar}$, $C_1^{ar}$ and $C_2^{ar}$, which include as full
subcategories some of the previous categories. Then we develop the
harmonic analysis for them. One needs this part of our theory to
work with  an arithmetical surface. The resulting formalism is
completely parallel to the constructions of the first part,
i.e.~\cite{OsipPar}. We refer to the introduction in \cite{OsipPar}
for a short description of the formalism. Also, the origin and the
history of the whole theory is explained there~\footnote{We have to
mention once more the paper \cite{K} (2001) by M.~Kapranov, where
the important trick with a space of virtual measures was introduced
(the trick was inspired by the construction of semi-infinite monomes
in the theory of Sato Grassmanian, see a short description in
appendix 2 to the electronic version of \cite{Par2}). Kapranov
studied just the case of the Archimedean local fields such as
$\dr((t))$ or $\dc((t))$. As was pointed out to us by D.A.~Kazhdan,
he has given another  approach to harmonic analysis on the
two-dimensional local fields such as $\Q_p((t))$ in the paper:
Kazhdan D., {\em  Fourier transform over local fields}, Milan J.
Math. 74 (2006), 213--225. He has introduced the spaces of functions
and Fourier transforms, first, in an elementary coordinate-dependent
way and then redefined them not as a space (and a map) but as a
gerbe.}.

Second, we describe several examples of adelic spaces and show how
the general categorical notions can be applied for them. Note that
the paper~\cite{OsipPar} contained no one concrete example. Here, we
carefully study the case of adelic groups for schemes of dimension
$1$ (see example~\ref{exad} from section~\ref{c0c1} and
example~\ref{exam88} from section~\ref{subsec4.2}) and dimension $2$
(see section~\ref{lastexample} and example~\ref{arsur} from
section~\ref{c_2-spaces}). In particularly, we prove a compactness
theorem of the adelic quotient groups $\da_X/(\da_{X, 01} + \da_{X,
02})$ in the case of an projective algebraic surface $X$ defined
over a finite field (see sections~\ref{sqr} and~\ref{mrc}). This
theorem is a generalization of the classical compactness theorem for
the adelic quotient group $\da_C/\da_{C, 0}$ in the case of an
projective algebraic curve $C$ defined over a finite field.

Concerning the development of harmonic analysis on the groups such
as $\da^*_X = \mbox{GL}(1, \da_X)$ or the multiplicative group of a local field on $X$,
one has to take into account that these groups are extensions of
discrete groups by some non-locally compact groups. The harmonic
analysis for the discrete part is already very non-trivial
(see~\cite{Par4}).

We use the opportunity to correct the following statements
from~\cite{OsipPar}: corollary~1 and corollary~2  in~\cite[\S
5.9]{OsipPar}. We  have to consider more narrow class of
automorphisms. Namely, one needs to consider only the elements of
automorphism group (or its central extension) which
interchange the
elements of the first filtration of an object from the category
$C_2$\footnote{We note that this correction was already done in the last
electronic "arXiv" version of this paper: arXiv:0707.1766v3
[math.AG].}. This condition is satisfied in many concrete situations.
For example, one can consider the standard action of the multiplicative group
$K^*$ on a local field $K$ or the action of $\da_X^*$ on $\da_X$.

\section{Notations and agreements} \label{not}
For any two abelian groups $A$, $B$ we denote by $\Hom (A, B)$ the
abelian group of all homomorphisms from the group $A$ to the group
$B$.

By a topological group we always mean a topological group with a
Hausdorff topology.

For any locally compact abelian group $G$ by $\hat{G}$ we denote the
Pontryagin dual group of $G$, i.e. $\hat{G} = \Hom_{cont} (G, \dt)$,
where the group  $\dt= \{z \in \dc : |z|=1 \}$ is a one-dimensional
torus.

For any topological vector space $V$ over the field of complex
numbers $\dc$ we denote by $V'$ its continuous dual $\dc$-vector
space, i.e., $V' = \Hom_{\sdc, cont}(V, \dc)$ is the space of all
continuous $\dc$-linear functionals.

\section{Categories $C_n^{ar}$}
\subsection{Categories $C_n^{fin}$}

We recall that categories $C_n$, $n \in \dn$ were defined
in~\cite{Osip}. (See also sections~4.1 and~5.1 of~\cite{OsipPar} for
additional properties of categories $C_0$, $C_1$ and $C_2$.)

To construct categories $C_n$ we started by induction on $n$ from
category $C_0$ of  finite-dimensional vector spaces over a fixed
field $k$. As it was remarked in introduction to~\cite{Osip}, we can
construct categories $C_n^{fin}$ in a similar way starting from
category $C_0^{fin}$ of finite abelian groups. If the field $k$ is a
finite field, then for any $n \in \dn $ the category $C_n$ is a full
subcategory of the category $C_n^{fin}$.

Now we will give the definition of the category $C_n^{fin}$ by
induction on $n \in \dn$.

\begin{defin} \label{def1}
We say that $(I, F, V)$ is a filtered abelian group if
\begin{enumerate}
\item $V$ is an abelian group,
\item $I$ is a
partially ordered set, such that for any $i,j \in  I$ there are $k,l
\in I$ with $k \le i \le l$ and $k \le j \le l$,
\item $F$ is a
function from $I$ to the set of  subgroups of $V$ such that if $i
\le j $ are any from $I$, then $F(i) \subset F(j)$,
\item $ \bigcap\limits_{i \in I}  F(i) = 0$  and $
\bigcup\limits_{i \in I} F(i) = V $.
\end{enumerate}
\end{defin}

\begin{defin} \label{defdom}
We say that a filtered abelian group $(I_1, F_1, V)$ dominates
another filtered abelian group $(I_2, F_2, V)$ when there is an
order-preserving function $\phi: I_2 \to I_1$ such that
\begin{enumerate}
\item for any $i \in I_2$ we have $F_1(\phi (i)) = F_2 (i)$
\item for any $j \in I_1$ there are $i_1, i_2 \in I_2$ such that $ \phi(i_1) \le j \le \phi(i_2)$.
\end{enumerate}
\end{defin}

Now we define by induction the category $C_n^{fin}$ and morphisms
between them.

\begin{defin}
\begin{enumerate}
 \item The category  $C_0^{fin}$ is  the category of finite abelian groups
 with morphisms of abelian groups.
 \item The triple from $C_0^{fin}$
 $$0 \lto V_0 \lto V_1  \lto V_2 \lto  0 $$ is admissible
 when it is an exact triple of abelian groups \mbox{.}
 \end{enumerate}
\end{defin}

Now we define the objects of the  category $C_n^{fin}$ by induction.
We suppose that we have already defined the objects of the category
$C_{n-1}^{fin}$ and the notion of admissible triple in
$C_{n-1}^{fin}$.

\begin{defin} \label{def4}
\begin{enumerate}
\item Objects of the category $C_n^{fin}$, i.e. $Ob(C_n^{fin})$, are filtered abelian groups $(I, F, V)$
with the following additional structures:
\begin{enumerate}
\item for any $i \le j \in I$  on the abelian group $F(j) / F(i)$
it is given a structure $E_{i,j} \in Ob(C_{n-1}^{fin})$,
\item for any $i \le j \le k \in I$
$$
0 \lto E_{i,j}  \lto E_{i,k}  \lto  E_{j,k} \lto 0
$$
is an admissible triple from $C_{n-1}^{fin}$.
\end{enumerate}
\item Let $E_1 = (I_1, F_1, V_1)$, $E_2 = (I_2, F_2, V_2)$ and $E_3=(I_3, F_3, V_3)$
be from  $Ob(C_n^{fin})$. Then we say that
$$
0 \lto E_1 \lto E_2 \lto E_3 \lto 0
$$
is an admissible triple from $C_n^{fin}$ when
 the following conditions are satisfied:
\begin{enumerate}
\item
$$
0 \lto V_1 \lto V_2 \lto V_3 \lto 0
$$
is an exact triple of abelian groups,
\item \label{itaa}
the filtration $(I_1, F_1, V_1)$ dominates the filtration $(I_2,
F'_1, V_1)$, where $F'_1 (i) = F_2(i) \cap V_1$  for any $i \in
I_2$,
\item \label{itbb}
the filtration $(I_3, F_3, V_3)$ dominates the filtration $(I_2,
F'_3, V_3)$, where $F'_3(i) =  F_2(i) / F_2(i) \cap V_1$,
\item for any $i \le j \in I_2$
\begin{equation} \label{trojkaa1}
0 \lto \frac{F'_1(j)}{F'_1(i)} \lto \frac{F_2(j)}{F_2(i)}  \lto
\frac{F'_3(j)}{F'_3(i)} \lto 0
\end{equation}
is an admissible triple from $C_{n-1}^{fin}$. (By definition of
$Ob(C_n^{fin})$, on every abelian group from triple~(\ref{trojkaa1})
it is given the structure of $Ob(C_{n-1}^{fin})$).
\end{enumerate}
\end{enumerate}
\end{defin}

By induction, we define now  the morphisms in the category
$C_n^{fin}$. We suppose that we have already defined  the morphisms
in $C_{n-1}^{fin}$.

\begin{defin} \label{d1}
Let $E_1 = (I_1, F_1, V_1)$ and $E_2 = (I_2, F_2, V_2)$ be from
$Ob(C_n^{fin})$. Then $\Mor_{C_n^{fin}}(E_1, E_2)$ consists of
elements $A \in \Hom (V_1, V_2)$ such that the following conditions
hold:
\begin{enumerate}
\item \label{i1} for any $i \in I_1$ there is an $j \in I_2$ such that $A (F_1(i)) \subset F_2(j)$,
\item \label{i2}  for any $j \in I_2$ there is an $i \in I_1$ such that $A (F_1(i)) \subset F_2(j)$,
\item \label{i3}  for any $i_1 \le i_2 \in I_1$ and $j_1 \le j_2 \in I_2$ such that $A (F_1(i_1)) \subset F_2(j_1)$
and $A (F_1(i_2)) \subset F_2(j_2)$ we have that the induced  map of
abelian groups
$$   \bar{A} : \frac{F_1(i_2)}{F_1(i_1)} \lto \frac{F_2(j_2)}{F_2(j_1)}
$$
is an element from
$$\Mor_{C_{n-1}^{fin}}(\frac{F_1(i_2)}{F_1(i_1)}, \frac{F_2(j_2)}{F_2(j_1)}) \mbox{.}$$
\end{enumerate}
\end{defin}

The correctness of the definition of morphisms in categories
$C_n^{fin}$ and some other properties follow from the following
proposition, which is an analog of proposition~2.1 from~\cite{Osip}.
\begin{prop} \label{propmor}
 Let $E_1= (I_1, F_1, V_1)$, $E_2= (I_2, F_2,
V_2)$, $E_1'$, $E_2'$ be from $Ob(C_n^{fin})$, and  $A$ is a  map
from $\Hom (V_1, V_2)$.
\begin{enumerate}
\item If the filtered abelian group $E_1$ dominates the filtered abelian group
$E'_1$, and the filtered abelian group $E_2$ dominates the filtered
abelian group $E'_2 $, then $A \in \Mor_{C_n^{fin}}(E_1, E_2)$ if
and only if $A \in \Mor_{C_n^{fin}}(E'_1, E'_2)$.
\item $\Mor_{C_n^{fin}}(E_1, E_2)$ is an abelian subgroup of $\Hom(V_1, V_2)$.
\item If $E_3$ is an object of $C_n^{fin}$, then
$$     \Mor_{C_n^{fin}}(E_2, E_3)   \circ   \Mor_{C_n^{fin}}(E_1, E_2)  \subset  \Mor_{C_n^{fin}}(E_1, E_3) \mbox{.}$$
\end{enumerate}
\end{prop}
\proof of this proposition is fully analogous to the proof of
proposition~2.1 from~\cite{Osip}.

\begin{example} \label{exam1}{\em
Let $X$ be an $n$-dimensional scheme of finite type over $\dz$.

Let $P(X)$ be the set of points of the scheme $X$. We consider
$\eta, \nu \in P(X)$. We define $\eta \ge \nu$ if $\nu \in \bar{\{
\eta\}}$. The relation $\ge$ is a partial order on $P(X)$. Let
$S(X)$ be the simplicial set induced by $(P(X), \ge)$, i.e.,
$$S(X)_m = \{ (\nu_0, \ldots, \nu_m) \mid \nu_i \in P(X); \nu_i \ge \nu_{i+1} \}$$
is the set of $m$-simplices of $S(X)$ with the $i$-th boundary map,
which is the deleting of the $i$-th point, and the $i$-th degeneracy
map, which is the duplication of the $i$-th point.  By
$S(X)_m^{red}$ we denote the set of all non-degenerated simplices
from $S(X)_m$, i.e.,
$$S(X)_m^{red} = \{ (\nu_0, \ldots, \nu_m) \mid \nu_i \in P(X); \; \nu_i < \nu_{i+1} \} \mbox{.}$$

In~\cite{B},~\cite{H} for any subset $K \subset S(X)_m$ it was
constructed the  adelic space $\da(K, \f)$ for any quasicoherent
sheaf $\f$ on $X$ such that
$$
\da(K, \f) \subset \prod_{\delta \in K} \da(\delta, \f) \mbox{.}
$$

In a completely analogous way as the proof of theorem~1
from~\cite{Osip}, we have that for any subset $K \subset
S(X)_n^{red}$ the adelic space $\da(K, \f)$ is an object of the
category $C_n^{fin}$.

In particularly, we obtain that if $K= S(X)_n^{red}$, $\f = \oo_X$,
then the whole adelic space $\da_X= \da(S(X)_n^{red}, \oo_X)$ is an
object of the category $C_n^{fin}$.

If $\delta = (\eta_0 > \ldots > \eta_n) \in S(X)_n^{red}$ such that
$\eta_n$ is a regular point on every subscheme $\bar{\eta_i}$, which
is closure of the point $\eta_i$, then $K_{\delta} = \da(\delta,
\oo_V)$ is an $n$-dimensional local field with the finite last
residue field, see~\cite{PF}, \cite{Osi}. We obtain that $K_{\delta}
\in \Ob(C_n^{fin})$.

In particularly, the fields $$\df_q((t_n)) \ldots ((t_1)) \in
\Ob(C_n^{fin}) \mbox{,}$$
$$E((t)) \in \Ob(C_2^{fin}) \mbox{,}$$
$$E\{\{t\}\} \in \Ob(C_2^{fin}) \mbox{,}$$
where $E \supset \dq_p$ is a finite field extension.
}\end{example}

\subsection{Categories $C_0^{ar}$ and $C_1^{ar}$} \label{c0c1}
We define the category $C_0^{ar}$ in the following way.
\begin{defin}
The category $C_0^{ar}$ is the full subcategory of the category of
commutative finite-dimensional smooth real  Lie groups such that
$$
G \in \Ob(C_0^{ar}) \qquad \mbox{iff} \qquad   \pi_0(G) \quad
\mbox{is finitely generated,}
$$
where the abelian group $\pi_0(G)$ is the group of connected
componets of the group $G$.
\end{defin}

We note the following well-known property.
\begin{lemma}
A group $G \in \Ob(C_0^{ar})$ if and only if there is the following
isomorphism:
\begin{equation} \label{decomp}
G \simeq A \times \dz^{r} \times  \dt^{p} \times \dr^{q} \mbox{,}
\end{equation}
where $p, q, r$ are non-negative integers, $A$ is a finite abelian
group, $\dt= S^1 = \dr/\dz$ is a one-dimensional torus.
\end{lemma}
\proof If a group $G$ satisfies  decomposition~(\ref{decomp}), then
it is obvious that $G \in \Ob(C_0^{ar})$. Now  let $G \in
\Ob(C_0^{ar})$, and $G^0 \subset G$ be the connected component of
the identity element $e \in G$. Since $G$ is an abelian group, and
$G^0$ is connected, we have that  the exponential map is a
surjective homomorphism from the tangent space of $e \in G$ (which
is a vector Lie group) to $G^0$. Therefore, we have that $G^0 \simeq
\dt^{p} \times \dr^{q}$ for some integers $p,q$. Now the discrete
group
$$G/G^0 \simeq \dz^r \times \prod\limits_{1 \le i \le l} \dz/n_i\dz
$$
for some integers $n_i$ $( 1 \le i \le l)$ and $r$. We choose any $1
\le i \le l$, and  let $e_i \in G$ be a preimage of a generator of
the group $\dz/n_i\dz$. Then an element $s_i=e_i^{n_i} \in G^0$. The
group $G^0 \simeq \dt^{p} \times \dr^{q}$, therefore there is an
element $t_i \in G^0$ such that $t_i^{n_i}=s_i$. Now for an element
$e_i'=e_it_i^{-1}$ we have $(e_i')^{n_i}=1$. Therefore the subgroup
in $G$ generated by the element $e_i'$ is isomorphic to the group
$\dz/n_i\dz$. We obtained the lifting of the group $\dz/n_i\dz$ to
the group $G$ for any $1 \le i \le l$. Now by lifting the generators
of the group $\dz^r$ to their preimages in the group $G$ we
isomorphically map the group $\dz^r$ into the group $G$. The lemma
is proved.

\vspace{0.5cm}

Decomposition~(\ref{decomp}) is not unique. But we have always the
following canonical filtration of the group $G$ by Lie subgroups:
\begin{equation} \label{dee}
G \supset G_{tor} \supset G^0 \supset K  \mbox{,}
\end{equation}
where the group $G^0$ is the connected component of the group $G$
which contains the identity element $e$, the group $K$ is a maximal
compact subgroup of the group $G^0$, the group $G_{tor}/G^0$ is the
maximal torsion subgroup of the discrete group $G / G^0$. Then we
have the following isomorphisms
$$
K \simeq \dt^{p} \quad \mbox{,} \quad G^0/ K \simeq \dr^q \mbox{,}
$$
$$
G/G^0 \quad \mbox{is finitely generated abelian group,}
$$
$$
G_{tor}/G^0 \simeq A \quad \mbox{,} \quad G/G_{tor} \simeq \dz^{r}
\mbox{.}
$$
Therefore the non-negative integers $p,q,r$ from
decomposition~(\ref{decomp}) are defined unique by the group $G$,
the finite group $\Psi$ from decomposition~(\ref{decomp}) is also
defined by the group $G$ unique up to an isomorphism. Besides, we
have that
$$\pi_0(G)=G/G^0 \quad \mbox{,} \quad \rank(\pi_1(G))= p \quad \mbox{,} \quad \dim G= p+q \mbox{.}$$

\begin{defin} \label{adtr}
We say that the triple from $C_0^{ar}$
\begin{equation} \label{adtreq}
0 \arrow{e} G_1 \arrow{e} G_2 \arrow{e,t}{\phi} G_3 \arrow{e} 0
\end{equation}
is an admissible triple from $C_0^{ar}$ if $G_1 = \Ker(\phi)$, and
$\phi$ is an epimorphism of Lie groups.
\end{defin}
\begin{nt} \label{Cartan} {\em
Definition~\ref{adtr} can be reformulated in the following way: the
triple~(\ref{adtreq}) is an admissible triple if and only if it is
isomorphic to the following triple:
$$
0 \arrow{e} G_1 \arrow{e} G_2 \arrow{e} G_2/G_1 \arrow{e} 0 \mbox{,}
$$
where $G_1$ is a closed Lie subgroup of a Lie group $G_2$. According
to Cartan's theorem, it is equivalent to the fact that $G_1$ is a
 closed subgroup of a Lie group $G_2$.}
\end{nt}

\begin{example} \label{example11} {\em
The category $C_0^{fin}$ is the full subcategory of the category
$C_0^{ar}$ such that a triple from $C_0^{fin}$ is an admissible
triple from $C_0^{fin}$ if and only if it is an admissible triple in
$C_0^{ar}$.

For any field $k$ by $C_n(k)$ we denote the category $C_n$ defined
over the field $k$ (see~\cite{Osip}). We have that the category
$C_0(\dr)$ is the full subcategory of the category $C_0^{ar}$ such
that a triple from $C_0(\dr)$ is an admissible triple from
$C_0(\dr)$ if and only if it is an admissible triple in $C_0^{ar}$.
}
\end{example}

Now we define the categories $C_1^{ar}$.
\begin{defin} \label{def4'}
Objects of the category $C_1^{ar}$, i.e. $Ob(C_1^{ar})$, are
filtered abelian groups $E=(I, F, V)$ (see definition~\ref{def1})
with the following additional structures and conditions
\begin{enumerate}
\item for any $i \le j \in I$  on the abelian group $F(j) / F(i)$
it is given a structure $E_{i,j} \in Ob(C_{0}^{ar})$;
\item for any $i \le j \le k \in I$
$$
0 \lto E_{i,j}  \lto E_{i,k}  \lto  E_{j,k} \lto 0
$$
is an admissible triple from $C_{0}^{ar}$;
\item  \label{condition3} there exist $i_E \le j_E \in I$ such that for any $k \le
i_E$ the group $F(i_E)/F(k)$ is a finite abelian group, and for any
$l \ge j_E$ the group $F(l) / F(j_E)$ is a finite abelian group.
\end{enumerate}
\end{defin}
\begin{nt} \label{remark2}{\em Condition~{\ref{condition3}} of
definition~\ref{def4'} is equivalent to the following condition:
there is $d_E \in \dn$ such that for any $i \le j \in I$ the group
$F(j)/F(i)$ as an object of the category $C_0^{ar}$ has the
following properties
$$
\dim F(j)/F(i) \le d_E  \qquad \mbox{,} \qquad \rank \pi_0
(F(j)/F(i)) \le d_E  \mbox{.}
$$
}
\end{nt}

\begin{defin} \label{dd1}
Let $E_1 = (I_1, F_1, V_1)$ and $E_2 = (I_2, F_2, V_2)$ be from
$Ob(C_1^{ar})$. Then $\Mor_{C_1^{ar}}(E_1, E_2)$ consists of
elements $A \in \Hom (V_1, V_2)$ such that the following conditions
hold:
\begin{enumerate}
\item \label{ii1} for any $i \in I_1$ there is an $j \in I_2$ such that $A (F_1(i)) \subset F_2(j)$,
\item \label{ii2}  for any $j \in I_2$ there is an $i \in I_1$ such that $A (F_1(i)) \subset F_2(j)$,
\item \label{ii3}  for any $i_1 \le i_2 \in I_1$ and $j_1 \le j_2 \in I_2$ such that $A (F_1(i_1)) \subset F_2(j_1)$
and $A (F_1(i_2)) \subset F_2(j_2)$ we have that the induced  map of
abelian groups
$$   \bar{A} : \frac{F_1(i_2)}{F_1(i_1)} \lto \frac{F_2(j_2)}{F_2(j_1)}
$$
is an element from
$$\Mor_{C_0^{ar}}(\frac{F_1(i_2)}{F_1(i_1)}, \frac{F_2(j_2)}{F_2(j_1)}) \mbox{.}$$
\end{enumerate}
\end{defin}

Obviously, we have the full analog of  proposition~\ref{propmor} for
categories $C_1^{ar}$.
\begin{prop} \label{propmor2}
 Let $E_1= (I_1, F_1, V_1)$, $E_2= (I_2, F_2,
V_2)$, $E_1'$, $E_2'$ be from $Ob(C_1^{ar})$, and  $A$ is a  map
from $\Hom (V_1, V_2)$.
\begin{enumerate}
\item If the filtered abelian group $E_1$ dominates the filtered abelian group
$E'_1$, and the filtered abelian group $E_2$ dominates the filtered
abelian group $E'_2 $, then $A \in \Mor_{C_1^{ar}}(E_1, E_2)$ if and
only if $A \in \Mor_{C_1^{ar}}(E'_1, E'_2)$.
\item $\Mor_{C_1^{ar}}(E_1, E_2)$ is an abelian subgroup of $\Hom(V_1, V_2)$.
\item If $E_3$ is an object of $C_1^{ar}$, then
$$     \Mor_{C_1^{ar}}(E_2, E_3)   \circ   \Mor_{C_1^{ar}}(E_1, E_2)  \subset  \Mor_{C_1^{ar}}(E_1, E_3) \mbox{.}$$
\end{enumerate}
\end{prop}
\proof of this proposition is fully analogous to the proof of
proposition~2.1 from~\cite{Osip}.

\begin{defin}
\begin{enumerate}
\item Let $E_1 = (I_1, F_1, V_1)$, $E_2 = (I_2, F_2, V_2)$ and $E_3=(I_3, F_3, V_3)$
be from  $Ob(C_1^{ar})$. Then we say that
$$
0 \lto E_1 \lto E_2 \lto E_3 \lto 0
$$
is an admissible triple from $C_1^{ar}$ when
 the following conditions are satisfied:
\begin{enumerate}
\item
$$
0 \lto V_1 \lto V_2 \lto V_3 \lto 0
$$
is an exact triple of abelian groups
\item \label{itaa1}
the filtration $(I_1, F_1, V_1)$ dominates the filtration $(I_2,
F'_1, V_1)$, where $F'_1 (i) = F_2(i) \cap V_1$  for any $i \in
I_2$,
\item \label{itbb1}
the filtration $(I_3, F_3, V_3)$ dominates the filtration $(I_2,
F'_3, V_3)$, where $F'_3(i) =  F_2(i) / F_2(i) \cap V_1$,
\item for any $i \le j \in I_2$
\begin{equation} \label{trojkaa1'}
0 \lto \frac{F'_1(j)}{F'_1(i)} \lto \frac{F_2(j)}{F_2(i)}  \lto
\frac{F'_3(j)}{F'_3(i)} \lto 0
\end{equation}
is an admissible triple from $C_{0}^{ar}$. (By definition of
$Ob(C_1^{ar})$, on every abelian group from triple~(\ref{trojkaa1'})
it is given the structure of $Ob(C_{0}^{ar})$).
\end{enumerate}
\end{enumerate}
\end{defin}

\begin{example} \label{examp1} {\em
The category $C_1^{fin}$ is a full subcategory of the category
$C_1^{ar}$ such that a triple from $C_1^{fin}$ is an admissible
triple from $C_1^{fin}$ if and only if it is an admissible triple in
$C_1^{ar}$.

The category $C_0^{ar}$ is a full subcategory of the category
$C_1^{ar}$ with respect to the following functor $I$: for any $G \in
C_0^{ar}$ we put $I(G)= (\{(0), (1)\}, F, G )$, where $(0) < (1)$,
$F((0))=e$, $F((1))=G$, $e$ is the trivial subgroup of $G$, and $I$
acts as identity map on the morphisms from $C_0^{ar}$. }
\end{example}

Let $G \in \Ob(C_0^{ar})$. Then $G$ is a locally compact abelian
group. Therefore its Pontryagin dual group $\hat{G}$ is also a
locally compact abelian group. Moreover, by formula~(\ref{decomp}),
we have $\hat{G} \in \Ob(C_0^{ar})$, because for any $n \in \dn$
 $$\widehat{\dz/n\dz} \simeq
\dz/n\dz \mbox{,} \qquad \hat{\dt} \simeq \dz \mbox{,} \qquad
\hat{\dz} \simeq \dt \mbox{,} \qquad \hat{\dr} \simeq \dr \mbox{.}$$

Let $E = (I, F, V) \in \Ob(C_1^{ar})$. We define  the dual object
$\check{E} = (I^0, F^0, \check{V}) \in \Ob(C_1^{ar})$ in the
following way. The abelian group $\check{V} \subset \Hom(V, \dt)$ is
defined as
\begin{equation} \label{f6}
\check {V} \eqdef \mathop{\Lim_{\rightarrow}}_{j \in I}
\mathop{\Lim_{\leftarrow}}_{i \ge j} \widehat{F(i)/F(j)} \mbox{.}
\end{equation}
The set $I^0$ is a partially ordered set, which has the same set as
$I$, but with the inverse order then $I$. For $j \in I^0$ we define
the subgroup
\begin{equation} \label{f7}
F^0(j) \eqdef \mathop{\Lim_{\leftarrow}}_{i \le j \in I^0}
\widehat{F(i)/F(j)} \quad \subset \quad \check{V} \mbox{.}
\end{equation}

If $E_1, E_2 \in \Ob(C_1^{ar})$ and $\theta \in
\Mor_{C_1^{ar}}(E_1,E_2)$, then, by definition, there is canonically
$\check{\theta} \in \Mor_{C_1^{ar}} (\check{E_2}, \check{E_1})$. If
$$
0 \arrow{e} E_1  \arrow{e,t}{\alpha}  E_2 \arrow{e,t}{\beta}   E_3
\arrow{e} 0
$$
is an admissible triple from $C_1^{ar}$, where $E_i = (I_i, F_i,
V_i)$, $ 1 \le i \le 3$, then  we have the following canonical
admissible triple from $C_1^{ar}$:
$$
0 \arrow{e} \check{E_3} \arrow{e,t}{\check{\beta}} \check{E_2}
\arrow{e,t}{\check{\alpha}}  \check{E_1} \arrow{e} 0 \mbox{.}
$$

\begin{defin}
We say that an object $(I,F,V) \in C_1^{ar}$ is a complete object if
the following condition is satisfied:
\begin{equation}
V = \mathop{\Lim_{\to}}_{i \in I} \mathop{\Lim_{\gets}}_{j \le i}
F(i)/F(j)  \mbox{.}
\end{equation}
\end{defin}

For any locally compact abelian group $G$ we have that
$\hat{\hat{G}} = G$. Therefore it follows from definitions that $E
\in \Ob(C_1^{ar})$ is a complete object  if and only if
$\check{\check{E}} = E$.

We have the obvious functor of completion  $\Psi : C_1^{ar} \lto
C_1^{ar}$, where for $E = (I,F,V) \in \Ob(C_1^{ar})$ we define
$$
\Psi(E) \eqdef (I, F', \mathop{\Lim_{\to}}_{i \in I}
\mathop{\Lim_{\gets}}_{j \le i} F(i)/F(j)) \mbox{,}
$$
and where for any $j \in I$ we define
$$
F'(j) \eqdef \mathop{\Lim_{\gets}}_{j \le i} F(i)/F(j) \mbox{.}
$$
$\Psi$ is easily extended to $\Mor_{C_1^{ar}}(E_1, E_2)$, $E_1, E_2
\in \Ob(C_1^{ar})$ from $\Mor_{C_0^{ar}}$.

It is clear from definition that for any $E \in \Ob(C_1^{ar})$ the
object $\Psi(E)$ is a complete object, and $\Psi^2 = \Psi$.
Moreover, $\Psi(E)= \check{\check{E}}$.

By $C_{1, compl}^{ar}$ we denote the full subcategory of the
category $C_1^{ar}$ that consists of complete objects of the
category $C_1^{ar}$.

Let $Loc$ be the category of locally compact abelian groups. For any
group $G \in \Ob(Loc)$ there is a pair of subgroups $V \subset U$ in
$G$ with the following properties (see~\cite[ch. II, \S 2.2]{Bo1}):
\begin{description}
\item[\rm (*)] $U$ is an open subgroup in $G$, and $U$ is generated by an
compact neighbourhood of the identity element $e$ (neighbourhood of
$e$  is a set which contains an open neighbourhood of $e$);
\item[\rm (**)] $V$ is a compact subgroup in $G$;
\item[\rm (***)] $U /V \in \Ob(C_0^{ar})$.
\end{description}
Moreover, this $G$ is an inductive limit of subgroups $U$, and every
$U$ is a projective limit of groups $U/V$, where $U$ and $V$ satisfy
properties (*) -- (***) above, and we consider on $G$ the topology
of inductive and projective limits starting from topological groups
$U/V$ (see~\cite[ch. II, \S 2.2]{Bo1}).

By $Loc^{ar}$ we denote the following full subcategory of the
category $Loc$. Let a group $G \in \Ob(Loc)$. Then the group $G \in
\Ob(Loc^{ar})$ if the following condition is satisfied: there is $d
\in \dn$ such that for any pairs of subgroups $V \subset U$ in $G$
satisfying properties (*) -- (***) above, we have   that
\begin{equation} \label{dimrank}
\dim (U/V) < d \qquad \mbox{,} \qquad \rank \pi_0(U/V) < d \mbox{.}
\end{equation}
(Compare this condition with
remark~\ref{remark2}).

Now we define the functor $\Phi \: : \:  C_{1}^{ar} \to Loc^{ar}$ in
the following way. Let $E = (I, F, V) \in \Ob(C_{1}^{ar})$. We
define
$$
\Phi(E) \eqdef  \mathop{\Lim_{\to}}_{i \in I}
\mathop{\Lim_{\gets}}_{j \le i} F(i)/F(j) \mbox{,}
$$
where on the group $\Phi(E)$ we consider the topology of inductive
and projective limits starting from topological groups $F(i)/F(j)
\in \Ob(C_0^{ar})$. (Besides, the condition from
remark~\ref{remark2} goes to  condition (\ref{dimrank}).) $\Phi$ is
defined also on morphisms, because we can glue compatible system of
morphisms from $C_0^{ar}$, which is given by definition~\ref{dd1}.
We have the following proposition.
\begin{prop} \label{propos9}
\begin{enumerate}
\item \label{st1} $\Phi \: |_{C_{1, compl}^{ar}}$ is an equivalence of categories $C_{1, compl}^{ar}$ and
$Loc^{ar}$.
\item \label{st2}  For any $E \in \Ob(C_{1}^{ar})$ we have
$$
\Phi(\check{E}) = \widehat{\Phi(E) } \mbox{.}
$$
\item \label{st3} $\Phi$ maps the set of admissible triples from $C_{1}^{ar}$
onto the following set of exact triples of locally compact abelian
groups from $Loc^{ar}$:
$$
0 \lto H \lto G \lto G/H \lto 0 \mbox{,}
$$
where $H$ is a closed subgroup in $G$.
\end{enumerate}
\end{prop}
\proof. To prove statement~\ref{st1} of this proposition we
construct a functor
$$\Lambda \; : \; Loc^{ar} \lto C_{1, compl}^{ar}$$
in the following way. Let a group $G \in Loc^{ar}$.  We define
$\Lambda (G) = (I, F, G)$. Here $I$ is the set of  subgroups of $G$
such that $K \in I$ if and only if the subgroup $K$ coincides either
with a subgroup $U$, or with a subgroup $V$ satisfying properties
(*) -- (***) above. The set $I$ is partially ordered by inclusions
of subgroups. The function $F$ maps $K \in I$ to the corresponding
subgroup $K \subset G$. From reasonings before this proposition it
is clear that $(I, F, G) \in \Ob(C_{1, compl}^{ar})$.

Now let $\phi \in \Mor_{Loc^{ar}}(G_1, G_2)$. We show that
$\Lambda(\phi)=\phi$ is well-defined as an element from
$\Mor_{C_1^{ar}}(\Lambda(G_1), \Lambda(G_2))$, i.e., $\phi$
satisfies conditions~\ref{ii1}--\ref{ii3} of definition~\ref{dd1}.
Let $U_1 \supset V_1$ be subgroups of the group $G_1$ satisfying
properties (*) -- (***) above for the group $G_1$. Let $U_2 \supset
V_2$ be subgroups of the group $G_1$ satisfying properties (*) --
(***) above for the group $G_2$. Due to condition~(\ref{dimrank}),
without loss of generality we can assume  that $G/U_2$ is a discrete
torsion group. We have that $\phi(V_1)$ is a compact group, and thus
$\frac{\phi(V_1) +U_2}{U_2}$ is finite group. Besides,
$\frac{\phi(U_1) + U_2}{\phi(V_1) +U_2}$ is a finite group by
formula~(\ref{decomp}), therefore $\frac{\phi(U_1) +U_2}{U_2}$ is a
finite group. Hence the group $\phi(U_1) + U_2$ is from the set $I$.
Thus the map $\phi$ satisfies condition~\ref{ii1} of
definition~\ref{dd1}. To show that the map $\phi$ satisfies
condition~\ref{ii2} of definition~\ref{dd1} it is enough to verify
condition~\ref{ii1} of definition~\ref{dd1} for the dual map
$\widehat{\phi} \: : \: \widehat{G_2} \to \widehat{G_1}$. Indeed, if
$U \supset V$ are subgroups of a group $G \in \Ob(Loc)$ satisfying
properties (*) -- (***) above for the group $G$, then $V^{\bot}
\supset U^{\bot}$ are subgroups of the group $\widehat{G}$
satisfying properties (*) -- (***) above for the group
$\widehat{G}$. Here for any subgroup $H \subset G$ we defined the
subgroup $H^{\bot} \subset \widehat{G}$:
$$H^{\bot} \eqdef \{ g \in \widehat{G} : g \mid_H \equiv 1    \} \mbox{.}$$
But condition~\ref{ii1} of definition~\ref{dd1} we have just checked
for any map. Condition~\ref{ii3} of definition~\ref{dd1} is
satisfied for the map $\phi$, since any continuous homomorphism
between Lie groups is a smooth homomorphism between these groups.

Now it is easy to see that the functor $ \Phi \circ \Lambda$ is
isomorphic to the identity functor  $Id_{Loc^{ar}}  $, and the
functor $\Lambda \circ \Phi$ is isomorphic to the identity functor
$Id_{C_{1, compl}^{ar}}$. Therefore we proved statement~\ref{st1} of
proposition~\ref{propos9}.

Statement~\ref{st2} of proposition~\ref{propos9} follows from the
following observation. If $U \supset V$ are subgroups of a group $G
\in \Ob(Loc)$ satisfying properties (*) -- (***) above for the group
$G$, then $V^{\bot} \supset U^{\bot}$ are subgroups of the group
$\widehat{G}$ satisfying properties (*) -- (***) above for the group
$\widehat{G}$. Due to condition~(\ref{dimrank}), without loss of
generality we can assume that $V$ is a profinite group. Besides,
$$
\widehat{V} = \widehat{G}/ V^{\bot} \; \mbox{,} \qquad
\widehat{U/V}= V^{\bot}/U^{\bot} \; \mbox{,} \qquad \widehat{G/U}=
U^{\bot} \mbox{.}
$$

Statement~\ref{st3} of proposition~\ref{propos9} follows from the
following observation. Let
$$
0 \arrow{e} H \arrow{e} G \arrow{e,t}{\psi} G/H \arrow{e} 0
$$
be an exact sequence of locally compact abelian groups, where $H$ is
a closed subgroup in $G$. If $U \supset V$ are subgroups of the
group $G$ satisfying properties (*) -- (***) above for the group
$G$, then $(U \cap H) \supset (V \cap H) $ are subgroups of the
group $H$ satisfying properties (*) -- (***) above for the group
$H$, and $\psi(U) \supset \psi(V)$ are subgroups of the group $G/H$
satisfying properties (*) -- (***) above for the group $G/H$.

Proposition~\ref{propos9} is proved.

\begin{defin} Let $E = (I, F, V) \in \Ob(C_1^{ar})$.
\begin{enumerate}
\item We say that $E$ is a compact object iff there is an element $i_0 \in I$
such that $F(i_0)=V$, and  for any $i \le i_0 \in I$ we have that
$F(i_0)/F(i)$ is a compact Lie group.
\item We say that $E$ is a discrete object iff there is an element $i_0 \in I$
such that $F(i_0)=\{0\}$, and for any  $i \ge i_0 \in I$ we have
that
 $F(i)/F(i_0)$ is a discrete Lie group.
\end{enumerate}
\end{defin}
\begin{nt} {\em
Due to formula~(\ref{decomp}), $F(i_0)/F(i)$ is a compact Lie group
if and only if it is isomorphic to $\dt^k \times A$ for some
non-negative integer $k$ and a finite abelian group $A$.

Due to formula~(\ref{decomp}), $F(i)/F(i_0)$ is a discrete Lie group
if and only if it is isomorphic to $\dz^l \times B$ for some
non-negative integer $l$ and a finite abelian group $B$.}
\end{nt}

\begin{prop}
Let $E = (I, F, V) \in \Ob(C_1^{ar})$. Then
\begin{enumerate}
\item $E$ is a compact object if and only if $\Phi(E)$ is a compact
group.
\item $E$ is a discrete object if and only if $\Phi(E)$ is a
discrete group.
\end{enumerate}
\end{prop}
\proof. The group $\Phi(E)$ is a compact group if and only if it is
a projective limit of compact groups. The group $\Phi(E)$ is a
discrete group if and only if it is an inductive limit of discrete
groups. The proposition is proved.

\begin{example}  \label{exad}  {\em
We note that for $E_1= (I_1, F_1, V_1) \in \Ob(C_1^{ar})$ and $E_2=
(I_2, F_2, V_2) \in \Ob(C_1^{ar})$ we have the following canonical
construction:
\begin{equation} \label{for1}
E_1 \times E_2 \eqdef (I_1 \times I_2, F_1 \times F_2, V_1 \times
V_2) \in \Ob(C_1^{ar}) \mbox{,}
\end{equation}
where  $(i_1, i_2) \le (j_1, j_2) \in I_1 \times I_2$ if and only if
$i_1 \le j_1$ and  $i_2 \le j_2$, and for any $(k_1, k_2) \in I_1
\times I_2$ we define
$$(F_1 \times F_2) ((k_1, k_2)) =
F_1(k_1) \times F_2(k_2) \mbox{.}$$

Now we consider a number field $K$ such that $[K : \dq] =n$. Let $E
\subset K$ be the ring of integers in $K$. Then we have an
isomorphism of abelian groups $E \simeq \dz^n$.

Let $p_1, \ldots, p_l$ be all Archimedean places of the field $K$.
For any $1 \le i \le l$ we denote by $K_{p_i}$ the field which is
the completion of the field $K$ with respect to the absolute value
$p_i$. Then either $K_{p_i} \simeq \dc$, or $K_{p_i} \simeq \dr$. We
have the following isomorphism of Lie groups:
$$
\prod_{1 \le i \le l} K_{p_i} \simeq \dr^n \mbox{.}
$$

By example~\ref{examp1}, we consider
\begin{equation} \label{for2}
 (\{(0), (1)\}, F,  \prod_{1 \le i \le l} K_{p_i}) \in \Ob(C_1^{ar}) \mbox{,}
\end{equation}
 where $(0) < (1)$,
$F((0))=e$, $F((1))=\prod\limits_{1 \le i \le l} K_{p_i}$, $e$ is
the trivial subgroup of the group $\prod\limits_{1 \le i \le l}
K_{p_i}$.

From the diagonal embedding of the field $K$ we have the embedding
$E \subset \prod\limits_{1 \le i \le l} K_{p_i}$ such that the
following triple is an admissible triple from $C_0^{ar}$:
\begin{equation} \label{seqtor}
0 \lto E \lto  \prod\limits_{1 \le i \le l} K_{p_i} \lto  \dt^n \lto
0 \mbox{.}
\end{equation}

Let $\da_K^{fin}$ be the "finite" adele ring of the field $K$, i.e.,
it is a restricted product of completions of $K$ on equivalence
classes of all non-Archimedean absolute values with respect to
valuation rings. (In notations of example~\ref{exam1}, $\da_K^{fin}=
\da_{\Spec E}$.) Then we consider
\begin{equation} \label{for3}
(D, H, \da_K^{fin} ) \in \Ob(C_1^{ar}) \mbox{,}
\end{equation}
where $D$ is the set of divisors of $\Spec E$, for $d \in D$ we
define
$$
H(d) = d \cdot \prod_{p} \oo_{K_p} \subset \da_K^{fin} \mbox{,}
$$
where the  product is over all non-Archimedean absolute values, and
$\oo_{K_p}$ is the corresponding valuation ring. (It is clear that
$(D, H, \da_K^{fin} ) \in \Ob(C_1^{fin})$.)

Let $\da_K= \da_K^{fin} \times \prod\limits_{1 \le i \le l} K_{p_i}$
be the "full" adele ring of the field $K$. By
formulae~(\ref{for1})-(\ref{for3}), we have
$$
Q = (D_0 \cup D_1,  H \times F, \da_K ) \in \Ob(C_1^{ar}) \mbox{,}
$$
where the sets $D_0 = D \times (0)$, $D_1= D \times (1)$. It is
clear that $Q$ is a complete object from $C_1^{ar}$.

Diagonal embedding of $K$ induces the following exact triple of
abelian groups:
\begin{equation} \label{trip}
0 \arrow{e} K \arrow{e} \da_K \arrow{e,t}{\phi} \da_K/K \arrow{e} 0
\mbox{.}
\end{equation}
For any $i \in D_0$ we have that $(H \times F)(i) \cap K = 0$. By
the strong approximatiom theorem (see~\cite[ch. II, \S15]{A}), for
any $j \in D_1$ we have that $\phi((H \times F)(j))= \da/K$.

Let $P = (X, Y, K) \in \Ob(C_1^{ar})$, where the set $X$ is the set
of subgroups of $K$ given by $(H \times F)(i) \cap K$ for all $i \in
D_0 \cup D_1$. The set $X$ is partially ordered by inclusions of
subgroups, and the function $Y$ maps an element from $X$ to the
corresponding subgroup. Then the trivial subgroup $e \subset K$
belongs to $X$. By construction, we have that for any $j \in X$ the
group $Y(j)/e \simeq Y(j)\simeq \dz^n$. For any $j \ge i > e \in X$
we have that $Y(j)/Y(i)$ is a finite abelian group. Therefore $P$ is
a discrete object from $C_1^{ar}$. It is clear also that $P$ is a
complete object from $C_1^{ar}$.

Let $R = (Z, W, \da_K/K) \in \Ob(C_1^{ar})$, where the set $Z$ is
the set of subgroups of $\da_K/K$ given by $\phi((H \times F)(i)) $
for all $i \in D_0 \cup D_1$. The set $Z$ is partially ordered by
inclusions of subgroups, and the function $W$ maps an element from
$Z$ to the corresponding subgroup. Then the whole group $g=\da_K/ K$
belongs to $Z$. By construction, we have that for any $j \in Z$ the
group $g/W(j) \simeq \dt^n$. For any $g > j \ge i  \in Z$ we have
that $Y(j)/Y(i)$ is a finite abelian group. Therefore $R$ is a
compact object from $C_1^{ar}$. It is clear also that $R$ is a
complete object from $C_1^{ar}$.

By construction, we obtain the following admissible triple from
$C_1^{ar}$:
\begin{equation} \label{aaddrr}
0 \lto P \lto Q \lto R \lto 0 \mbox{.}
\end{equation}
This triple induces exact triple~(\ref{trip}) of abelian groups.

Besides, it is not hard to see from constructions that $\check{Q}
\simeq Q$, $\check{P} \simeq R$, $\check{R} \simeq P$.
 }
\end{example}

\begin{nt}{ \em
By formula~(\ref{for2}) we introduced the coarsest filtration of
$\prod\limits_{1 \le i \le l} K_{p_i} \simeq \dr^n $ when $l >1$. If
we introduce any other filtration of the space $\prod\limits_{1 \le
i \le l} K_{p_i}$ by $\dr$-vector subspaces $\prod\limits_{ i \in I}
K_{p_i}$ ($I$ is a subset of the finite set $\{1, \ldots, l\}$, for
$I =\emptyset$ we take the zero subspace) which dominates the
previous filtration and involves other $\dr$-vector subspaces than
zero subspace and $\dr^n$, then triple~(\ref{aaddrr}) will not be an
admissible triple from $C_1^{ar}$. Indeed, we have that $E \bigcap
(\prod\limits_{ i \in I} K_{p_i}) = 0$ for any $I \ne \{1, \ldots, l
\}$, where the intersection is taken inside of $\prod\limits_{1 \le
i \le l} K_{p_i}$. Hence, the image of $\prod\limits_{ i \in I}
K_{p_i}$ ($I \ne \emptyset$, $I \ne \{1, \ldots, l\}$) in $\dt^n$ is
not closed, see sequence~(\ref{seqtor}). (In such a way we obtain a
dense winding of a torus.)}
\end{nt}

\vspace{0.3cm}

In the sequel we will need the following technical lemma.
\begin{lemma} \label{lbc}
Let
$$
0 \arrow{e} E_1  \arrow{e,t}{\alpha}  E_2 \arrow{e,t}{\beta}   E_3
\arrow{e} 0
$$
be an admissible triple from $C_1^{ar}$.
\begin{enumerate}
\item Let $D \in \Ob(C_1^{ar})$, and $\gamma \in
\Mor_{C_1^{ar}} (D, E_3)$. Then there is the following admissible
triple from $C_1^{ar}$
\begin{equation} \label{adm1}
0 \arrow{e} E_1  \arrow{e,t}{\gamma_{\alpha}} E_2
\mathop{\times}_{E_3} D \arrow{e,t}{\gamma_{\beta}} D \arrow{e} 0
\end{equation}
and morphism $\beta_{\gamma} \in \Mor_{C_1^{ar}}(E_2
\mathop{\times}\limits_{E_3} D, E_2)$ such that the following
diagram is commutative:
\begin{equation} \label{adm2}
\begin{diagram}
\node{0} \arrow{e} \node{E_1}  \arrow{e,t}{\gamma_{\alpha}}
\arrow{s,=} \node{E_2 \mathop{\times}_{E_3} D}
\arrow{e,t}{\gamma_{\beta}}
\arrow{s,l}{\beta_{\gamma}} \node{D} \arrow{e} \arrow{s,r}{\gamma} \node{0} \\
\node{0} \arrow{e} \node{E_1}  \arrow{e,t}{\alpha}  \node{E_2}
\arrow{e,t}{\beta} \node{E_3} \arrow{e} \node{0}
\end{diagram}
\end{equation}
\item Let $A \in \Ob(C_1^{ar})$, and $\theta \in
\Mor_{C_1^{ar}}(E_1, A)$.
 Then there is the following admissible triple from
$C_1^{ar}$
\begin{equation}  \label{adm'1}
0 \arrow{e} A  \arrow{e,t}{\theta_{\alpha}} A  \mathop{\amalg}_{E_1}
E_2 \arrow{e,t}{\theta_{\beta}} E_3 \arrow{e} 0
\end{equation}
and morphism $\alpha_{\theta} \in \Mor_{C_1^{ar}}(E_2, A
\mathop{\amalg}\limits_{E_1} E_2)$ such that the following diagram
is commutative:
\begin{equation}  \label{adm'2}
\begin{diagram}
\node{0} \arrow{e} \node{E_1}  \arrow{e,t}{\alpha}
\arrow{s,l}{\theta} \node{E_2} \arrow{e,t}{\beta}
\arrow{s,r}{\alpha_{\theta}} \node{E_3} \arrow{e} \arrow{s,=}  \node{0} \\
\node{0} \arrow{e} \node{A}  \arrow{e,t}{\theta_{\alpha}} \node{A
\mathop{\amalg}\limits_{E_1} E_2} \arrow{e,t}{\theta_{\beta}}
\node{E_3} \arrow{e} \node{0}
\end{diagram}
\end{equation}
\end{enumerate}
\end{lemma}
\proof. We prove the first statement of the lemma. Let $E_i = (I_i,
F_i, V_i)$, $ 1 \le i \le 3$, and
 $D = (K, H, T)$.

 We construct $E_2 \mathop{\times}\limits_{E_3}  D = (J, G, W) \in \Ob(C_1^{ar})$ in the
 following way. We define an abelian group
 $$
 W = V_2 \mathop{\times}_{V_3} T \eqdef \{ (e,d) \in V_2 \times T \quad \mbox{such that} \quad  \beta(e) = \gamma(d)
 \} \mbox{.}
 $$
We define a partially ordered set
$$
J = \{ (i,j) \in I_2 \times K  \quad \mbox{such that} \quad
\gamma(H(j)) \subset \beta(F_2(i)) \}  \mbox{,}
$$
where $(i_1, j_1) \le (i_2, j_2)$ iff $i_1 \le i_2$ and $j_1 \le
j_2$. We define a function $G$ from $J$ to the set of  subgroups of
$W$ as
$$
G((i,j)) = (F_2(i) \times H(j) ) \cap W  \mbox{,}
$$
where the intersection is inside $V_2 \times T$.

Let $(i_1, j_1), (i_2, j_2) \in J$, and $(i_1, j_1) \le (i_2, j_2)$.
Then we have the following commutative diagram of morphisms between
abelian Lie groups:
$$
\begin{diagram}
\dgARROWLENGTH=1.3em
 \node{0} \arrow{e} \node{\frac{F_2(i_2) \cap V_1} {F_2(i_1) \cap
V_1} } \arrow{e}  \arrow{s,=} \node{\frac{G((i_2, j_2))}
{G((i_1,j_1))}} \arrow{e} \arrow{s} \node{\frac{H(j_2)}{H(j_1)}}
\arrow{e} \arrow{s} \node{0}
\\ \node{0} \arrow{e} \node{\frac{F_2(i_2) \cap V_1}{F_2(i_1) \cap
V_1}} \arrow{e}  \node{\frac{F_2(i_2)}{F_2(i_1)}} \arrow{e}
\node{\frac{F_2(i_2)/ (F_2(i_2) \cap V_1)}{F_2(i_1)/ (F_2(i_1)\cap
V_1)}} \arrow{e} \node{0}
\end{diagram}
$$
where
$$
\frac{G((i_2, j_2))} {G((i_1,j_1))} = \frac{F_2(i_2)}{F_2(i_1)}
 \mathop{\times}_{\frac{F_2(i_2)/ (F_2(i_2) \cap
V_1)}{F_2(i_1)/ (F_2(i_1)\cap V_1)}} \frac{H(j_2)}{H(j_1)}  \mbox{,}
$$
and horizontal triples are admissible triples from $C_0^{ar}$.

 Hence, $E_2 \mathop{\times}\limits_{E_3}
D$ is well-defined as an object from $C_1^{ar}$, and the maps
$\beta_{\gamma}$ and $\gamma_{\beta}$ are projections. Besides, from
this commutative diagram we obtain that triple~{(\ref{adm1})} is an
admissible triple from $C_1^{ar}$, and diagram~{(\ref{adm2})} is a
commutative diagram. The first statement of the lemma is proved.

We prove the second statement of the lemma, which is the dual
statement to the first statement of the lemma. Let $A= (K', H',
T')$.

We construct $A  \mathop{\amalg}\limits_{E_1} E_2 = (J', G', W') \in
\Ob(C_1)$ in the following way. We remark that $T' \amalg V_2 = T'
\times V_2$.
 We define an abelian group
$$
W' = T' \mathop{\amalg}_{V_1}   V_2 \eqdef (T' \times V_2) / E \quad
\mbox{, where} \quad E = (\theta \times \alpha)  V_1  \mbox{.}
$$
We define a partially ordered set
$$
J' = \{(i,j) \in K' \times I_2   \; : \; \theta(F_2(j) \cap V_1)
\subset H'(i) \} \mbox{,}
$$
where $(i_1, j_1) \le (i_2, j_2)$ iff $i_1 \le i_2$ and $j_1 \le
j_2$. We define a function $G'$ from $J'$ to the set of  subgroups
of $W'$ as
$$
G'((i,j)) = ((H'(i) \times F_2(j)) + E / E )      \quad \subset
\quad W'  \mbox{.}
$$

Let $(i_1, j_1), (i_2, j_2)  \in J'$, and $(i_1, j_1) \le (i_2,
j_2)$. Then we have the following commutative diagram of morphisms
between abelian Lie groups:
$$
\begin{diagram}
\dgARROWLENGTH=0.9em
 \node{0}  \arrow{e} \node{\frac{F_2(j_2) \cap V_1} {F_2(j_1) \cap
V_1} } \arrow{e}  \arrow{s} \node{\frac{F_2(j_2)}{F_2(j_1)}}
\arrow{e} \arrow{s} \node{\frac{F_2(j_2)/ (F_2(j_2) \cap
V_1)}{F_2(j_1)/ (F_2(j_1)\cap V_1)}} \arrow{e} \arrow{s,=} \node{0}
\\ \node{0}  \arrow{e} \node{\frac{H'(i_2)}{H'(i_1)}} \arrow{e}
\node{\frac{G'((i_2,j_2))}{G'((i_1,j_1))}} \arrow{e}
\node{\frac{F_2(j_2)/ (F_2(j_2) \cap V_1)}{F_2(j_1)/ (F_2(j_1)\cap
V_1)}}
 \arrow{e} \node{0}
\end{diagram}
$$
where
$$
\frac{G'((i_2, j_2))} {G'((i_1,j_1))} = \frac{H'(i_2)}{H'(i_1)}
 \mathop{\amalg}_{ \frac{F_2(j_2) \cap V_1} {F_2(j_1) \cap
V_1} } \frac{F_2(j_2)}{F_2(j_1)}  \mbox{,}
$$
and horizontal triples are admissible triples from $C_0^{ar}$.

Hence, $A  \mathop{\amalg}\limits_{E_1} E_2$ is well defined as an
object from category $C_1^{ar}$. Besides,  from this commutative
diagram we obtain that triple~{(\ref{adm'1})} is an admissible
triple from $C_1^{ar}$, and diagram~{(\ref{adm'2})} is a commutative
diagram. The second statement of the lemma is proved. The lemma is
proved.

\section{Functions and distributions on objects of $C_0^{ar}$ and $C_1^{ar}$}
\subsection{Functions and distributions on objects of $C_0^{ar}$}
\subsubsection{Definitions of basic spaces}

For any $G \in \Ob(C_0^{ar})$ we define the Schwartz space $\s(G)$
of rapidly decreasing functions on $G$ in the following way
(see~\cite[p.~138]{Br} and~\cite[ch.~VII]{S}).

Let $G \in \Ob(C_0^{ar})$ such that
\begin{equation} \label{poldec}
G \simeq \dz^r \times \dr^q \mbox{.}\end{equation} We say that a
function $p : G \to \dc $ is a polynomial if and only if $p$ is
polynomial of coordinate functions with respect to
decomposition~(\ref{poldec}). It is not difficult to see that this
definition does not depend on the choice of
isomorphism~(\ref{poldec}). Now let $G \in \Ob(C_0^{ar})$ be any
object. Then we have canonical filtration~(\ref{dee}) on $G$. The
subgroup $K$ is an element of this filtration. We denote by
$(G/K)_t$ the subgroup in the group $G/K$ which consists of torsion
elements of this group. Then the group $(G/K)/ (G/K)_t \simeq \dz^r
\times \dr^q$ (see decomposition~(\ref{decomp})). Let $u : G \to
(G/K)/ (G/K)_t$ be the natural map. Then a set $P$ of polynomials on
$G$ is defined as
$$
P \eqdef \{p : p =u^*(p') \} \mbox{,}
$$
where $p'$ runs over the set of all polynomials on $(G/K)/ (G/K)_t$,
which was defined above.

By $D$ we denote the set of all invariant differential operators on
an object $G \in \Ob(C_0^{ar})$, which is a commutative Lie group.

Let $G \in \Ob(C_0^{ar})$. We define
$$
\s(G)= \{\mbox{all} \quad \dc-\mbox{valued smooth functions}\quad f
\quad \mbox{on} \quad G \quad \mbox{such that} $$
$$ \sup_{x \in G}
\mid p(x) df(x) \mid < \infty \quad \mbox{for any} \quad p \in P
\mbox{,} \quad d \in D  \} \mbox{.}
$$

The topology on $\s(G)$ is defined by the following system of
seminorms $\{s_{p,d}\}_{p \in P, d \in D}$:
$$
s_{p,d}(f) \eqdef \sup_{x \in G} \mid p(x) df(x)  \mid \mbox{,}
\qquad \mbox{where} \quad p\in P \mbox{,} \quad d \in D \mbox{,}
\quad f \in \s(G) \mbox{.}
$$
It means that the base of this topology is the set of open balls of
these seminorms
$${\mathop{B}^{\circ}}_{s_{p,d}}(x,r)= \{y \in G : s_{p,d}(x-y) < r \} \mbox{,} \quad
x \in G \mbox{,} \quad r > 0 \mbox{,} \quad p \in P \mbox{,} \quad d
\in D$$
 and their finite intersections. Therefore $\s(G)$ is a locally
 convex topological vector space.

 If we take only monomials of coordinate functions (after a fixed decomposition~(\ref{decomp}))
 instead of all polynomials
 in definition of the set of seminorms above,
 and if we take only monomials of
 partial derivations with respect to coordinates (after the fixed decomposition~(\ref{decomp}))
 instead of all
 invariant differential operators,
 then we obtain the new system of seminorms $\{s_i\}_{i \in \sdn}$ which will be countable
 and equivalent to the previous system of seminorms $\{s_{p,d}\}_{p \in P, d \in D}$.
 (It means  that this system of
seminorms generates the same topology on $\s(G)$ as the  system of
seminorms $\{s_{p,d}\}_{p \in P, d \in D}$ generates.)
   Therefore the space $\s(G)$  is a metrizable
 vector space, where a translation invariant metric $d(\,,\,)$ can be
 defined in the following way:
$$
d(x,y)= \sum_{i=1}^{\infty} 2^{-i} \frac{s_i(x-y)}{1+s_i(x-y)}
\mbox{,} \qquad x,y \in \s(G) \mbox{.}
$$
(Here, without loss of generality, by change to an equivalent system
of seminorms, we assumed that $s_1 \le s_2 \le \ldots \le s_m \le
\ldots$.)

It can be also proved that $\s(G)$ is a complete vector space.
Therefore $\s(G)$ is a Fr\'echet space.

For a locally convex  topological $\dc$-vector space $V$ we can
consider a lot of topologies on its continuous dual $\dc$-vector
space $V'$ (see~\cite{RR}). We will always consider the weak-${}^*$
topology on $V'$. This topology is generated by the following system
of seminorms: $\{ s_v\}_{v \in V}$, where for $v \in V$ the seminorm
$s_v$ is defined as $s_v(f)= \mid f(v )\mid$ for any $f \in V'$. The
topology so defined on $V'$ is a locally convex topology. Besides,
we always have that the natural map $V \to (V')'$ is an isomorphism
of $\dc$-vector spaces. (But they are not isomorphic as topological
vector spaces.)

For any $G \in \Ob(C_0^{ar})$ we define the Schwartz space $\s'(G)$
of tempered distributions on $G$ as $\s'(G) \eqdef \s(G)'$. Since
the topological vector space $\s(G)$ is a complete and metrizable
space, the $\dc$-vector topological space $\s'(G)$ is a sequentially
complete space.

Now we consider some examples of Schwartz spaces for various types
of objects  $ G \in \Ob(C_0^{ar})$
\begin{example} \label{ex5} {\em
Let $G \simeq \dr^n$. Then the space $\s(\dr^n)$ consists of all
smooth $\dc$-valued functions $f$ on $\dr^n$ such that
$$
s_{\alpha, \beta}(f) \eqdef \sup_{x \in \sdr^n} \mid x^{\alpha}
\partial^{\beta}f(x) \mid < \infty \mbox{,}
$$
where multi-indices $\alpha= (\alpha_1, \ldots, \alpha_n)$,
$\beta=(\beta_1, \ldots, \beta_n)$ consist of non-negative integers,
$x^{\alpha}= x_1^{\alpha_1} \ldots x_n^{\alpha_n}$,
$\partial^{\beta}= \frac{\partial^{\beta_1}}{\partial x_1^{\beta_1}}
\ldots \frac{\partial^{\beta_n}}{\partial x_n^{\beta_n}}$. The
system of norms $s_{\alpha, \beta}$ (for all multi-indices $\alpha$
and $\beta$) on the space $\s(\dr^n)$ is equivalent to the system of
seminorms $\{s_{p,d}\}_{p\in P, d\in D}$, which was defined above
for any $G \in \Ob(C_0^{ar})$. We note also that a system of norms
$$
s'_{\alpha, \beta}(f) \eqdef \int_{\sdr^n} \mid x^{\alpha}
\partial^{\beta} f (x) \mid dx
$$
(for all multi-indices $\alpha$  and $\beta$) is also equivalent to
the previous systems of seminorms on $\s(\dr^n)$.}
\end{example}

\begin{example}{\em
Let $G \simeq \dt$. Then the space $\s(\dt)$ consists of all smooth
$\dc$-valued functions  on $\dt$. The topology on $\s(\dt)$ is given
by the following system of seminorms $\{s_{\alpha}\}_{\alpha \in
\sdn}$:
$$
s_{\alpha}(f) \eqdef \sup_{t \in \sdt} \mid
\partial^{\alpha}f(t) \mid < \infty \mbox{,}
$$
where $\alpha \in \dn$,
$\partial^{\alpha}=\frac{\partial^{\alpha}}{\partial t^{\alpha}}$,
$f \in \s(\dt)$. The system of seminorms $\{s_{\alpha}\}_{\alpha \in
\sdn}$ on the space $\s(\dt)$ is equivalent to the system of
seminorms $\{s_{p,d}\}_{p\in P, d\in D}$,  which was defined above
for any $G \in \Ob(C_0^{ar})$. We note also that a system of
seminorms
$$
s'_{\alpha}(f) \eqdef \int_{\sdt} \mid
\partial^{\alpha} f (t) \mid dt
$$
(where $\alpha \in \dn$) is also equivalent to the previous systems
of seminorms on $\s( \dt )$.

The same reasonings can be  done analogously for $G \simeq \dt^n$,
$n
>1$. }
\end{example}

\begin{example}{\em
Let $G \simeq \dz$. Then $\s(\dz)$ is the space of two-sided
sequences with some condition:
$$
\s(\dz)=\{ a_n \; : \; n \in\dz, \; a_n \in \dc \quad \mbox{such
that} \quad \mid a_n\mid = O(\mid n^{-k}\mid) \quad \mbox{for any}
\quad k \in \dn \}
 \mbox{.}
$$
The topology in $\s(\dz)$ is given by the following system of norms
$\{s_k\}_{k \in \sdn}$:
$$
s_k(\{a_n\})= \sup_{n \in \sdz} \mid  n^k a_n\mid \mbox{.}
$$
(This system of norms is equivalent to the system of seminorms
$\{s_{p,d}\}_{p\in P, d\in D}$ on $\s(G)$ which was introduced
above.) It is clear that $\s(\dz) \subset l^2(\dz)$.

The space $\s'(\dz)$ is defined as the space of the following
sequences:
$$
\s'(\dz)=\{ b_n \; : \; n \in\dz, \; b_n \in \dc \quad \mbox{such
that} \quad \mid b_n\mid = O(\mid n^{k}\mid) \quad \mbox{for some}
\quad k \in \dn \}
 \mbox{.}
$$
It is clear that $\s'(\dz) \supset l^2(\dz)$. Now if $a = \{ a_n\}
\in \s(\dz)$, $b = \{ b_n  \} \in \s'(\dz)$, then $b(a) \eqdef
\sum_{n \in \sdz} b_n a_n $ is an absolutely convergent series.
Moreover, by definition, linear functional $b$ is continuous with
respect to some norm $s_k$ on $\s(\dz)$. Then we have $\mid b_n\mid
= O(\mid n^{k}\mid)$. Conversely, if $\mid b_n\mid = O(\mid
n^{k}\mid)$ for some $k \in \dn$, then $b=\{ b_n\}$ is a continuous
linear functional with respect to the norm $s_{k+2}$ on $\s(\dz)$.

The same reasonings can be  done analogously for $G \simeq \dz^n$,
$n>1$.
}
\end{example}

\subsubsection{Fourier transform} \label{ft00}
Let $G \in \Ob(C_0^{ar})$. For any $a \in G$, for any $ f \in \s(G)$ we define $T_a(f) \in \s(G)$ as
$$T_a(f)(b) \eqdef f(b+a) \mbox{,} \qquad b \in G \mbox{.}$$
For any $a \in G$, for any $ H \in \s'(G)$ we define $T_a(H) \in \s'(G)$ as
$$T_a(H)(f) \eqdef H(T_{-a}(f)) \mbox{,} \qquad  f \in G \mbox{.}$$

We define the following $\dc$-vector space:
$$
\s'(G)^G \eqdef \{H \in \s'(G) \; : \; T_a(H)=H \quad \mbox{for any} \quad a \in G \mbox{.}    \}
$$

Any $G \in \Ob(C_0^{ar})$ is a locally compact Abelian
group. Therefore there exists a Haar measure $\nu$  on $G$, which is defined up to a
multiplication on positive real number. We define the following element ${\bf 1}_{\nu} \in \s'(G)$ as
$$
\de_{\nu} (f) \eqdef \int_{G} f(x) d \nu(x) \mbox{,} \quad f \in \s(G) \mbox{.}
$$
Besides, we have a well-defined map $\I_{\nu} : \s(G) \to \s'(G) $ by the rule:
$$
\I_{\nu}(f)(g)  \eqdef \de_{\nu} (f g)  \mbox{,} \quad f \mbox{,} \; g \in \s(G) \mbox{.}
$$

\begin{prop} \label{hm}
Let $G \in \Ob(C_0^{ar})$, and   $\nu$ be a Haar measure on $G$, then  the following properties are satisfied.
\begin{enumerate}
\item  $ \dim_{\sdc} \s'(G)^G = 1$.
\item  $\s'(G)^G = \dc \cdot \de_{\nu}$ inside the space $\s'(G)$.
\end{enumerate}
\end{prop}
\proof.
It is enough to suppose that $G$ is a connected Lie group. By definition of a Haar measure, we have $\de_{\nu} \in \s'(G)^G$.

Now let a vector $v \in T_{G,e}=Lie G$, $v \ne 0$,  then the vector $v$
defines an invariant  vector field $X_v$ on $G$, i.e, a first order differential operator $X_v : \s(G) \to \s(G)$.
Since the operator $X_v$ is continuous, it can be extended to the operator $X_v : \s'(G) \to \s'(G)$ by the rule:
$$
X_v (H)(f) \eqdef -H (X_v(f)) \mbox{.}
$$

On the other hand, the vector $v$ defines the one-parametric subgroup
$$g_v(t) \eqdef exp(tv) \subset G \mbox{, } \quad t \in \dr \mbox{.}$$
By virtue of mean value theorem and Taylor's formula we have  for any $f \in \s(G)$
$$
\lim_{t \to 0} T_{g_v(t)}(f) =f \quad \mbox{, and}
$$
$$
\lim_{t \to 0} \frac{T_{g_v(t)}(f) -f}{t} = X_v (f)
$$
in the space $\s(G)$
for any seminorm $s_{p,d}$, $p \in P$, $d \in D$ on this space.

Hence, if $H \in \s'(G)^G$, then $X(H)=0$ for any invariant vector field $X$ on $G$. Then, by induction on $\dim G$ and using the standard arguments
of  theory of generalized functions, we obtain that $H = c \cdot \de_{\nu}$ for some $c \in \dc$. (See~\cite[ch. 1, \S 2.6]{GSh} for the case $G = \dr$.) The proposition is proved.

\vspace{0.5cm}

For any $G \in \Ob(C^{ar}_0)$ we define the one-dimensional
$\dc$-vector space $\mu(G)$ as
$$
\mu(G) \eqdef \s'(G)^G \mbox{.}
$$

We note that by the properties of Haar measures and by proposition~\ref{hm} we have that if
$$
0 \lto G_1 \lto G_2 \lto G_3 \lto 0
$$
is an admissible triple from $C^{ar}_0$, then there is a natural isomorphism
\begin{equation} \label{isotr}
\mu(G_1) \otimes_{\sdc} \mu(G_3) \simeq \mu(G_2) \mbox{.}
\end{equation}

Let $G \in \Ob(C^{ar}_0)$, and $\hat{G} \in \Ob(C^{ar}_0)$ be its Pontryagin dual group. Let $\mu \in \mu(G)$.
Then the Fourier transform
$\F_{\mu} : \s(G) \to \s(\hat{G})$ is defined as
\begin{equation} \label{defftr}
\F_{\mu}(\chi) \eqdef \mu (f \overline{\chi}) \mbox{,}  \quad \chi \in \hat{G} \mbox{.}
\end{equation}
The map $\F_{\mu}$ depends linearly on $\mu \in \mu(G)$. Therefore
the following  map $\F : \s(G) \otimes_{\sdc} \mu(G) \to
\s(\hat{G})$ is well-defined as
$$
\F (f \otimes \mu) \eqdef \F_{\mu}(f)  \mbox{.}
$$

As the conjugate map to the continuous map $\F$ we define the Fourier transform on distributions:
$$
\F : \s'(G)\otimes_{\sdc} \mu(\hat{G})  \lto \s'(\hat{G})
$$
such that $\F_{\nu} (H) \eqdef \F(H \otimes \nu)$, $\nu \in \mu(\hat{G})$, $H \in \s'(G)$.

From properties of the Fourier transform it follows that for any $\mu \in \mu(G)$, any $\nu \in \mu(\hat{G})$ we have
$$
\F_{\nu}(\mu)= c \delta
$$
for some $c \in \dc$ (which linearly depends on $\mu$ and $\nu$),
and $\delta \in \s'(\hat{G})$ is the Dirac delta function: $\delta
(f) \eqdef f(0)$ for $f \in \s(\hat{G})$. Hence we obtain a
well-defined isomorphism:
\begin{equation} \label{isodu}
\mu(G) \otimes_{\sdc} \mu(\hat{G}) \simeq \dc  \mbox{.}
\end{equation}

Let an element $\mu \in \mu(G)$, $\mu \ne 0$. Then  the element
$\mu^{-1} \in \mu(\hat{G})$ is well-defined such that $\mu \otimes
\mu^{-1} = 1$ with respect to isomorphism~(\ref{isodu}).

Let $f \in \s(G)$. Then we define $\check{f} \in \s(G)$ as
\begin{equation} \label{ch1}
\check{f}(x) \eqdef f(-x)\mbox{,} \quad x \in G \mbox{.}
\end{equation}

Let $H \in \s'(G)$. Then we define $\check{H} \in \s'(G)$ as
\begin{equation} \label{ch2}
\check{H}(f) \eqdef H(\check{f}) \mbox{,} \quad f \in \s(G) \mbox{.}
\end{equation}

We collect the other well-known properties of the Fourier transform in the following proposition.
\begin{prop} \label{fourtra}
Let $G \in \Ob(C_0^{ar})$, $\mu \in \mu(G)$, $\mu \ne 0$. The following properties are satisfied.
\begin{enumerate}
\item The map $\F_{\mu}$ is an isomorphism of topological $\dc$-vector spaces $\s(G)$ and $\s(\hat{G})$.
\item The map $\F_{\mu^{-1}}$ is an isomorphism of topological $\dc$-vector spaces $\s'(G)$ and $\s'(\hat{G})$.
\item  \label{it3ft} $ \F_{\mu^{-1}} \F_{\mu}(f) = \check{f}$ for any $f \in \s(G)$.
\item   $ \F_{\mu}   \F_{\mu^{-1}} (H) =\check{H}$ for any $H \in \s'(G)$.
\end{enumerate}
\end{prop}

\subsubsection{Direct and inverse images.} \label{diim}
We consider an admissible triple from $C_0^{ar}$:
\begin{equation}  \label{admtr}
0 \arrow{r} G_1 \arrow{r,t}{\alpha} G_2 \arrow{r,t}{\beta} G_3 \arrow{r} 0 \mbox{.}
\end{equation}

There is the following map $\beta_* : \s(G_2) \otimes_{\sdc} \mu(G_1)  \to \s(G_3)$ which is defined as
$$
\beta_*(f \otimes \mu)(x) \eqdef \mu(f_{y}) \mbox{,}
$$
where $f \in \s(G_2)$, $\mu \in \mu(G_1)$, $x \in G_3$, $y \in f^{-1}(x)$, $f_{y} \in \s(G_1)$ is defined as $f_y(z)  \eqdef f(z+y) $
for $z \in G_1$. The function $\beta_*(f \otimes \mu)$ does not depend on the choice of $y \in f^{-1}(x)$, because  $\mu \in \s'(G_1)$
is an invariant element.

It is not difficult to see that the map $\beta_*$ is continuous, therefore there is the continuous map $\beta^*: \s'(G_3) \otimes_{\sdc} \mu(G_1) \to \s'(G_2)$ which is defined as the conjugate map to the map $\beta_*$.

There is the map $\alpha^* : \s(G_2) \to \s(G_1)$ defined as
$$
\alpha^*(f)(x) \eqdef f(\alpha(x)) \mbox{,} \quad \mbox{where} \quad f \in \s(G_2) \mbox{,} \quad x \in G_1 \mbox{.}
$$

The map $\alpha^*$ is continuous, therefore there is the continuous map $\alpha_*: \s'(G_1)  \to \s'(G_2)$  defined as the conjugate map to the map $\alpha^*$.

If $G_1$ is a compact Lie group, and $G_3$ is any object in sequence~(\ref{admtr}), then there is the following map
$\beta^* : \s(G_3) \to \s(G_2) $ defined as
$$
\beta^*(f)(x) \eqdef f(\beta(x)) \mbox{,} \quad \mbox{where} \quad f \in \s(G_3) \mbox{,} \quad x \in G_2 \mbox{.}
$$

The map $\beta^*$ is continuous, therefore there is the continuous map $\beta_*: \s'(G_2)  \to \s'(G_3)$  defined as the conjugate map to the map $\beta^*$.

If $G_3$ is a discrete Lie group, and $G_1$ is any object in sequence~(\ref{admtr}), then there is the following map
$\alpha_* : \s(G_1) \to \s(G_2) $ defined as
$$
\alpha_*(f)(x) \eqdef
\left\{
\begin{array}{ccc}
0  & \mbox{if} & x \notin f(G_1) \mbox{;}\\
f(y) & \mbox{if} & x = f(y) \mbox{,}
\end{array}
\right.
\quad \mbox{where} \quad  f  \in \s(G_1) \mbox{,} \quad  x \in G_2 \mbox{.}
$$

The map $\alpha_*$ is continuous, therefore there is the continuous map $\alpha^*: \s'(G_2)  \to \s'(G_1)$  defined as the conjugate map to the map $\alpha_*$.

\vspace{0.5cm}

Let
$$
0 \arrow{e} E_1  \arrow{e,t}{\alpha}  E_2 \arrow{e,t}{\beta}   E_3
\arrow{e} 0
$$
be an admissible triple from $C_0^{ar}$.
Let
$$
0 \arrow{e} D \arrow{e,t}{\gamma} E_3 \arrow{e,t}{\delta} B
 \arrow{e} 0 \mbox{.}
$$
be another admissible triple from $C_0^{ar}$.

There is the following commutative diagram:
$$
\begin{diagram}
\node{E_2 \mathop{\times}\limits_{E_3} D}
\arrow{s,l}{\beta_{\gamma}} \arrow{e,t}{\gamma_{\beta}} \node{D} \arrow{s,r}{\gamma} \\
\node{E_2} \arrow{e,b}{\beta} \node{E_3}
\end{diagram}
$$
where $\gamma_{\beta}$ is an admissible epimorphism (i.e., it is an epimorphism part of some admissible triple from $C_0^{ar}$), and
$\beta_{\gamma}$ is an admissible monomorphism (i.e., it is a monomorphism part of some admissible triple from $C_0^{ar}$).

We note that  $E_3 \simeq E_2  \mathop{\amalg}\limits_{G} D$, where $G = E_2 \mathop{\times}\limits_{E_3} D$.

There is the following proposition.
\begin{prop} \label{pr18}
We have the following formulas.
\begin{enumerate}
\item For any $f \in \s(E_2)$, $H \in \s'(D)$, $\mu \in
\mu(E_1)$
\begin{equation}  \label{e''1}
\gamma^* \beta_{*} (f \otimes \mu) = (\gamma_{\beta})_*
(\beta_{\gamma}^* (f) \otimes \mu)  \mbox{.}
\end{equation}
\begin{equation} \label{e''2}
\beta^* (\gamma_*(H) \otimes \mu)=(\beta_{\gamma})_*
\gamma_{\beta}^* (H \otimes \mu)  \mbox{.}
\end{equation}
\item If $E_1$ is a compact Lie group, then  for any $f \in \s(E_3)$,  $H \in \s'(G)$
\begin{equation} \label{e''3}
\gamma_{\beta}^* \gamma^* (f) = \beta_{\gamma}^* \beta^* (f)
\mbox{.}
\end{equation}
\begin{equation} \label{e''4}
\beta_* (\beta_{\gamma})_* (H) = \gamma_* (\gamma_{\beta})_* (H)
\mbox{.}
\end{equation}
\item If $B$ is a discrete Lie group, then for any $f \in \s(G)$, $H \in \s'(E_3)$, $\mu \in \mu(E_1)$
\begin{equation} \label{pr18f4c}
\beta_* ((\beta_{\gamma})_* (f) \otimes \mu) = \gamma_*
(\gamma_{\beta})_* (f \otimes \mu) \mbox{.}
\end{equation}
\begin{equation} \label{pr18f4d}
\gamma_{\beta}^* (\gamma^* (H) \otimes \mu ) = \beta_{\gamma}^*
\beta^* (H \otimes \mu) \mbox{.}
\end{equation}
\item If $E_1$ is a compact
Lie group,  and $B$ is a discrete Lie group, then  for any  $f \in \s(D)$, for any $H \in \s'(E_2)$
\begin{equation} \label{pr18f5}
\beta^* \gamma_*(f) =(\beta_{\gamma})_* \gamma_{\beta}^* (f)
\mbox{.}
\end{equation}
\begin{equation}  \label{pr18f6}
\gamma^* \beta_{*} (H) = (\gamma_{\beta})_* \beta_{\gamma}^* (H)
\mbox{.}
\end{equation}
\end{enumerate}
\end{prop}
\proof for functions follows from corresponding definitions. Proof for distributions follows from conjugate formulas to formulas for functions. The proposition is proved.

\vspace{0.5cm}

Let
$$
0 \arrow{e} E_1  \arrow{e,t}{\alpha}  E_2 \arrow{e,t}{\beta}   E_3
\arrow{e} 0
$$
be an admissible triple
from $C_0^{ar}$.
Let
$$
0 \arrow{e} L  \arrow{e,t}{\alpha'}  H \arrow{e,t}{\beta'}   E_2
\arrow{e} 0
$$
be another admissible triple from $C_0^{ar}$.
Then we have  the following admissible triple  from $C_0^{ar}$:
$$
0 \arrow{e} H \mathop{\times}_{E_2} E_1 \arrow{e,t}{\beta'_{\alpha}}
H \arrow{e,t}{\beta \beta'} E_3 \arrow{e} 0  \mbox{.}
$$
We denote $E = H \mathop{\times}\limits_{E_2} E_1 \in \Ob(C_0^{ar})$.

\begin{prop}. \label{pr16} We have the following formulas.
\begin{enumerate}
\item For any $f \in \s(H)$, $G \in \s'(E_3)$, $\nu \in \mu(L)$, $\mu \in
\mu(E_1)$
\begin{equation}
\label{e1} (\beta \beta')_* (f \otimes (\nu \otimes \mu)) =  \beta_*
(\beta'_* (f \otimes \nu) \otimes \mu) \mbox{.}
\end{equation}
\begin{equation}
\label{e2} (\beta \beta')^* (G \otimes (\nu \otimes \mu)) =
(\beta')^*  (\beta^*  (G \otimes \mu) \otimes \nu)  \mbox{.}
\end{equation}
\item If $E_1$ and $L$ are compact Lie groups, then $E$ is a compact
 Lie group, and the following formulas are satisfied for any $f \in \s(E_3)$, $G \in \s'(H)$
\begin{equation} \label{pr16f5}
(\beta \beta')^*(f) = (\beta')^* \beta^* (f) \mbox{.}
\end{equation}
\begin{equation} \label{pr16f6}
(\beta \beta')_*(G)= \beta_* (\beta')_* (G) \mbox{.}
\end{equation}
\end{enumerate}
\end{prop}
\proof.
We have the following admissible triple from $C_0^{ar}$:
\begin{equation} \label{rastr}
0 \arrow{e} L \arrow{e,t}{\alpha_{\alpha'}}
 E  \arrow{e,t}{\alpha_{\beta'}} E_1 \arrow{e} 0  \mbox{.}
\end{equation}
Therefore, by formula~(\ref{isotr}), for  $\nu \in \mu(L)$, $\mu \in \mu(E_1)$ we have
canonically that $\nu \otimes \mu \in \mu(E))$.
Now formula~(\ref{e1}) follows from Fubini's theorem. Formula~(\ref{e2}) is the conjugate
formula to formula~(\ref{e1}).

From triple~(\ref{rastr}) we obtain that if  $E_1$ and $L$ are
compact  Lie groups, then $E$ is
also a compact  Lie group.
Now formula~(\ref{pr16f5}) follows from definitions.
Formula~(\ref{pr16f6}) is the conjugate formula to
formula~(\ref{pr16f5}). The proposition is proved.

\vspace{0.5cm}

Let
$$
0 \arrow{e} E_1  \arrow{e,t}{\alpha}  E_2 \arrow{e,t}{\beta}   E_3
\arrow{e} 0
$$
be an admissible triple from $C_0^{ar}$.
Let
$$
0 \arrow{e} E_2  \arrow{e,t}{\alpha'}  H' \arrow{e,t}{\beta'}   L'
\arrow{e} 0  \mbox{.}
$$
be another admissible triple from $C_0^{ar}$.
Then
we have the following admissible triple from $C_0^{ar}$:
$$
0 \arrow{e} E_1 \arrow{e,t}{\alpha' \alpha} H'
\arrow{e,t}{\alpha'_{\beta}} E_3 \mathop{\amalg}\limits_{E_2} H'
\arrow{e} 0  \mbox{.}
$$

\begin{prop}. \label{pr17}
We have the following formulas.
\begin{enumerate}
\item  For any $f \in \s(H')$, $G \in \s'(E_1)$
\begin{equation} \label{e'1}
 (\alpha' \alpha)^* (f) = \alpha^* (\alpha')^* (f) \mbox{.}
\end{equation}
\begin{equation}
\label{e'2} (\alpha' \alpha)_* (G) = (\alpha')_* \alpha_* (G)
\mbox{.}
\end{equation}
\item
If $E_3$ and $L'$ are discrete
Lie groups, then $E_3 \mathop{\amalg}\limits_{E_2} H'$ is a
discrete  Lie group, and the following formulas are satisfied for any $f \in \s(E_1)$, $G \in \s'(H')$ :
\begin{equation} \label{pr17f5}
(\alpha' \alpha)_* (f) = (\alpha')_* \alpha_* (f) \mbox{.}
\end{equation}
\begin{equation} \label{pr17f6}
(\alpha' \alpha)^* (G) = \alpha^* (\alpha')^*(G) \mbox{.}
\end{equation}
\end{enumerate}
\end{prop}
\proof. Formula~(\ref{e'1}) follows from definitions of the
corresponding maps. Formula~(\ref{e'2}) is the conjugate formula to
 formula~(\ref{e'1}).

From an admissible triple from $C_0^{ar}$
$$
0 \arrow{e} E_3 \arrow{e,t}{\beta_{\alpha'}} E_3
\mathop{\amalg}\limits_{E_2} H' \arrow{e,t}{\beta_{\beta'}} L'
 \arrow{e} 0
$$
we obtain that if $E_3$ and $L'$ are discrete Lie groups, then
$E_3 \mathop{\amalg}\limits_{E_2} H'$ is also a discrete Lie group.
Now formula~(\ref{pr17f5})  follows at once from
definition. Formula~(\ref{pr17f6}) is the conjugate formula to
formula~(\ref{pr17f5}). The proposition is proved.

\vspace{0.5cm}

Let
\begin{equation}  \label{aaa}
0 \arrow{e} E_1  \arrow{e,t}{\alpha}  E_2 \arrow{e,t}{\beta}   E_3
\arrow{e} 0
\end{equation}
be an admissible triple from $C_0^{ar}$.
Then we have the following admissible triple from $C_0^{ar}$ of the Pontragin dual groups:
\begin{equation}  \label{bbb}
0 \arrow{e} \widehat{E_3} \arrow{e,t}{\hat{\beta}} \widehat{E_2}
\arrow{e,t}{\hat{\alpha}}  \widehat{E_1} \arrow{e} 0 \mbox{.}
\end{equation}

\begin{prop} \label{prpr}
We have the following commutative diagrams.

\begin{equation} \label{eq1}
\begin{diagram}
\node{\s(E_2) \otimes_{\sdc} \mu(E_2)} \arrow[3]{e,t}{\beta_*
\otimes \Id_{\mu(E_3)}}
\arrow{s,l}{\F} \node[3]{\s(E_3) \otimes_{\sdc} \mu(E_3)} \arrow{s,r}{\F} \\
\node{\s(\widehat{E_2})} \arrow[3]{e,b}{\hat{\beta}^*}
\node[3]{\s(\widehat{E_3})}
\end{diagram}
\end{equation}

\begin{equation}  \label{eq2}
\begin{diagram}
\node{\s(E_2)}  \arrow[3]{e,t}{\alpha^*}
\arrow{s,l}{\F \otimes \Id_{\mu(\widehat{E_2})}} \node[3]{\s(E_1)}
\arrow{s,r}{\F \otimes \Id_{\mu(\widehat{E_1})}}
\\
\node{\s(\widehat{E_2})\otimes_{\sdc} \mu(\widehat{E_2})}
\arrow[3]{e,t}{\hat{\alpha}_* \otimes \Id_{\mu(\widehat{E_1})}}
\node[3]{\s(\widehat{E_1}) \otimes_{\sdc} \mu(\widehat{E_1})}
\end{diagram}
\end{equation}

\begin{equation}  \label{eq3}
\begin{diagram}
\node{\s'(E_3) \otimes_{\sdc} \mu(\widehat{E_3})}
\arrow[3]{e,t}{\beta^* \otimes \Id_{\mu(\widehat{E_2})}}
\arrow{s,l}{\F} \node[3]{\s'(E_2) \otimes_{\sdc} \mu(\widehat{E_2})}
\arrow{s,r}{\F}
\\
\node{\s'(\widehat{E_3})} \arrow[3]{e,b}{\hat{\beta}_* }
\node[3]{\s'(\widehat{E_2})}
\end{diagram}
\end{equation}

\begin{equation}  \label{eq4}
\begin{diagram}
\node{\s'(E_1)}  \arrow[3]{e,t}{\alpha_*} \arrow{s,l}{\F \otimes \Id_{\mu(E_1)} } \node[3]{\s'(E_2)}  \arrow{s,r}{\F \otimes \Id_{\mu(E_2)}}
\\
\node{\s'(\widehat{E_1}) \otimes_{\sdc} \mu(E_1)}
\arrow[3]{e,b}{\hat{\alpha}^* \otimes \Id_{\mu(E_2)}}
\node[3]{\s'(\widehat{E_2}) \otimes_{\sdc} \mu(E_2)}
\end{diagram}
\end{equation}
\end{prop}
\proof. We use the natural isomorphisms~(\ref{isotr}) and~(\ref{isodu}) for admissible
triples~(\ref{aaa}), (\ref{bbb}) and pairs of objects $E_i$, $\widehat{E_i}$ $(1 \le i \le 3)$.
Now we have
\begin{equation} \label{ccc}
\mu_3(\beta_*(f \otimes \mu_1) \overline{\chi}) = \mu_3 (\beta_*(f \beta^*(\overline{\chi})\otimes \mu_1))
= \mu_2(f \overline{\beta^*(\chi)})
\mbox{,}
\end{equation}
where $\mu_i \in \mu(E_i)$ $(1 \le i \le 3)$, $\mu_3 = \mu_1 \otimes \mu_2$, $\chi \in \widehat{E_3}$.
The first equality in formula~(\ref{ccc})   follows
from the definition of map $\beta_*$, the second equality in formula~(\ref{ccc})
follows from Fubini's theorem. Formula~(\ref{ccc}) is equivalent to diagram~(\ref{eq1}).
Now diagram~(\ref{eq2}) follows from diagram~(\ref{eq1}) by using statement~\ref{it3ft}  of
 proposition~\ref{fourtra}.
Diagram~(\ref{eq3}) is the conjugate diagram to diagram~(\ref{eq1}).
Diagram~(\ref{eq4}) is the conjugate diagram to diagram~(\ref{eq2}).
The proposition is proved.

\begin{cons}[Poisson formula]
Let $\mu_1 \in \mu(E_1)$, $\mu_1 \ne 0$, and $\mu_3 \in \mu(\widehat{E_3})$.
Let $\delta_{E_1, \mu_1}  \eqdef \alpha_*(\mu_1) \in \s'(E_2)$, and analogously $\delta_{\widehat{E_3}, \mu_3} \eqdef
\hat{\beta}_*(\mu_3) \in \s'(\widehat{E_2})$. Then $\F_{\mu_1^{-1} \otimes \mu_3} (\delta_{E_1, \mu_1})= \delta_{\widehat{E_3}, \mu_3}$.
\end{cons}
\proof. Using isomorphisms~(\ref{isotr}) and~(\ref{isodu}) we obtain that $\mu_1^{-1} \otimes \mu_3  \in \mu(\widehat{E_2})$.
By definitions, we have $\delta_{E_1, \mu_1}= \beta^*(\delta_{E_3} \otimes \mu_1)$, where for any $E \in \Ob(C_0^{ar})$ we define $\delta_E \in \s'(E)$ as
$\delta_E(f) \eqdef f(0)$ for any $f \in \s(E)$. Analogously, by definitions, we have
$\delta_{\widehat{E_3}, \mu_3}= \hat{\alpha}^*(\delta_{\widehat{E_1}} \otimes \mu_3)$. Besides, for any $E \in \Ob(C_0^{ar})$,
$\mu \in \mu(E)$, $\mu \ne 0$ we have $\F_{\mu^{-1}}(\mu) = \delta_{\widehat{E}}$, $\F_{\mu^{-1}}(\delta_E) = \mu^{-1}$. Now the Poisson formula follows from diagram~(\ref{eq3}).

\begin{prop} \label{comdia}
We have the following commutative diagrams.
\begin{enumerate}
\item If $E_1$ is a compact Lie group, then
\begin{equation}  \label{eqq1}
\begin{diagram}
\node{\s(E_3) \otimes_{\sdc} \mu(E_3)}   \arrow[4]{e,t}{\beta^* \otimes \Id_{\mu(E_3)} \otimes { \bf 1_{E_1}} } \arrow{s,l}{\F}
\node[4]{\s(E_2) \otimes_{\sdc} \mu(E_2)}
 \arrow{s,r}{\F}
\\
\node{\s(\widehat{E_3})} \arrow[4]{e,b}{\hat{\beta}_* }
\node[4]{\s(\widehat{E_2})}
\end{diagram}
\end{equation}
where ${\bf 1_{E_1}} \in \mu(E_1)$, and for $1 \in \s(E_1)$, by definition, ${\bf 1_{E_1}}(1) = 1 \in \dc$;
\begin{equation} \label{eqq4}
\begin{diagram}
\node{\s'(E_2) \otimes_{\sdc} \mu(\widehat{E_2})} \arrow[4]{e,t}{\beta_* \otimes \Id_{\mu(\widehat{E_3})} \otimes {\bf 1_{E_1}}(\cdot) }
\arrow{s,l}{\F} \node[4]{\s'(E_3) \otimes_{\sdc} \mu(\widehat{E_3})} \arrow{s,r}{\F} \\
\node{\s'(\widehat{E_2})} \arrow[4]{e,b}{\hat{\beta}^*}
\node[4]{\s'(\widehat{E_3})}
\end{diagram}
\end{equation}
where the map ${\bf 1_{E_1}}(\cdot)  : \mu(\widehat{E_1}) \to \dc$ is given by the element ${\bf 1_{E_1}} \in \mu(E_1) \simeq \mu(\widehat{E_1})^*$.
\item
If $E_3$ is a discrete Lie group, then
\begin{equation}  \label{eqq2}
\begin{diagram}
\node{\s(E_1) \otimes_{\sdc} \mu(E_1)}  \arrow[4]{e,t}{\alpha_* \otimes \Id_{\mu(E_1)} \otimes {\bf \sum_{E_3}}} \arrow{s,l}{ \F }
\node[4]{\s(E_2) \otimes_{\sdc} \mu(E_2)}  \arrow{s,r}{\F}
\\
\node{\s(\widehat{E_1}) } \arrow[4]{e,b}{\hat{\alpha}^*}
\node[4]{\s(\widehat{E_2})}
\end{diagram}
\end{equation}
where ${\bf \sum_{E_3}} \in \mu(E_3)$, and for $f \in \s(E_3)$, by definition, ${\bf \sum_{E_3}}(f) = \sum\limits_{x \in E_3} f(x)$;
\begin{equation}  \label{eqq3}
\begin{diagram}
\node{\s'(E_2) \otimes_{\sdc} \mu(\widehat{E_2})}  \arrow[4]{e,t}{\alpha^*\otimes \Id_{\mu(\widehat{E_1})} \otimes {\bf \sum_{E_3}}(\cdot)} \arrow{s,l}{\F} \node[4]{\s'(E_1) \otimes_{\sdc} \mu(\widehat{E_1})}
\arrow{s,r}{\F}
\\
\node{\s'(\widehat{E_2})} \arrow[4]{e,b}{\hat{\alpha}_* }
\node[4]{\s'(\widehat{E_1}) }
\end{diagram}
\end{equation}
where the map ${\bf \sum_{E_3}}(\cdot)  : \mu(\widehat{E_3}) \to \dc$ is given by the element ${\bf \sum_{E_3}} \in \mu(E_3) \simeq \mu(\widehat{E_3})^*$.
\end{enumerate}
\end{prop}
\proof.
If $E \in \Ob(C_0^{ar})$ is a compact Lie group, then  $\widehat{E}$ is a discrete Lie group, and
${\bf 1_{E}} \otimes {\bf \sum_{\widehat{E}}} =1$ with respect to isomorphism~(\ref{isodu}). We use also isomorphism~(\ref{isotr}).
Diagram~(\ref{eqq1}) follows from
Fubini's theorem and the formula $\F_{{\bf 1_{E_1}}}(c)= c \delta$, where $c \in\dc$, $\delta \in \s(\widehat{E_1})$,
and $\delta(0)=1$, $\delta(x)=0$ if $x \in \widehat{E_1}$, $x \ne 0$.
Diagram~(\ref{eqq4}) is the conjugate diagram to  diagram~(\ref{eqq1}).
Diagram~(\ref{eqq2}) follows from diagram~(\ref{eqq1}) and statement~\ref{it3ft}  of proposition~\ref{fourtra}.
Diagram~(\ref{eqq3}) is the conjugate diagram to diagram~(\ref{eqq2}). The proposition is proved.

\subsection{Functions and distributions on objects of $C_1^{ar}$}
\label{subsec4.2}
 Let $E=(I,F,V) \in \Ob(C_1^{ar})$. By definition,
for any $i,j,k \in I$ such that $i \le j \le k$ we have an
admissible triple from $C_0^{ar}$
\begin{equation} \label{extri}
0 \arrow{e} F(j) / F(i) \arrow{e,t}{\alpha_{ijk}} F(k) / F(i)
\arrow{e,t}{\beta_{ijk}} F(k)/F(j) \arrow{e} 0 \mbox{.}
\end{equation}
Hence, we have that for any $i,j,k, l \in I$, $i \le j \le k \le l$
$$
F(k)/ F(i) \simeq (F(l)/ F(i)) \mathop{\times}\limits_{F(l)/F(j)}
(F(k)/F(j)) $$
$$F(l)/ F(j)  \simeq (F(l)/F(i))
\mathop{\amalg}\limits_{F(k)/F(i)} (F(k)/F(j))
$$
as objects of $C_0^{ar}$. Besides, from definition~\ref{def4'},
there exist $i_E \le j_E \in I$ such that for any $k \le i_E$ the
group $F(i_E)/F(k)$ is a finite abelian group, and for any $l \ge
j_E$ the group $F(l) / F(j_E)$ is a finite abelian group.

Therefore, using subsection~\ref{diim}, we have the following
well-defined  $\dc$-vector spaces:
\begin{equation} \label{s}
\s(E) \eqdef \mathop{\Lim_{\lto}}_{m \ge j_E}
\mathop{\Lim_{\lto}}_{n \le i_E} \s(F(m) / F(n)) =
\mathop{\Lim_{\lto}}_{n \le i_E} \mathop{\Lim_{\lto}}_{m \ge j_E}
\s(F(m) / F(n)) \mbox{,}
\end{equation}
where limits are taken with respect to the maps $\beta_{ijk}^*$ and
$(\alpha_{ijk})_*$;
\begin{equation} \label{s'}
\s'(E) \eqdef \mathop{\Lim_{\longleftarrow}}_{m \ge j_E}
\mathop{\Lim_{\longleftarrow}}_{n \le i_E} \s'(F(m) / F(n)) =
\mathop{\Lim_{\longleftarrow}}_{n \le i_E} \mathop{\Lim_{\longleftarrow}}_{m
\ge j_E} \s'(F(m) / F(n)) \mbox{,}
\end{equation}
where limits are taken with respect to the maps $(\beta_{ijk})_*$
and $\alpha_{ijk}^*$.

From this definition we have that the space $\s(E)$ is a
$\dc$-subalgebra of $\dc$-algebra $\f(V)$ of all $\dc$-valued
functions on $V$.

\begin{nt} {\em
A similar construction was used by F.~Bruhat to define the corresponding spaces of functions and distributions on locally compact abelian groups
(see~\cite{Br1}).
}
\end{nt}

\begin{lemma} \label{l2}
In notations of formula~(\ref{extri}) we have that for any  $i \le j
\le k \in I$ when the following maps are defined
$$(\alpha_{ijk})_* \; : \; \s(F(j)/F(i)) \to
\s(F(k)/F(i)) \mbox{,} \quad \beta_{ijk}^* \; : \; \s(F(k)/F(j)) \to
\s(F(k)/F(i)) \mbox{,}$$ these maps  are injective maps, and the
corresponding  maps
$$\alpha_{ijk}^* \; : \;
\s'(F(k)/F(i)) \to \s'(F(j)/F(i)) \mbox{,} \quad (\beta_{ijk})_* \;
: \; \s'(F(k)/F(i)) \to \s'(F(k)/F(j))$$ are surjective maps.
\end{lemma}
\proof. The injectivity of corresponding maps follows from
definitions. Besides, we have the map: $\alpha_{ijk}^*:
\s(F(k)/F(i)) \to \s(F(j)/F(i))$ such that $ \alpha_{ijk}^* \cdot
(\alpha_{ijk})_* = \Id_{\s(F(j)/F(i))}$. Therefore, as conjugate
equality we have $\alpha_{ijk}^* \cdot (\alpha_{ijk})_* =
\Id_{\s'(F(j)/F(i))}$. Hence, the map $\alpha_{ijk}^* :
\s'(F(k)/F(i)) \to \s'(F(j)/F(i))$ is a surjective map. Analogously,
we have the map:
$$(\beta_{ijk})_{*, {\bf 1_{F(j)/F(i)}}} :
\s(F(k)/F(i)) \to \s(F(k)/F(j)) $$
which is defined as
$$(\beta_{ijk})_{*, {\bf 1_{F(j)/F(i)}}}(f) \eqdef
(\beta_{ijk})_* (f \otimes {\bf 1_{F(j)/F(i)}}) \mbox{,}$$ where $f
\in \s(F(k)/F(i))$, ${\bf 1_{F(j)/F(i)}} \in \mu(F(j)/F(i))$, ${\bf
1_{F(j)/F(i)}}(1) =1$, because $F(j)/F(i)$ is a compact Lie group.
And we have that $(\beta_{ijk})_{*, {\bf 1_{F(j)/F(i)}}}  \cdot
\beta_{ijk}^* = \Id_{F(k)/F(j)}$. Taking the conjugate map and the
conjugate equality to the last equality, we obtain that  the map
$(\beta_{ijk})_* : \s'(F(k)/F(i)) \to \s'(F(k)/F(j))$ is a
surjective map. The lemma is proved.

\vspace{0.5cm}

For any $i \le j \in I$ we have the natural non-degenerate
$\dc$-bilinear  pairing $$<\cdot, \cdot>_{i,j} \quad : \quad
\s(F(j)/F(i)) \times \s'(F(j)/F(i)) \lto \dc \mbox{.}$$ Using
lemma~\ref{l2} and conjugate properties of direct and inverse images
on Schwartz spaces of $C_0^{ar}$, we obtain that there is a
non-degenerate $\dc$-bilinear pairing:
\begin{equation} \label{par}
<\cdot, \cdot> \quad : \quad \s(E) \times \s'(E) \lto \dc \mbox{.}
\end{equation}

\begin{example} \label{exam88} {\em
We consider the situation of example~\ref{exad}. We have a number
field $K$ such that $[K : \dq] =n$, and  we consider the "full"
adele ring of the field $K$
$$\da_K= \da_K^{fin} \times \prod\limits_{1 \le i \le l} K_{p_i} \mbox{,}$$
where $\da_K^{fin}$ is the "finite" adele ring of the field $K$, and
 $p_1, \ldots, p_l$ are all Archimedean places of the field $K$. We have
 $
 \prod\limits_{1 \le i \le l} K_{p_i} \simeq \dr^n
 $. We constructed in example~\ref{exad} the following object in the category~$C_1^{ar}$:
 $$
Q = (D_0 \cup D_1,  H \times F, \da_K )  \mbox{.}
$$
Then, by definitions, we have that
$$
\s(Q) = \s(\da_K^{fin} ) \otimes_{\sdc} \s(\prod\limits_{1 \le i \le l} K_{p_i}) \mbox{,}
$$
where $\s(\prod\limits_{1 \le i \le l} K_{p_i}) = \s(\dr^n)$ (see example~\ref{ex5}), and $\s(\da_K^{fin} ) = \D(\da_K^{fin})$
is the space of all locally constant functions with compact support on~$\da_K^{fin}$ (see~\cite[4.2]{OsipPar} for analogous situation of the category $C_1(\df_q)$.) Besides, any $f \in  \s(\da_K^{fin})$ is a finite linear combination of functions $\bigotimes\limits_{p} f_p$,
where the product is taken over all non-Archimedean places of the field $K$, $f_p \in \s(K_p)=\D(K_p)$ (in other words, $f_p$ is from the space of all locally constant functions with compact support on~$K_p$), and $f_p = \delta_{\oo_{K_p}} $ for almost all $p$, where $\delta_{\oo_{K_p}}(x) =1$ if $x \in \oo_{K_p}$, and $\delta_{\oo_{K_p}}(x) =0$ if $x \notin \oo_{K_p}$.
}
\end{example}

\vspace{0.5cm}

Let $E=(I,F,V) \in \Ob(C_1^{ar})$. For any $ i \le j \in I$ we
defined in section~\ref{ft00} the space $\mu(F(i)/F(j)) \in
\s'(F(i)/F(j))$ such that $\dim_{\sdc} \mu(F(i)/F(j)) =1$.

For any $i \le j \le k \in I$ we have that the map $(\beta_{ijk})_*$
defines an isomorphism between the spaces $\mu(F(k)/F(i))$ and
$\mu(F(k)/F(j))$ if $j \le i_E$. Indeed, if $\mu \in
\mu(F(k)/F(j))$, then, by formula~(\ref{isotr}), ${\bf
1_{F(j)/F(i)}} \otimes \mu \in \mu(F(k)/F(i))$ (see notation ${\bf
1_{F(j)/F(i)}}$ in proposition~\ref{comdia}), and $(\beta_{ijk})_*
({\bf 1_{F(j)/F(i)}} \otimes \mu) = \mu$.

Analogously,  for any $i \le j \le k \in I$ we have that the map
$\alpha_{ijk}^*$ defines an isomorphism between the spaces
$\mu(F(k)/F(i))$ and $\mu(F(j)/F(i))$ if $j \ge j_E$. Indeed, if
$\mu \in \mu(F(j)/F(i))$, then, by formula~(\ref{isotr}), $\mu
\otimes {\bf \sum_{F(k)/F(j)}} \in \mu(F(k)/F(i))$ (see notation
${\bf \sum_{F(k)/F(j)}}$ in proposition~\ref{comdia}), and
$\alpha_{ijk}^* ( \mu  \otimes {\bf \sum_{F(k)/F(j)}} ) = \mu$.

From these reasonings it follows that it is well-defined the following space:
\begin{equation} \label{mc1}
\mu(E)  \eqdef \mathop{\Lim_{\longleftarrow}}_{m \ge j_E}
\mathop{\Lim_{\longleftarrow}}_{n \le i_E} \mu(F(m) / F(n)) =
\mathop{\Lim_{\longleftarrow}}_{n \le i_E} \mathop{\Lim_{\longleftarrow}}_{m
\ge j_E} \mu(F(m) / F(n)) \mbox{,}
\end{equation}
where limits are taken with respect to the maps $(\beta_{ijk})_*$
and $\alpha_{ijk}^*$. And we have
\begin{equation} \label{dmc1}
\dim_{\sdc} \mu(E) =1 \mbox{.}
\end{equation}
\begin{nt} \label{meass}{\em
Let $E =(I,F,V) \in \Ob(C_1^{ar})$. Using formulas~(\ref{s})
and~(\ref{s'}) we can define operators $T_a$ ($a \in V$) on the
spaces $\s(E)$ and $\s'(E)$ starting from the corresponding
operators on $\s(F(j)/F(i))$ and $\s'(F(j)/F(i))$, where $i \le j
\in I$ (see the beginning of subsection~\ref{ft00}). Then it is easy
to see that $\mu(E) = \s'(E)^V$. }
\end{nt}

Let $E = (I,F,V) \in \Ob(C_1^{ar})$. Then for any $i \le j \in I$ we have the natural isomorphism
$\mu(F(j)/F(i)) \otimes_{\sdc} \mu(\widehat{F(j)/F(i)}) \simeq \dc$ (see~(\ref{isodu})). Besides, if  $j \le i_E$  then
${\bf 1_{F(j)/F(i)}} \otimes {\bf \sum_{\widehat{F(j)/F(i)}}}=1$  under this isomorphism. Now from formulas~(\ref{f6})-(\ref{f7})
and formulas~(\ref{mc1})-(\ref{dmc1}) it follows that there is the following natural isomorphism:
\begin{equation} \label{isodu1}
\mu(E) \otimes_{\sdc} \mu(\check{E}) \lto \dc \mbox{.}
\end{equation}

Let
$$
0 \arrow{e} E_1 \arrow{e,t}{\alpha} E_2 \arrow{e,t}{\beta} E_3 \arrow{e} 0
$$
be an admissible triple from $C_1^{ar}$, where $E_i=(I_i, F_i, V_i)$ ($1 \le i \le 3$). Then by~(\ref{isotr}) we have the natural isomorphism for any
$i \le j \in I_2$:
\begin{equation} \label{equ}
\mu \left(\frac{F_2(j) \cap V_1}{  F_2(i) \cap V_1} \right) \otimes_{\sdc} \mu\left(\frac{\beta(F_2(j))}{\beta(F_2(i)}\right) \lto
\mu\left(\frac{F_2(j)}{F_2(i)}\right)
\mbox{.}
\end{equation}
Besides, if $j \le i_E$ then
$$
{\bf 1_{\frac{F_2(j) \cap V_1}{  F_2(i) \cap V_1}}} \otimes {\bf 1_{\frac{\beta(F_2(j))}{\beta(F_2(i)}}} = {\bf 1_\frac{F_2(j)}{F_2(i)}} \mbox{,}
$$
and if $i \ge j_E$ then
$$
{\bf \sum\nolimits_{\frac{F_2(j) \cap V_1}{  F_2(i) \cap V_1}}} \otimes {\bf \sum\nolimits_{\frac{\beta(F_2(j))}{\beta(F_2(i)}}} = {\bf \sum\nolimits_\frac{F_2(j)}{F_2(i)}}
$$
under isomorphism~(\ref{equ}).

Therefore from formulas~(\ref{mc1})-(\ref{dmc1}) it follows that there is the following natural isomorphism:
\begin{equation} \label{isotr1}
\mu(E_1) \otimes_{\sdc} \mu(E_3) \lto \mu(E_2) \mbox{.}
\end{equation}

\vspace{0.5cm}

Let $E=(I,F,V) \in \Ob(C_1^{ar})$. For any $i \le j \in I$ we have the Fourier transforms:
$$ \F_{i,j} \; : \; \s(F(j)/F(i)) \otimes_{\sdc} \mu(F(j)/F(i)) \lto \s(\widehat{F(j)/F(i)})  \quad  \mbox{and}$$
$$ \F_{i,j} \; : \; \s'(F(j)/F(i)) \otimes_{\sdc} \mu(\widehat{F(j)/F(i)}) \lto \s'(\widehat{F(j)/F(i)}) \mbox{.}$$
Therefore, after taking the corresponding limits, using
formulas~(\ref{s})-(\ref{s'}), (\ref{mc1})-(\ref{dmc1}), and
proposition~\ref{comdia} we  obtain a well-defined Fourier
transforms ({\em one-dimensional Fourier transforms})  from the maps
$\F_{i,j}$:
$$
\F \; : \; \s(E) \otimes_{\sdc} \mu(E) \lto \s(\check{E}) \quad \mbox{and}
$$
$$
\F \; : \; \s'(E) \otimes_{\sdc} \mu(\check{E}) \lto \s'(\check{E}) \mbox{.}
$$

\begin{nt}{\em
It is easy to see that this definition of Fourier transform $\F : \s(E) \otimes_{\sdc} \mu(E) \to \s(\check{E})$
coincides with the analog of formula~(\ref{defftr}):
$$
\F (f \otimes \mu)= \F_{\mu}(f) \mbox{,}
$$
$$
\F_{\mu} (f)(b) = \mu (f \overline{b}) \mbox{,}
$$
where $f \in \s(E)$, $\mu \in \mu(E)$, $b \in \check{V} \subset \Hom (V, \dt) $,
and we consider $\s(\check{E})$ as $\dc$-subspace of all $\dc$-valued functions on $\check{V}$.
}
\end{nt}

Let $E =(I,F,V) \in \Ob(C_1^{ar})$, let $f \in \s(E)$, $H \in
\s'(E)$. Using formulas~(\ref{ch1}-\ref{ch2}) and taking
corresponding limits in formulas~(\ref{s})-(\ref{s'}), we can define
$\check{f} \in \s(E)$, $\check{H} \in \s'(E)$ such that
$\check{\check{f}}=f$ and $\check{\check{H}} =H$. For any $i \le j
\in I$ the spaces $\s(F(j)/F(i))$ and $\s'(F(j)/F(i))$ are
topological $\dc$-vector spaces. Therefore the spaces $\s(E)$ and
$\s'(E)$ are also topological $\dc$-vector spaces which have
topologies of inductive and projective limits correspondingly. Now
directly from proposition~\ref{fourtra}  and taking corresponding
limits in formulas~(\ref{s})-(\ref{s'}) we obtain the generalization
of proposition~\ref{fourtra}.
\begin{prop} \label{fourtra1}
Let $E \in \Ob(C_1^{ar})$, $\mu \in \mu(E)$, $\mu \ne 0$. The
following properties are satisfied.
\begin{enumerate}
\item The map $\F_{\mu}$ is an isomorphism of topological $\dc$-vector spaces $\s(E)$ and $\s(\check{E})$.
\item The map $\F_{\mu^{-1}}$ is an isomorphism of topological $\dc$-vector spaces $\s'(E)$ and $\s'(\check{E})$.
\item For any $f \in \s(E)$ and $G \in \s'(\check{E})$
$$
<\F(f), G> = <f, \F(G)> \mbox{.}
$$
If $E$ is a complete object, then
\item   $ \F_{\mu^{-1}} \F_{\mu}(f) = \check{f}$ for any $f \in \s(E)$,
\item   $ \F_{\mu}   \F_{\mu^{-1}} (H) =\check{H}$ for any $H \in \s'(E)$.
\end{enumerate}
\end{prop}

\vspace{0.5cm}
We consider an admissible triple from~$C_1^{ar}$:
$$
0 \arrow{e} E_1 \arrow{e,t}{\alpha} E_2 \arrow{e,t}{\beta} E_3 \arrow{e} 0 \mbox{,}
$$
 where $E_i=(I_i, F_i, V_i)$ ($1 \le i \le 3$). In subsection~\ref{diim} we constructed direct and inverse images
 on functions and distributions of objects from $C_0^{ar}$. Using these definitions, formulas~(\ref{s})-(\ref{s'}), formulas~(\ref{mc1})-(\ref{dmc1}), and taking the corresponding inductive or projective limits, we obtain directly the following well-defined maps.
 \begin{enumerate}
 \item \label{I1} $\beta_* : \s(E_2) \otimes_{\sdc} \mu(E_1)  \to \s(E_3)$.
 \item \label{I2} $\beta^*: \s'(E_3) \otimes_{\sdc} \mu(E_1) \to \s'(E_2)$.
 \item \label{I3} $\alpha^* : \s(E_2) \to \s(E_1)$.
 \item \label{I4}  $\alpha_*: \s'(E_1)  \to \s'(E_2)$.
 \item \label{I5}  $\beta^* : \s(E_3) \to \s(E_2) $ if $E_1$ is a compact object.
 \item \label{I6}   $\beta_*: \s'(E_2)  \to \s'(E_3)$ if $E_1$ is a compact object.
 \item \label{I7}   $\alpha_* : \s(E_1) \to \s(E_2) $ if $E_3$ is a discrete object.
 \item \label{I8} $\alpha^*: \s'(E_2)  \to \s'(E_1)$  if $E_3$ is a discrete object.
 \end{enumerate}

We note that, by constructions, the maps from items~\ref{I1} and~\ref{I2}
are conjugate maps with respect to pairing~(\ref{par}). The maps from items~\ref{I3} and~\ref{I4}
are conjugate maps with respect to pairing~(\ref{par}). The maps from items~\ref{I5} and~\ref{I6}
are conjugate maps with respect to pairing~(\ref{par}). The maps from items~\ref{I7} and~\ref{I8}
are conjugate maps with respect to pairing~(\ref{par}).

\begin{nt}
\label{remm8} {\em Directly from constructions above of direct and
inverse images on functions and distributions of objects of the
category $C_1^{ar}$ and from propositions of subsection~\ref{diim}
we obtain that analogous propositions are valid for objects of
category $C_1^{ar}$: propositions~\ref{pr18}, \ref{pr16},
\ref{pr17}, \ref{prpr}, \ref{comdia} and Poisson formula. We have to
preserve the formulations of these propositions, but change in them
objects of $C_0^{ar}$ to objects of $C_1^{ar}$, compact Lie groups
to compact objects from $C_1^{ar}$, and discrete Lie groups to
discrete objects from $C_1^{ar}$.

Besides, if $E=(I,F,V) $ is a compact object from $C_1^{ar}$, then there is a natural element ${\bf 1_{E}} \in \mu(E)$
which is the projective limit of elements~${\bf 1_{F(j)/F(i)}} \in \mu(F(j)/F(i))$ ($i \le j \in I$). If $E=(I,F,V) $ is a discrete object from $C_1^{ar}$, then there is a natural element ${\bf \sum_{E}} \in \mu(E)$
which is the projective limit of elements~${\bf \sum_{F(j)/F(i)}} \in \mu(F(j)/F(i))$ ($i \le j \in I$).
}
\end{nt}

\begin{nt} \label{rem9}{\em
Let $E_2=(I_2,F_2,V_2) \in \Ob(C_1^{ar}) $. We consider another
object $E_1=(I_1,F_1,V_1) \in \Ob(C_1^{ar})$ which will dominate
object $E_2$ as filtered abelian group and such that for any $i \le
j \in I_2 $ we have the same structure of object of $C_0^{ar}$ on
$F_2(j)/F_2(i)$ and on $F_1(\phi(j))/F_1(\phi(i))$ (see
definition~\ref{defdom}.) In this case we will say that the object
$E_1$ {\em dominates} the object $E_2$.

\begin{defin} \label{defff}
Let $(I, F, V), (J, G, V) \in \Ob(C_1^{ar})$. We say that there is
an {\it equivalence} $(I, F, V) \sim (J, G, V)$ iff there is a
collection of objects $(I_l, F_l, V) \in \Ob(C_1^{ar})$, $1 \le l
\le n$
 such that
\begin{itemize}
\item[(i)] $I_1= I$, $F_1 = F$, $I_n = J$, $F_n =G$,
\item[(ii)] for any $1 \le l \le n-1$ either the object
$(I_l, F_l, V)$ dominates the object $(I_{l+1}, F_{l+1}, V)$, or the
object $(I_{l+1}, F_{l+1}, V)$ dominates the object $(I_l, F_l, V)$.
\end{itemize}
\end{defin}

Immediately from definitions we obtain that if $E_1 \sim E_2$, then
$\s(E_1) =\s(E_2)$, $\s'(E_1)=\s'(E_2)$,  $\mu(E_1)=\mu(E_2)$.

These isomorphisms  are well-defined with respect to
isomorphisms~(\ref{isotr}) and~(\ref{isodu}), constructions of
direct and inverse images, and Fourier transform. (More exactly, it
is possible to write the corresponding commutative diagrams.)}
\end{nt}

\section{Categories $C_2^{ar}$}
\label{c_2-spaces}

\begin{defin}
 Objects of the category $C_2^{ar}$, i.e. $\Ob(C_2^{ar})$, are filtered
abelian groups $(I, F, V)$ (see definition~\ref{def1}) with the
following additional structures
\begin{itemize}
\item[(i)] for any $i \le j \in I$  on abelian groups $F(j) / F(i)$
it is given a structure $E_{i,j} \in \Ob(C_1^{ar})$,
\item[(ii)] for any $i \le j \le k \in I$
$$
0 \lto E_{i,j}  \lto E_{i,k}  \lto  E_{j,k} \lto 0
$$
is an admissible triple from $C_1^{ar}$.
\end{itemize}
\end{defin}

\begin{defin} \label{morc2ar}
Let $E_1 = (I_1, F_1, V_1)$ and $E_2 = (I_2, F_2, V_2)$ be from
$\Ob(C_2^{ar})$. Then, by definition, $\Mor_{C_2^{ar}}(E_1, E_2)$
consists of elements $A \in \Hom_k (V_1, V_2)$ such that the
following conditions hold:
\begin{itemize}
\item[(i)] \label{i1_2} for any $i \in I_1$ there is an $j \in I_2$ such that  $A (F_1(i)) \subset F_2(j)$,
\item[(ii)] \label{i2_2}  for any $j \in I_2$ there is an $i \in I_1$ such that  $A (F_1(i)) \subset F_2(j)$,
\item[(iii)] \label{i3_2}  for any $i_1 \le i_2 \in I_1$ and $j_1 \le j_2 \in I_2$ such that $A (F_1(i_1)) \subset F_2(j_1)$
and $A (F_1(i_2)) \subset F_2(j_2)$ we have that the induced
$k$-linear map
$$   \bar{A} : \frac{F_1(i_2)}{F_1(i_1)} \lto \frac{F_2(j_2)}{F_2(j_1)}
$$
is an element from
$$\Mor\nolimits_{C_1^{ar}}(\frac{F_1(i_2)}{F_1(i_1)}, \frac{F_2(j_2)}{F_2(j_1)})  \mbox{.}$$
\end{itemize}
\end{defin}

The notion of morphisms in $C_2^{ar}$ is well-defined as it follows
from the following proposition.
\begin{prop} \label{propmor3}
 Let $E_1= (I_1, F_1, V_1)$, $E_2= (I_2, F_2,
V_2)$, $E_1'$, $E_2'$ be from $Ob(C_2^{ar})$, and  $A$ is a  map
from $\Hom (V_1, V_2)$.
\begin{enumerate}
\item If the filtered abelian group $E_1$ dominates the filtered abelian group
$E'_1$, and the filtered abelian group $E_2$ dominates the filtered
abelian group $E'_2 $, then $A \in \Mor_{C_2^{ar}}(E_1, E_2)$ if and
only if $A \in \Mor_{C_2^{ar}}(E'_1, E'_2)$.
\item $\Mor_{C_2^{ar}}(E_1, E_2)$ is an abelian subgroup of $\Hom(V_1, V_2)$.
\item If $E_3$ is an object of $C_2^{ar}$, then
$$     \Mor_{C_2^{ar}}(E_2, E_3)   \circ   \Mor_{C_2^{ar}}(E_1, E_2)  \subset  \Mor_{C_2^{ar}}(E_1, E_3) \mbox{.}$$
\end{enumerate}
\end{prop}
\proof of this proposition is fully analogous to the proof of
proposition~2.1 from~\cite{Osip} and can be done by induction, see
propostion~\ref{propmor2}. The proposition is proved.

\begin{defin}
Let $E_1 = (I_1, F_1, V_1)$, $E_2 = (I_2, F_2, V_2)$ and $E_3=(I_3,
F_3, V_3)$ be from $\Ob(C_2^{ar})$. Then we say that
\begin{equation} \label{adt'}
0 \arrow{e} E_1 \arrow{e,t}{\alpha} E_2 \arrow{e,t}{\beta} E_3
\arrow{e} 0
\end{equation}
is an admissible triple from $C_2^{ar}$ when
 the following conditions are satisfied:
\begin{itemize}
\item[(i)]
$$
0 \arrow{e} V_1 \arrow{e,t}{\alpha} V_2 \arrow{e,t}{\beta} V_3
\arrow{e} 0
$$
is an exact triple of abelian groups,
\item[(ii)] \label{itaa'}
the filtration $(I_1, F_1, V_1)$ dominates the filtration $(I_2,
F'_1, V_1)$, where  $F'_1 (i) = F_2(i) \cap V_1$  for any $i \in
I_2$,
\item[(iii)] \label{itbb'}
the filtration $(I_3, F_3, V_3)$ dominates the filtration $(I_2,
F'_3, V_3)$, where  $F'_3(i) =  F_2(i) / F_2(i) \cap V_1$ for any $i
\in I_2$,
\item[(iv)] for any $i \le j \in I_2$
\begin{equation} \label{trojkaa3}
0 \lto \frac{F'_1(j)}{F'_1(i)} \lto \frac{F_2(j)}{F_2(i)}  \lto
\frac{F'_3(j)}{F'_3(i)} \lto 0
\end{equation}
is an admissible triple from $C_1^{ar}$. (By definition of
$\Ob(C_2^{ar})$, on every abelian group from triple~(\ref{trojkaa3})
it is given the structure of element from $\Ob(C_1^{ar})$).
\end{itemize}
\end{defin}

We say that $\alpha$ from an admissible triple~(\ref{adt'}) is {\em
an admissible monomorphism}, and $\beta$ from an admissible
triple~(\ref{adt'}) is {\em an admissible epimorphism.}

\begin{example} \label{examp4} {\em
The category $C_2^{fin}$ is a full subcategory of the category
$C_2^{ar}$ such that a triple from $C_2^{fin}$ is an admissible
triple from $C_2^{fin}$ if and only if it is an admissible triple in
$C_2^{ar}$.

The category $C_1^{ar}$ is a full subcategory of the category
$C_2^{ar}$ with respect to the following functor $I$: for any $G \in
C_1^{ar}$ we put $I(G)= (\{(0), (1)\}, F, G )$, where $(0) < (1)$,
$F((0))=e$, $F((1))=G$, $e$ is the trivial subgroup of $G$, and $I$
acts as identity map on the morphisms from $C_1^{ar}$. }
\end{example}

We say that an object $E_1 = (I_1, F_1, V) \in \Ob(C_2^{ar})$ {\em
dominates} another object  $E_2 = (I_2, F_2, V) \in  \Ob(C_2^{ar})$
if the following conditions are satisfied:
\begin{itemize}
\item[(i)] the filtration $(I_1, F_1, V)$ dominates the filtration
$(I_2, F_2, V)$,
\item[(ii)] for any $i \le j \in I_2$  the object ${E_1}_{i,j}=F_2(j)/F_2(i) \in \Ob(C_1^{ar})$
with filtration induced by $C_2^{ar}$-structure on $E_1$ dominates
the object ${E_2}_{i,j}= F_2(j)/F_2(i) \in \Ob(C_1^{ar})$ with
filtration induced by $C_2^{ar}$-structure on $E_2$ (see
remark~\ref{rem9}).
\end{itemize}

We consider the group $\Aut_{C_2^{ar}}(E) \eqdef \Mor_{C_2^{ar}}(E,
E)^*$ (i.e., invertible elements in the ring $\Mor_{C_2^{ar}}(E,
E)$) for any  $E \in \Ob(C_2^{ar})$. By proposition~\ref{propmor3},
if  an object $E_1 \in \Ob(C_2^{ar})$ dominates an object  $E_2 \in
\Ob(C_2^{ar})$, then canonically
$$\Aut\nolimits_{C_2^{ar}}(E_1)
=\Aut\nolimits_{C_2^{ar}}(E_2) \mbox{.}$$

Let $E = (I, F, V) \in \Ob(C_2^{ar})$. We define the dual object
$\check{E} = (I^0, F^0, \check{V}) \in \Ob(C_2^{ar})$ in the
following way. The abelian group $\check{V} \subset \Hom(V, \dt )$
is defined as
\begin{equation}
\check {V} \eqdef \mathop{\Lim_{\rightarrow}}_{j \in I}
\mathop{\Lim_{\leftarrow}}_{i \ge j} \check{E}_{j,i} \mbox{,}
\end{equation}
where the object $\check{E}_{j,i} \in \Ob(C_1^{ar})$ is constructed
from the object $E_{j,i}= F(i)/F(j) \in \Ob(C_1^{ar})$ (see
subsection~\ref{c0c1}). The set $I^0$ is a partially ordered set,
which has the same set as $I$, but with the inverse order than $I$.
For $j \in I^0$ the subgroup
\begin{equation}
F^0(j) \eqdef \mathop{\Lim_{\leftarrow}}_{i \le j \in I^0}
\check{E}_{j,i} \quad \subset \quad \check{V} \mbox{.}
\end{equation}

If $E_1, E_2 \in \Ob(C_2^{ar})$ and $\theta \in
\Mor_{C_2^{ar}}(E_1,E_2)$, then there is canonically $\check{\theta}
\in \Mor_{C_2^{ar}} (\check{E_2}, \check{E_1})$. If
$$
0 \arrow{e} E_1  \arrow{e,t}{\alpha}  E_2 \arrow{e,t}{\beta}   E_3
\arrow{e} 0
$$
is an admissible triple from $C_2^{ar}$, then there is canonically
the following admissible triple from $C_2^{ar}$:
$$
0 \arrow{e} \check{E_3} \arrow{e,t}{\check{\beta}} \check{E_2}
\arrow{e,t}{\check{\alpha}}  \check{E_1} \arrow{e} 0 \mbox{.}
$$

\begin{lemma} \label{lbc2}
Let
$$
0 \arrow{e} E_1  \arrow{e,t}{\alpha}  E_2 \arrow{e,t}{\beta}   E_3
\arrow{e} 0
$$
be an admissible triple from $C_2^{ar}$.
\begin{enumerate}
\item Let $D \in \Ob(C_2^{ar})$, and $\gamma \in
\Mor_{C_2^{ar}} (D, E_3)$. Then there is the following admissible
triple from $C_2^{ar}$
\begin{equation} \label{adm12}
0 \arrow{e} E_1  \arrow{e,t}{\gamma_{\alpha}} E_2
\mathop{\times}_{E_3} D \arrow{e,t}{\gamma_{\beta}} D \arrow{e} 0
\end{equation}
and $\beta_{\gamma} \in \Mor_{C_2^{ar}}(E_2
\mathop{\times}\limits_{E_3} D, E_2)$ such that the following
diagram is commutative:
\begin{equation} \label{adm22}
\begin{diagram}
\node{0} \arrow{e} \node{E_1}  \arrow{e,t}{\gamma_{\alpha}}
\arrow{s,=} \node{E_2 \mathop{\times}_{E_3} D}
\arrow{e,t}{\gamma_{\beta}}
\arrow{s,l}{\beta_{\gamma}} \node{D} \arrow{e} \arrow{s,r}{\gamma} \node{0} \\
\node{0} \arrow{e} \node{E_1}  \arrow{e,t}{\alpha}  \node{E_2}
\arrow{e,t}{\beta} \node{E_3} \arrow{e} \node{0}
\end{diagram}
\end{equation}
\item Let $A \in \Ob(C_2^{ar})$, and $\theta \in
\Mor_{C_2^{ar}}(E_1, A)$.
 Then there is the following admissible triple from
$C_2^{ar}$
\begin{equation}  \label{adm'12}
0 \arrow{e} A  \arrow{e,t}{\theta_{\alpha}} A  \mathop{\amalg}_{E_1}
E_2 \arrow{e,t}{\theta_{\beta}} E_3 \arrow{e} 0
\end{equation}
and $\alpha_{\theta} \in \Mor_{C_2^{ar}}(E_2, A
\mathop{\amalg}\limits_{E_1} E_2)$ such that the following diagram
is commutative:
\begin{equation}  \label{adm'22}
\begin{diagram}
\node{0} \arrow{e} \node{E_1}  \arrow{e,t}{\alpha}
\arrow{s,l}{\theta} \node{E_2} \arrow{e,t}{\beta}
\arrow{s,r}{\alpha_{\theta}} \node{E_3} \arrow{e} \arrow{s,=}  \node{0} \\
\node{0} \arrow{e} \node{A}  \arrow{e,t}{\theta_{\alpha}} \node{A
\mathop{\amalg}\limits_{E_1} E_2} \arrow{e,t}{\theta_{\beta}}
\node{E_3} \arrow{e} \node{0}
\end{diagram}
\end{equation}
\end{enumerate}
\end{lemma}
\proof~follows by repeating of reasonings from the proof of
lemma~\ref{lbc} when in this proof the fibered product (or the
amalgamed sum) of quotients of  filtration is understood as the
fibered product (or the amalgamed sum) of objects from $C_1^{ar}$,
constructed in lemma~\ref{lbc}. The lemma is proved. \vspace{0.3cm}

\begin{defin}
We say that  $(I_1, F_1, V_1) \in \Ob(C_2^{ar})$ is a $c$-object iff
there is an $i \in I_1$ such that $ F_1(i) = V_1$.

We say that  $(I_2, F_2, V_2) \in \Ob(C_2^{ar})$ is an $d$-object
iff there is an $j \in I_2$ such that $F_2(j) = \{0\}$.
\end{defin}

It follows from definition that if $E_1 \in \Ob(C_2^{ar})$ is a
$c$-object, then $\check{E_1} \in \Ob(C_2)$ is a $d$-object. If $E_2
\in \Ob(C_2^{ar})$ is a $d$-object, then $\check{E_2} \in
\Ob(C_2^{ar})$ is a $c$-object.

\begin{defin}
We say that  $(I_1, F_1, V_1) \in \Ob(C_2^{ar})$ is a $cf$-object if
for any $i_1 \ge j_1 \in I_1$ the object $F(i_1)/F(j_1) \in
\Ob(C_1^{ar})$ is a compact object.

We say that  $(I_2, F_2, V_2) \in \Ob(C_2^{ar})$ is a $df$-object if
for any $i_2 \ge j_2 \in I_2$ the object $F(i_2)/F(j_2) \in
\Ob(C_1^{ar})$ is a discrete object.
\end{defin}

It follows from definition that if $E_1 \in \Ob(C_1^{ar})$ is a
$cf$-object, then $\check{E_1} \in \Ob(C_1^{ar})$ is a $df$-object.
If $E_2 \in \Ob(C_1^{ar})$ is a $df$-object, then $\check{E_2} \in
\Ob(C_2^{ar})$ is a $cf$-object.

\begin{defin}
We say that $(I,F,V) \in \Ob(C_2^{ar})$ is a complete object if the
following two conditions are satisfied.
\begin{enumerate}
\item
For any $i \ge j \in I$ the object $E_{j,i}= F(i)/F(j) \in C_1^{ar}$
is a complete object.
\item
$$
V = \mathop{\Lim_{\to}}_{i \in I}  \mathop{\Lim_{\gets}}_{j \le i}
F(i)/F(j)  \mbox{.}
$$
\end{enumerate}
\end{defin}

It follows from definitions that $E \in \Ob(C_2^{ar})$ is a complete
object iff $\check{\check{E}} = E$.

We have the obvious functor of completion  $\Omega : C_2^{ar} \to
C_2^{ar}$, where for $E = (I,F,V) \in \Ob(C_2^{ar})$ we define
$$
\Omega(E) \eqdef (I, F', \mathop{\Lim_{\to}}_{i \in I}
\mathop{\Lim_{\gets}}_{j \le i} \Psi(F(i)/F(j))) \mbox{,}
$$
where for any $j \in I$ we define
$$
F'(j) \eqdef \mathop{\Lim_{\gets}}_{j \le i} \Psi(F(i)/F(j))
\mbox{,}
$$
and $\Psi$ is a completion functor on $C_1^{ar}$ (see
subsection~\ref{c0c1}). The map $\Omega$ is easily extended to
$\Mor_{C_2^{ar}}(E_1, E_2)$, $E_1, E_2 \in \Ob(C_2^{ar})$ from
$\Mor_{C_1^{ar}}$.

It is clear from definition that for any $E \in \Ob(C_2^{ar})$ the
object $\Omega(E)$ is a complete object, and $\Omega^2 = \Omega$.
Moreover, $\Omega(E)= \check{\check{E}}$.

\begin{example} \label{eexx} {\em
We consider the field $K = \dr((t))$ (or the field $K = \dc((t))$).
Since any finite-dimensional $\dr$ (or $\dc$)-vector space is an
object of category $C_0^{ar}$, which is an object of category
$C_1^{ar}$ by example~\ref{examp1}, we have that $K$ has the natural
structure of object of category $C_2^{ar}$ with filtration given by
fractional ideals in the discrete valuation field $K$.
 }
\end{example}

\begin{example} \label{arsur} {\em
Let $X$ be an integral scheme of finite type over $\dz$, and $\dim X
=2 $. We suppose that there is a projective surjective morphism $X
\to \Spec E$, where $E$ is the ring of integers in a number field
$K$, $[K : \dq] = n$. It means that $X$ is an "arithmetic surface".

Now according to examples~\ref{exam1} and~\ref{examp4} the ring
$\da_X \in \Ob(C_2^{fin}) \subset \Ob(C_2^{ar})$ (compare also with
example~2.9 from~\cite{Osi}).

Let $p_1, \ldots, p_l$ be all Archimedean places of the field $K$.
Then for every $1 \le i \le l$ the completion field $K_{p_i}$ is
isomorphic either the field $\dc$, or the field $\dr$. We consider a
 curve
$$X_K \eqdef X
\otimes_{\Spec E} \Spec K \mbox{,}$$ which is  the generic  fibre of
the morphism from $X$ to $\Spec E$. We define the following  ring
\begin{equation} \label{eedd}
\da_{X, \infty} \eqdef \da_{X_K}
\mathop{\widehat{\otimes}}\nolimits_{K} \; ( \prod\limits_{1 \le i
\le l} K_{p_i} ) = \mathop{\mathop{\Lim}\limits_{\lto}}\limits_{D_1}
\mathop{\mathop{\Lim}\limits_{\longleftarrow}}\limits_{D_2 \le D_1}
\left(\frac{\da_{X_K}(D_1)}{\da_{X_K}(D_2)} \otimes_K (
\prod\limits_{1 \le i \le l} K_{p_i} ) \right) \mbox{,}
\end{equation}
where $\da_{X_K}$ is the adelic ring of the curve $X_K$,  $D_2 \le
D_1$ runs over all pairs of corresponding Cartier divisors on the
curve $X_K$, and for a Cartier divisor $D$ on $X_K$ the $K$-vector
space
$$
\da_{X_K}(D) \eqdef \prod\limits_{q \in X_K} \hat{\oo}_q \otimes_{
\oo_{X_K}} \oo(D) \mbox{.}
$$
Here the product is taken over all closed points of $X_K$, the ring
$\hat{\oo}_q$ is the completion of the local ring of point $q \in
X_K$ with respect to the maximal ideal, the invertible sheaf
$\oo(D)$ on $X_K$ corresponds to the Cartier divisor $D$ on $X_K$.

We note that
$$
\da_{X_K} =\mathop{\mathop{\Lim}\limits_{\lto}}\limits_{D_1}
\mathop{\mathop{\Lim}\limits_{\longleftarrow}}\limits_{D_2 \le D_1}
\frac{\da_{X_K}(D_1)}{\da_{X_K}(D_2)} \mbox{,}
$$
and the space $\frac{\da_{X_K}(D_1)}{\da_{X_K}(D_2)}$ is a
finite-dimensional vector space over the field $K$. Using the
isomorphism $\prod\limits_{1 \le i \le l} K_{p_i} \simeq \dr^n$ we
obtain that $ \da_{X, \infty} \in \Ob(C_2^{ar})$, because we can
take the filtration of $\da_{X, \infty}$ by Cartier divisors on the
curve $X_K$ (see formula~(\ref{eedd})), and quotients of this
filtration are finite-dimensional vector spaces over $\dr$ and
belong to~$\Ob(C_1^{ar})$ according to examples~\ref{examp1}
and~\ref{examp4}.

 Now we define  the arithmetic  adelic ring of $X$:
$$
\da_X^{ar} \: \eqdef \: \da_X \times \da_{X, \infty} \mbox{.}$$ We
have that $\da_X^{ar} \in \Ob(C_2^{ar})$, because we take the
product filtration on $\da_X^{ar}$, which is  induced by filtrations
on $\da_X$ and $\da_{X, \infty}$.

We note that
$$
\da_X \subset \prod\limits_{x \in C} K_{x,C} \mbox{,}
$$
where this product is taken over all pairs:  integer one-dimensional
subschemes $C$ of $X$ and  closed points $p$ on $C$. The ring
$K_{x,C}$ is a finite product of two-dimensional local fields. If
$x$ is a regular point of $C$ and of $X$, then $K_{x,C}$ is a
two-dimensional local field. (See, for example,~\cite[th.~2.10,
prop.~3.4]{Osi}).

We note also that the closed points of the curve $X_K$ are in
one-to-one  correspondence with irreducible "horizontal" curves on
$X$, i.e., with integer one-dimensional subschemes of $X$ that are
surjectively mapped onto $\Spec E$. (Such an irreducible
"horizontal" curve is the closure of a point of generic fibre.)

From this point of view we consider the definition of the ring
$\da_X^{ar}$ as addition  of  local fields $\dr((t))$ or $\dc((t))$
to the ring $\da_X$ such that these fields correspond to pairs:
irreducible "horizontal" curves on $X$ and infinite points (or
Archimedean places) on these curves. It follows from the following
well-known formula:
$$
L \otimes_K K_p = \prod\limits_q L_q \mbox{,}
$$
where the field $L$ is a finite extension of the field $K$, $p$ is
an Archimedean place of $K$, and the product is taken over all
Archimedean places $q$ of $L$ which are over $p$.}
\end{example}

\begin{nt}{\em
In a way, which is similar to the definition of categories
$C_n^{fin}$ (see definitions~\ref{def4},~\ref{d1}), we could define
categories $C_n^{ar}$, starting by induction from categories
$C_1^{ar}$. Further in the article we will apply categories
$C_0^{ar}$, $C_1^{ar}$ and $C_2^{ar}$. Therefore we do not give here
the definition of categories $C_n^{ar}$ for any $n \in \dn$.}
\end{nt}

\section{Virtual measures.} \label{vm}
Let $E = (I, F, V) \in \Ob(C_2^{ar})$.

For any $i,j \in I$ we define a one-dimensional $\dc$-vector space
of virtual measures
\begin{equation} \label{virm}
\mu(F(i) \mid F(j)) \eqdef \mathop{\mathop{\Lim_{\lto}}_{l \in
I}}_{l \le i, l \le j} \Hom\nolimits_{\sdc} (\mu(F(i)/ F(l)) \, , \,
\mu (F(j)/ F(l))) \mbox{,}
\end{equation}
where  to take  the inductive limit we need the identities which
follow from formula~(\ref{isotr1}): for $l' \le l \in I$
$$
\mu(F(i)/ F(l')) = \mu(F(i)/ F(l)) \otimes_{\sdc} \mu(F(l)/ F(l'))
$$
$$
\mu(F(j)/ F(l')) = \mu(F(j)/ F(l)) \otimes_{\sdc} \mu(F(l)/ F(l'))
\mbox{.}
$$
And the maps in this limit are given in the following way:
$$
\textstyle f \in \Hom\nolimits_{\sdc} (\mu(F(i)/ F(l)), \mu (F(j)/
F(l))) \mapsto f' \in \Hom\nolimits_{\sdc} (\mu(F(i)/ F(l')), \mu
(F(j)/ F(l')))  \mbox{,}$$
 where $f'(a \otimes c) \eqdef f(a) \otimes c$, $a$ is any
from $\mu(F(i)/ F(l))$, $c$ is any from $\mu(F(l)/ F(l'))$. These
maps are  isomorphisms.

\begin{prop} \label{pci}
For any $i,j,l \in I$ there is a canonical isomorphism
$$
\gamma \quad : \quad \mu(F(i) \mid F(j)) \otimes_{\sdc} \mu(F(j)
\mid F(l))  \arrow{e} \mu(F(i) \mid F(l))
$$
such that the following diagram of associativity is commutative for
any $i,j,l,n \in I$:
$$
\begin{diagram}
 \node{\scriptstyle \mu(F(i) \mid F(j)) \otimes_{\ssdc}
\mu(F(j) \mid F(l)) \otimes_{\ssdc} \mu(F(l) \mid F(n)) }  \arrow{e}
\arrow{s}
\node{\scriptstyle \mu(F(i) \mid F(l)) \otimes_{\ssdc} \mu(F(l) \mid F(n))} \arrow{s} \\
\node{\scriptstyle \mu(F(i) \mid F(j)) \otimes_{\ssdc} \mu(F(j) \mid
F(n))} \arrow{e} \node{\scriptstyle \mu(F(i) \mid F(n))}
\end{diagram}
$$
\end{prop}
\proof. We have a canonical map:
\begin{equation} \label{eqt}
\begin{diagram}
\node{ \Hom\nolimits_{\sdc} (\mu(F(i)/ F(l')), \mu (F(j)/ F(l')))
\otimes_{\sdc} \Hom\nolimits_{\sdc} (\mu(F(j)/ F(l')), \mu
(F(l)/F(l')))} \arrow{s} \\
\node{ \Hom\nolimits_{\sdc} (\mu(F(i)/ F(l')), \mu (F(l)/ F(l')))}
\end{diagram}
\end{equation}
which satisfies the associativity diagram. And this map commutes
with the inductive limit from the definition of $\mu (\cdot \mid
\cdot)$. We obtain the map $\gamma$ after the taking the inductive
limit in~(\ref{eqt}). The proposition is proved.
\vspace{0.3cm}

\begin{nt} \label{nci} { \em
We have the following canonical isomorphisms. For any $i, j \in I$
$$
\mu(F(i) \mid F(i)) =\dc \quad , \quad \mu(F(i) \mid F(j)) =
\mu(F(j) \mid F(i))^*  \mbox{.}
$$
Let $i \le j \in I$. Then
$$
\mu(F(i) \mid F(j)) = \mu (F(j)/F(i)) \quad , \quad \mu(F(j) \mid
F(i))= \mu (F(j)/F(i))^*  \mbox{.}
$$ }
\end{nt}

\section{Basic spaces} \label{bs}
Let $E = (I, F, V) \in \Ob(C_2^{ar})$.
Let $i \ge j \ge l \ge n
\in I$. Then $F(i) \supset F(j) \supset F(l)
\supset F(n)$, and from the definition of objects in category
$C_2^{ar}$ we have
$$
F(j)/ F(n) \sim (F(i)/ F(n)) \mathop{\times}\limits_{F(i)/F(l)}
(F(j)/F(l)) $$
$$F(i)/ F(l)  \sim (F(i)/F(n))
\mathop{\amalg}\limits_{F(j)/F(n)} (F(j)/F(l))
$$
as objects of category $C_1^{ar}$.

For any $i \ge j \ge l \in I$ we have the following admissible
triple from $C_1^{ar}$:
$$
0 \arrow{e} F(j)/ F(l)  \arrow{e,t}{\alpha_{lji}} F(i)/F(l)
\arrow{e,t}{\beta_{lji}} F(i)/F(j) \arrow{e} 0 \mbox{.}
$$

We fix some $o \in I$. Then from subsection~\ref{diim} we have that
the following spaces over the field $\dc$ are well-defined.

\begin{equation} \label{2d}
\begin{array}{c}
\s_{F(o)}(E) \eqdef \mathop{\Lim\limits_{\longleftarrow}}\limits_{p \in I}
\mathop{\Lim\limits_{\longleftarrow}}\limits_{q \le p} \s(F(p) / F(q))
\otimes_{\sdc} \mu(F(q) \mid F(o)) =  \qquad \qquad \qquad \qquad \qquad \\
\qquad \qquad \qquad \qquad \qquad =
\mathop{\Lim\limits_{\longleftarrow}}\limits_{q \in I}
\mathop{\Lim\limits_{\longleftarrow}}\limits_{p \ge q} \s(F(p) /
F(q)) \otimes_{\sdc} \mu(F(q) \mid F(o))
\end{array}
\end{equation}
with respect to the maps $(\beta_{lji})_* \otimes \Id_{\mu(F(j) \mid
F(o))}$:
$$
 \s(F(i)/F(l)) \otimes_{\sdc} \mu(F(l) \mid F(o)) \arrow{e} \s(F(i)/
F(j)) \otimes_{\sdc} \mu(F(j) \mid F(o)) \mbox{,}
$$
and the maps $\alpha^*_{lji} \otimes \Id_{\mu(F(l) \mid F(o))}$:
$$ \s(F(i)/ F(l)) \otimes_{\sdc} \mu(F(l) \mid F(o)) \arrow{e}
\s(F(j)/ F(l)) \otimes_{\sdc} \mu(F(l) \mid F(o))
  \mbox{.}$$

\begin{equation}  \label{2d'}
\begin{array}{c}
\s'_{F(o)}(E) \eqdef \mathop{\Lim\limits_{\lto}}\limits_{p \in I}
\mathop{\Lim\limits_{\lto}}\limits_{q \le p} \s'(F(p) / F(q))
\otimes_{\sdc} \mu(F(o) \mid F(q)) =  \qquad \qquad \qquad \qquad \qquad \\
\qquad \qquad \qquad \qquad \qquad =
\mathop{\Lim\limits_{\lto}}\limits_{q \in I}
\mathop{\Lim\limits_{\lto}}\limits_{p \ge q} \s'(F(p) / F(q))
\otimes_{\sdc} \mu(F(o) \mid F(q))
\end{array}
\end{equation}
with respect to the maps $\beta^*_{lji} \otimes \Id_{\mu(F(o) \mid
F(l))}$:
$$
\s'(F(i)/F(j)) \otimes_{\sdc} \mu(F(o) \mid F(j)) \arrow{e}
\s'(F(i)/F(l)) \otimes_{\sdc} \mu(F(o) \mid F(l)) \mbox{,}
$$
and the maps $(\alpha_{lji})_* \otimes \Id_{\mu(F(o) \mid F(l))}$:
$$ \s'(F(j)/ F(l)) \otimes_{\sdc} \mu(F(o) \mid F(l)) \arrow{e}
\s'(F(i)/ F(l)) \otimes_{\sdc} \mu(F(o) \mid F(l))
  \mbox{.}
$$

\vspace{0.5cm}

From definitions of these spaces we have the following canonical
isomorphism for any $o_1 \in I$:
$$
\s_{F(o)}(E)  \otimes_{\sdc} \mu(F(o) \mid F(o_1))  \arrow{e}
\s_{F(o_1)}(E)
$$
such that the following diagram is commutative for any $o_2 \in I$:
$$
\begin{diagram}
 \node{\scriptstyle   \s_{F(o)}(E)  \otimes_{\ssdc}
\mu(F(o) \mid F(o_1)) \otimes_{\ssdc} \mu(F(o_1) \mid F(o_2))}
\arrow{e} \arrow{s} \node{\scriptstyle \s_{F(o_1)}(E)
\otimes_{\ssdc} \mu(F(o_1) \mid F(o_2))}
\arrow{s} \\
\node{\scriptstyle   \s_{F(o)}(E) \otimes_{\ssdc} \mu(F(o) \mid
F(o_2))} \arrow{e} \node{\scriptstyle   \s_{F(o_2)} (E)}
\end{diagram}
$$

Dually, we have the following canonical isomorphism for any $o_1 \in
I$:
$$
\mu(F(o_1) \mid F(o)) \otimes_{\sdc}  \s'_{F(o)}(E) \arrow{e}
\s'_{F(o_1)}(E)
$$
such that the following diagram is commutative for any $o_2 \in I$:
$$
\begin{diagram}
 \node{\scriptstyle  \mu(F(o_2) \mid F(o_1)) \otimes_{\ssdc} \mu(F(o_1) \mid F(o)) \otimes_{\ssdc}
 \s'_{F(o)}(E) } \arrow{e} \arrow{s} \node{\scriptstyle \mu(F(o_2) \mid F(o_1)) \otimes_{\ssdc} \s'_{F(o_1)}(E)}
\arrow{s} \\
\node{\scriptstyle  \mu(F(o_2) \mid F(o)) \otimes_{\ssdc}
\s'_{F(o)}(E)} \arrow{e} \node{\scriptstyle   \s'_{F(o_2)} (E)}
\end{diagram}
$$

\begin{nt} \label{nntt} { \em
Let $\{B_l\}$, $l \in J$ be a projective system of vector spaces
over a field $k$ with the surjective transition maps $\phi_{l_1,l_2}
: B_{l_1} \to B_{l_2}$ for any $l_1 \ge l_2 \in J$. Let $J^0$ be a
partially ordered set with the same set as $J$ but with inverse
order then $J$. Let $\{A_l\}$, $l \in J^0$ be an inductive system of
$k$-vector spaces with the injective transition maps $\psi_{l_2,
l_1} : A_{l_2} \to A_{l_1}$. Suppose for any $l \in J$ we have a
non-degenerate $k$-linear pairing
$$
<\cdot , \cdot >_l \quad B_l \times A_l  \arrow{e}  k
$$
such that for any $l_1 \ge l_2 \in I$, $x \in B_{l_1}$, $y \in
A_{l_2}$ we have
$$
<\phi_{l_1, l_2}(x), y>_{l_2} = <x, \psi_{l_2, l_1}(y)>_{l_1}
\mbox{.}
$$
Then we have canonically non-degenerate $k$-linear pairing between
$k$-vector spaces $\mathop{\Lim\limits_{\longleftarrow}}\limits_{l
\in J} B_l$ and $\mathop{\Lim\limits_{\lto}}\limits_{l \in J^0} A_l$
which is induced by pairings $<\cdot , \cdot >_l$, $l \in J$. }
\end{nt}

For  $E = (I, F, V) \in \Ob(C_2^{ar})$, for any $i \ge j \in I$
there are non-degenerate $\dc$-linear pairings between
$\s(F(i)/F(j))$ and $\s'(F(i)/F(j))$ (see formula~(\ref{par})).
Besides, it is not difficult to see from definitions and
corresponding facts for direct an inverse images in the category
$C_0^{ar}$ that transitions maps in projective system~(\ref{2d}) are
surjective maps, and transitions maps in inductive
system~(\ref{2d'}) are injective maps.

Therefore applying remark~\ref{nntt} twice to
formulas~(\ref{2d})-(\ref{2d'}) we obtain that there is the
following non-degenerate pairing:
$$
<\cdot, \cdot >_{\s_{F(o)}(E)} \quad : \quad \s_{F(o)}(E) \times
\s'_{F(o)}(E) \arrow{e}  \dc \mbox{.}
$$

\begin{nt}
\label{remark9}
 { \em
Let $E = (I,F,V) \in \Ob(C_2^{ar})$ be a $cf$-object. Then,
according to remark~\ref{remm8}, for any $i \ge j \in I$ we have a
canonical element
$$1_{ij}  \eqdef {\bf 1_{F(i)/F(j)}}\in \mu
(F(i)/F(j))  \mbox{.}$$ (This element can be also defined via
non-degenerate pairing~(\ref{par}): $<1,1_{ij}>  = 1$.) Therefore
for any $k, l \in I$ we have a canonical element $1_{kl} \in \mu
(F(k) \mid F(l))$ such that for any $k,l,n \in I$ according to
proposition~\ref{pci}:
$$
\gamma( 1_{kl} \otimes 1_{ln})  = 1_{kn} \mbox{.}
$$
Taking into account these elements $1_{lo}$ and $1_{ol}$ for $l, o
\in I$ we can do not write $1$-dimensional $\dc$-spaces $\mu(F(l)
\mid F(o) )$ and $\mu(F(o) \mid F(l) )$ in formulas~(\ref{2d})
and~(\ref{2d'}) which define $\dc$-spaces $\s_{F(o)(E)}$ and
$\s'_{F(o)}(E)$. These $\dc$-spaces {\it do not depend} on the
choice of $o \in I$.

Similarly, let $E = (I,F, V) \in \Ob(C_2^{ar})$ be a $df$-object.
Then, according to remark~\ref{remm8}, for any $i \ge j \in I$ we
have a canonical element
$$\delta_{ij} \eqdef {\bf \sum\nolimits_{F(i)/F(j)} } \in \mu
(F(i)/F(j)) \mbox{.}$$ (This element can be also defined via
non-degenerate pairing~(\ref{par}): $<f,\delta_{ij}> =
\sum\limits_{x \in V_{ij}} f(x)$, where $f \in \s(E)$, $F(i)/F(j)=
(I_{ij}, F_{ij}, V_{ij}) \in \Ob(C_1^{ar})$.) Therefore for any $k,
l \in I$ we have a canonical element $\delta_{kl} \in \mu (F(k) \mid
F(l))$ such that for any $k,l,n \in I$ according to
proposition~\ref{pci}:
$$
\gamma( \delta_{kl} \otimes \delta_{ln})  = \delta_{kn} \mbox{.}
$$
Taking into account these elements $\delta_{lo}$ and $\delta_{ol}$
for $l, o \in I$ we can do not write $1$-dimensional $\dc$-spaces
$\mu(F(l) \mid F(o) )$ and $\mu(F(o) \mid F(l) )$ in
formulas~(\ref{2d}) and~(\ref{2d'}) which define $\dc$-spaces
$\s_{F(o)}(E)$ and $\s'_{F(o)}(E)$. These $\dc$-spaces {\it do not
depend} on the choice of $o \in I$.
 }
\end{nt}

\section{Fourier transform} \label{ft} Let $E = (I, F, V) \in \Ob(C_2^{ar})$ be a
complete object. Then we have $\check{\check{E}} =E$. We fix some $o
\in I$.

\subsection{} Let $i \ge j \in I$ be any. By subsection~\ref{subsec4.2} there
are the following $\cc$-linear maps:
$$
f \otimes \mu \in \s(F(i)/F(j)) \otimes_{\cc} \mu(F(j) \mid F(o))
\mapsto \check{f} \otimes \mu \in \s(F(i)/F(j)) \otimes_{\cc}
\mu(F(j) \mid F(o))
$$
$$
G \otimes \mu' \in \s'(F(i)/F(j)) \otimes_{\cc} \mu(F(o) \mid F(j))
\mapsto \check{G} \otimes \mu' \in \s'(F(i)/F(j)) \otimes_{\cc}
\mu(F(o) \mid F(j))  \mbox{.}
$$
These maps commute with direct and inverse images when we change $i
,j \in I$ to $i', j' \in I $, $i' \ge i$, $j' \le j$.

Therefore using formulas~(\ref{2d})-(\ref{2d'}) we have that the
following $\cc$-linear maps are well defined:
$$
f \in \s_{F(o)}(E)   \longmapsto \check{f} \in \s_{F(o)}(E)
$$
$$
G \in \s'_{F(o)}(E)   \longmapsto \check{G} \in \s'_{F(o)}(E)
\mbox{.}
$$

The square of any from these maps is the indentity map. From the
evident corresponding formulas for category $C_1^{ar}$ we have now
the following formulas:
$$
<\check{f}, G>_{\s_{F(o)}(E)} = <f, \check{G}>_{\s_{F(o)}(E)} \quad
\mbox{for any} \quad f \in \s_{F(o)}(E), \quad G \in \s'_{F(o)}(E)
\mbox{.}
$$

\subsection{} \label{2Ft} According to section~\ref{c_2-spaces} we
have the dual object $\check{E} = (I^0, F^0, \check{V}) \in
\Ob(C_2^{ar})$. For any $i \ge j \in I$, for the object $E_{j,i} =
F(i)/F(j) \in \Ob(C_1^{ar})$ we have the object $\check{E}_{j,i} =
F^0(j) / F^0(i) \in \Ob(C_1^{ar})$.

For any $l,n \in I$ we have
\begin{equation} \label{fff}
\mu(F(l) \mid F(n)) = \mu (F^0(l) \mid F^0(n))  \mbox{.}
\end{equation}
Indeed, let  $k \in I$ such that $k \le l$, $k \le n$. Then
according to proposition~\ref{pci}, remark~\ref{nci} and
formula~(\ref{isodu1}) we have
$$
\mu (F(l) \mid F(n)) = \mu(F(l) \mid F(k)) \otimes_{\cc} \mu(F(k)
\mid F(n)) = $$
$$ = \mu(F(l)/F(k))^* \otimes_{\cc} \mu(F(n)/F(k)) =
 \mu(F^0(k) / F^0(l)) \otimes_{\cc} \mu(F^0(k) / F^0(n))^* =$$
$$
= \mu(F^0(l) \mid F^0(k)) \otimes_{\cc} \mu(F^0(k) \mid F^0(n)) =
\mu(F^0(l) \mid F^0(n)) \mbox{.}
$$

Therefore we have the following maps for any $i \ge j \in I$
(Fourier trasnsforms):
\begin{equation}
\begin{array}{cc}
\F  \otimes  \mu(F(i) \mid F(o)) \quad :
\\
  \s(F(i)/ F(j)) \otimes_{\cc} \mu (F(j) \mid F(o)) \arrow{e}
\s(F^0(j)/ F^0(i)) \otimes_{\cc} \mu(F^0(i) \mid F^0(o))
\end{array}
\end{equation}
\begin{equation}
\begin{array}{c}
\F \: \otimes \: \mu(F(o) \mid F(i)) \quad : \\  \s'(F(i)/ F(j)) \oc
\mu (F(o) \mid F(j)) \arrow{e} \s'(F^0(j)/ F^0(i)) \otimes_{\cc}
\mu(F^0(o) \mid F^0(i))
\end{array}
\end{equation}

Now we use remark~\ref{remm8}, more exactly: the analog of
proposition~\ref{prpr}  for category $C_1^{ar}$, which connect
Fourier transform with direct and inverse images. We use it when we
change $i ,j \in I$ to $i', j' \in I $, $i' \ge i$, $j' \le j$, and
take the limits according to formulas~(\ref{2d})-(\ref{2d'}). We
obtain that the following $\cc$-linear maps ({\it two-dimensional
Fourier transforms}) are well defined:
\begin{equation} \label{ae1}
\F \quad : \quad \s_{F(o)}(E) \arrow{e} \s_{F^0(o)}(\check{E})
\end{equation}
\begin{equation} \label{ae2}
\F \quad : \quad \s'_{F(o)}(E) \arrow{e} \s'_{F^0(o)}(\check{E})
\end{equation}

\begin{nt} {\em
Since for any $i \ge j \in I$ the spaces $\s(F(i)/F(j))$ and
$\s'(F(i)/F(j))$ are topological $\dc$-vector spaces, and
$\dim_{\dc} \mu (F(j) \mid F(o)) =1 $, we have that the spaces
$\s(E)$ and $\s'(E)$ are also topological $\dc$-vector spaces, which
have the topologies of projective and inductive limits
correspondingly (see~\cite{RR}) according to
formulas~(\ref{2d})-(\ref{2d'}). }
\end{nt}

\begin{prop}
The following statements are satisfied.
\begin{enumerate}
\item The maps $\F$ given by formulas~(\ref{ae1})-(\ref{ae2}) are
isomorphisms of topological $\cc$-vector spaces.
\item For any $f \in \s_{F(o)}(E)$ and $G \in \s'_{F(o)}(E)$
$$
\F \circ \F (f) = \check{f}  \mbox{,} \quad  \F \circ \F (G) =
\check{G} \mbox{.}
$$
\item For any $f \in \s_{F(o)}(E)$ and $G \in \s'_{F^0(o)}(\check{E})$
$$
<\F(f), G>_{\s_{F^0(o)}(\check{E})} = <f, \F(G)>_{\s_{F(o)}(E)}
\mbox{.}
$$
\end{enumerate}
\end{prop}
\proof follows from the construction of two-dimensional Fourier
transform given in this section and the properties of  Fourier
transform for category $C_1^{ar}$(see subsection~\ref{subsec4.2} and
proposition~\ref{fourtra1}) applied to functions and distributions
on objects $F(i)/F(j) \in \Ob(C_1^{ar})$ for any $i \ge j \in I$.
The proposition is proved.
 \vspace{0.3cm}

\section{Central extension and its representations} \label{ce} In
this section we construct a central extension of some subgroup of
group of automorphisms of an object of the category $C_2^{ar}$ and
consider representations of this group on the spaces of functions
and distributions of the given object.
\subsection{} \label{ss}
 Let an object $E_1= (I_1, F_1, V) \ic$ dominates an object
 $E_2= (I_2, F_2, V) \ic $ (see section~\ref{c_2-spaces}). We fix some $o \in
I_2$. Then we obtain the following canonical isomorphisms:
$$ \s_{F(o)}(E_1) = \s_{F(o)}(E_2) \; \mbox{,} \qquad \qquad
  \s'_{F(o)}(E_1) = \s'_{F(o)}(E_2)  \mbox{.}
$$
Indeed, from formulas~(\ref{2d})-(\ref{2d'}) we obtain that
inductive and projective limits used in the definitions of these
spaces will be the same on  sets of indices depending on $I_1$ or
$I_2$. Moreover, by the same reasons, the two-dimensional Fourier
transform coincides for the spaces depending on $E_1$ or $E_2$.
(Compare it with remark~\ref{rem9} and reasonings after
definition~\ref{defff}.)

\subsection{} \label{ssec} Let $E= (I, F, V) ) \ic $.
We define the following subgroup of the group $\Aut_{C_2^{ar}}(E)$.

\begin{defin} \label{grdef}
A set $\Aut_{C_2^{ar}}(E)'$ consists of elements $g \in
\Aut_{C_2^{ar}}(E)$ such that
\begin{enumerate}
\item \label{us1} for any $i \in I$ there is $j \in I$
such that $g F(i)= F(j)$;
\item \label{us2}  for any $i \ge j \in I$ let $g(F(i))=F(p)$,
$g(F(j))=F(q)$, then $p \ge q \in I$;
\item the element $g^{-1}\in
\Aut_{C_2^{ar}}(E)$ satisfies conditions which are analogous to
conditions~\ref{us1}-\ref{us2}  of this definition.
\end{enumerate}
\end{defin}

It is clear that the set  $\Aut\nolimits_{C_2^{ar}}(E)'$ is a subgroup of the group
$\Aut\nolimits_{C_2^{ar}}(E)$.

Now
we construct an object $\tilde{E} \ic$, which dominates the object $E \ic$.

At first, for any fixed $i \ge j \in I$ we consider the natural
object $E_{j,i}= (I_{j,i}, F_{j,i}, F(i)/F(j)) \in \Ob(C_1^{ar})$
which appears from $C_2^{ar}$-structure on $E$. We define an object
$\widetilde{E_{j,i}} \in \Ob(C_1^{ar})$, which dominates the object
$E_{j,i}$, in the following way. Let a filtered abelian group
$\widetilde{E_{j,i}} =(G, Id, F(i)/F(j))$, where $G$ is a set of the
following subgroups of the group $F(i)/F(j)$:
 $g F_{l, k}(n)$, where
$g \in \Aut\nolimits_{C_2^{ar}}(E)'$, $k \ge l \in I$, $n \in  I_{l,k} $, $g(F(k))= F(i)$, $g(F(l))=F(j)$, $F_{l,k}(n) \subset F(k)/F(l)$.
The set $G$ is partially ordered by inclusions of subgroups. The function $Id$ maps an element from $G$ to the corresponding subgroup of $F(i)/F(j)$.
From definition~\ref{morc2ar} it follows that for any such subgroup $g F_{l, k}(n)$ there are $p \ge q \in I_{j,i}$
such that
$$
F_{j,i}(p) \supset g F_{l, k}(n)   \supset  F_{j,i}(q)  \mbox{,}
$$
and the subgroup $g F_{l, k}(n)$ is a closed subgroup in the Lie
group $F_{j,i}(p)/F_{j,i}(q)$. Therefore, according to
remark~\ref{Cartan} we obtain that the filtered abelian group
$\widetilde{E_{j,i}}$ is well-defined as an object of $C_1^{ar}$
such that the object $\widetilde{E_{j,i}}$ dominates the
object~$E_{j,i}$.

Immediately from these constructions we have that $\tilde{E} = (I, F, V)$, where on the group $F(i)/F(j)$ ($i \ge j \in I$) we put a $C_1^{ar}$-structure
$\widetilde{E_{j,i}}$, is well-defined as an object from $C_2^{ar}$. Besides, the object   $\tilde{E}$ dominates the object
$E$ in the category $C_2^{ar}$. Therefore, $\Aut_{C_2^{ar}}(\tilde{E})= \Aut_{C_2^{ar}}(E)$.
By definition we have that
\begin{equation} \label{izom}
\Aut\nolimits_{C_2^{ar}}(\tilde{E})'= \Aut\nolimits_{C_2^{ar}}(E)'
\end{equation}
as subgroups of $\Aut_{C_2^{ar}}(E)$.

\subsection{} Let $E=(I, F, V) \ic$. Then for any $i_1 , i_2 \in I$ we have the virtual
measures space $\mu(F(i_1) \mid F(i_2))$.  For any $g \in
\Aut_{C_2^{ar}}(E)'$, for any $i \in I$ we have from the definition
of the group $\Aut_{C_2^{ar}}(E)'$ that $g F(i) = F(p)$ for some $p
\in {I}$. Therefore for any $g , h \in \Aut_{C_2^{ar}}(E)'$, for any
$i_1 , i_2 \in I$ we have the well-defined one-dimensional
$\cc$-vector space  $  \mu(g F(i_1) \mid h F(i_2))$.

Using isomorphism~(\ref{izom}) and subsection~\ref{ss}, for any $g
\in \Aut_{C_2^{ar}}(E)'$, for any $p \ge q \in I$ we have the
following isomorphism:
\begin{equation} \label{form1}
r_g \quad : \quad    \s({F}(p) / {F}(q)) \arrow{e} \s(g{F}(p) /
g{F}(q)) \mbox{,}
\end{equation}
where $r_g(f) (v) \eqdef f (g^{-1} v)$ for any $v \in g F(p) /
g{F}(q)$ and any $f$ from $\s({F}(p) / {F}(q))$. (We used that it is
enough to consider $\tilde{E}$ instead of $E$, and $\s(F(i)/F(j))$
is a $\dc$-subalgebra of the space $\f(F(i)/F(j))$ of all
$\dc$-valued functions on $V$, where $i \ge j \in I$ are any
elements.)

Dually to formula~(\ref{form1}), for any $g \in
\Aut_{C_2^{ar}}(E)'$, for any $p \ge q \in {I}$ we have the
following isomorphism:
\begin{equation}  \label{form3'}
r'_g \quad : \quad    \s'({F}(p) / {F}(q)) \arrow{e} \s'(g{F}(p) /
g{F}(q)) \mbox{.}
\end{equation}
(In other words, the map $r'_g$ is the conjugate map to the map
$r_{g^{-1}}$ with respect to the pairing~(\ref{par})).

Using remark~\ref{meass}, we obtain that the map $r_g'$ induces the
following isomorphism of one-dimensional $\dc$-vector spaces:
$$
n_g \quad : \quad \mu ({F}(p) /{F}(q)) \arrow{e} \mu (g{F}(p)
/g{F}(q)) \mbox{.}
$$

For any $g \in \Aut_{C_2^{ar}}(E)'$, for any $p, q, s \in {I}$,  $s
\le p$, $s \le q$ we have the following $\cc$-isomorphism $m_g$:
$$
 \Hom\nolimits_{\cc} (\mu({F}(p)/ {F}(s))  ,  \mu({F}(q)/ {F}(s)))
\arrow{e} \Hom\nolimits_{\cc} (\mu(g{F}(p)/ g{F}(s))  , \mu(g{F}(q)/
g{F}(s))) \mbox{,}
$$
where for any $f \in \Hom\nolimits_{\cc} (\mu({F}(p)/ {F}(s)) ,
\mu({F}(q)/ {F}(s)))$:
$$
m_g (f) \eqdef n_g \circ f \circ n_{g^{-1}}  \mbox{.}
 $$

Now we apply isomorphisms $m_g$ to the inductive limit  in
formula~(\ref{virm}), then we obtain the following $\cc$-linear
isomorphism for any $p, q \in {I}$:
$$
l_g \quad : \quad \mu ({F}(p) \mid {F}(q)) \arrow{e} \mu(g{F}(p)
\mid g{F}(q))  \mbox{.}
$$

We have for any $g_1, g_2 \in \Aut_{C_2^{ar}}(E)'$ that $l_{g_1 g_2}
= l_{g_1} l_{g_2}$. Besides, for any $a \in \mu({F}(p) \mid
{F}(q))$, for any $b \in \mu({F}(q) \mid {F}(s))$, for any $p,q,s
\in {I}$ we have that
$$\gamma (l_g (a) \otimes l_g(b) )= l_g \gamma (a \otimes
b)$$ (see proposition~\ref{pci} for definition of the
isomorphism~$\gamma$).

Let $p,q \in {I}$. Then for any $\mu \in \mu ({F}(p) \mid {F}(q) )$,
$\mu \ne 0$ we define canonically $\mu^{-1} \in \mu ({F}(q) \mid
{F}(p) )$
 such that $\mu \otimes \mu^{-1} =1$ with respect to the following
 canonical isomorphism:
 $$
\mu ({F}(p) \mid {F}(q) )\otimes_{\cc} \mu ({F}(q) \mid {F}(p) ) =
\cc \mbox{.}
 $$

 Let $E= (I,F,V) \ic$. We fix some $o \in I$.
Then there is the following {\em central extension of groups}:
\begin{equation} \label{cextm}
1 \arrow{e} {\cc}^* \arrow{e}
\widehat{\Aut\nolimits_{C_2^{ar}}(E)'}_{F(o)} \arrow{e,t}{\Lambda}
\Aut\nolimits_{C_2^{ar}}(E)' \arrow{e} 1 \mbox{,}
\end{equation}
where
$$
\widehat{\Aut\nolimits_{C_2^{ar}}(E)'}_{F(o)} \eqdef \{ (g, \mu) \:
\; : \; \:
 g \in  \Aut\nolimits_{C_2^{ar}}(E)', \; \mu \in \mu(F(o) \mid g F(o)),
\; \mu \ne 0 \} \mbox{.}
$$
Here $\Lambda((g, \mu)) = g$. The operations $(g_1, \mu_1) \cdot
(g_2, \mu_2) = (g_1 g_2, \gamma(\mu_1 \otimes l_{g_1}(\mu_2)))$ and
$ (g, \mu)^{-1} = (g^{-1}, l_{g^{-1}}(\mu^{-1}))$ define the
structure of a group on the set
$\widehat{\Aut\nolimits_{C_2^{ar}}(E)'}_{F(o)}$. (The unit element
of this group is $(e,1)$, where $e$ is the unit element of the group
$\Aut\nolimits_{C_2^{ar}}(E)'$).

\begin{nt} {\em
For any $o_1 \in I$ there is a canonical isomorphism
$$
\alpha_{o, o_1} \; : \;
\widehat{\Aut\nolimits_{C_2^{ar}}(E)'}_{F(o)} \lto
\widehat{\Aut\nolimits_{C_2^{ar}}(E)'}_{F(o_1)} \mbox{.}
$$
Indeed, we fix any $\nu \in \mu(F(o_1) \mid  F(o))$, $\nu \ne 0$.
Then
$$
\alpha_{o, o_1} ((g, \mu)) \eqdef (g, \gamma( \gamma(\nu \otimes
\mu) \otimes l_g(\nu^{-1} ))) \mbox{.}
$$
The map $\alpha_{o, o_1}$ does not depend on the choice of $\nu \in
\mu(F(o_1) \mid  F(o))$, $\nu \ne 0$.}
\end{nt}

\subsection{} Let $E=(I,F,V) \ic$. In subsection~\ref{ssec} we constructed the
object $ \tilde{E} \ic$ which dominates the object $E$.

For any $g \in \Aut_{C_2^{ar}}(E)'$, for any  $p \ge q \in {I}$ we
have isomorphism~(\ref{form1}):
$$
r_g \quad : \quad    \s({F}(p) / {F}(q)) \arrow{e} \s(g{F}(p) /
g{F}(q)) \mbox{.}
$$

We fix some $o \in I$. Using formula~(\ref{form1}) we construct the
map $R_{\tilde{g}}$ for any $\tilde{g}= (g, \mu) \in
\widehat{\Aut\nolimits_{C_2^{ar}}(E)'}_{F(o)}$ for any $p \ge q \in
{I}$:
\begin{equation} \label{rf}
R_{\tilde{g}} \quad : \quad  \s({F}(p) / {F}(q)) \otimes_{\cc}
\mu({F}(q) \mid F(o)) \arrow{e} \s(g {F}(p) / g {F}(q))
\otimes_{\cc} \mu( g {F}(q) \mid F(o))
\end{equation}
as composition of the map $r_g \otimes l_g$ with multiplication by
$\mu^{-1} \in \mu(g F(o) \mid F(o))$.

Applying this formula to formula~(\ref{2d}) we obtain the map for
any $\tilde{g} \in \widehat{\Aut\nolimits_{C_2^{ar}}(E)'}_{F(o)}$:
$$
R_{\tilde{g}} \quad : \quad \s_{F(o)}(\tilde{E}) \arrow{e}
\s_{F(o)}(\tilde{E})
$$
such that $R_{\tilde{g}} R_{\tilde{h}} = R_{\tilde{g} \tilde{h}}$.

Therefore using subsection~\ref{ss} we obtain {\em the
representation of the group} $
\widehat{\Aut\nolimits_{C_2^{ar}}(E)'}_{F(o)}$ by maps
$R_{\tilde{g}}$ on the $\cc$-space $\s_{F(o)}(E)$.

Let $o \in I$. Dually to formula~(\ref{rf}) and using the
isomorphisms $r_g'$ given by formula~(\ref{form3'}),
 we construct the map $R'_{\tilde{g}}$ for any $\tilde{g}= (g, \mu)
\in \widehat{\Aut\nolimits_{C_2^{ar}}(E)'}_{F(o)}$ for any $p \ge q
\in {I}$:
\begin{equation} \label{rf'}
R'_{\tilde{g}} \quad : \quad  \s'({F}(p) / {F}(q)) \otimes_{\cc}
\mu(F(o) \mid {F}(q) ) \arrow{e} \s'(g {F}(p) / g {F}(q))
\otimes_{\cc} \mu( F(o) \mid  g {F}(q))
\end{equation}
as composition of the map $r_g' \otimes l_g$ with multiplication by
$\mu \in \mu( F(o) \mid gF(o))$.

Applying formula~(\ref{rf'}) to formula~(\ref{2d'}), and using
subsection~\ref{ss}, we obtain {\em the representation of the group}
$ \widehat{\Aut\nolimits_{C_2^{ar}}(E)'}_{F(o)}$  on the $\cc$-space
$\s'_{F(o)}(E)$ by maps $R'_{\tilde{g}}$.

By constructions, we have the following formula for any $H \in
\s'_{F(o)}(E)$, any $f \in \s_{F(o)}(E)$, any $\tilde{g} \in
\widehat{\Aut\nolimits_{C_2^{ar}}(E)'}_{F(o)} $:
\begin{equation} \label{mf3}
< R'_{\tilde{g}}(H)  \:,\:  R_{\tilde{g}}(f)  >_{\s_{F(o)}(E)} = < H
 \:,\: f>_{\s_{F(o)}(E)} \mbox{.}
\end{equation}

\subsection{} \label{sec5.5.5} Let $E = (I,F,V) \ic$ be a complete
object. Then there is the dual object $\check{E} = (I^0, F^0,
\check{V})$, and $\check{\check{E}}=E$. For any $g \in
\Aut_{C_2^{ar}}(E)$ we have canonically $\check{g} \in
\Aut_{C_2^{ar}}(\check{E})$. And we have the following isomorphism
of groups:
$$
g \in \Aut\nolimits_{C_2^{ar}}(E) \longmapsto \check{g}^{-1} \in
\Aut\nolimits_{C_2^{ar}}(\check{E}) \mbox{,}
$$
which induces the following isomorphism of groups:
$$
\Aut\nolimits_{C_2^{ar}}(E)' \lto
\Aut\nolimits_{C_2^{ar}}(\check{E})' \mbox{.}
$$

Hence, for any $o \in I$ we have the following isomorphism of
groups:
\begin{equation} \label{eee2}
(g, \mu) \in \widehat{\Aut\nolimits_{C_2^{ar}}(E)'}_{F(o)} \longmapsto
(\check{g}^{-1}, \mu)  \in
\widehat{\Aut\nolimits_{C_2^{ar}}(\check{E})'}_{F^{0}(o)} \mbox{,}
\end{equation}
where we use isomorphism~(\ref{fff}) to obtain that for any $g \in
\Aut_{C_2^{ar}}(E)$ we have canonically
$$
\mu(F(o) \mid g F(o)) = \mu(F^0(o) \mid \check{g}^{-1} F^0(o) )
\mbox{.}
$$

Using isomorphism~(\ref{eee2}) we obtain representations of the
group $\widehat{\Aut\nolimits_{C_2^{ar}}(E)'}_{F(o)}$ on $\cc$-spaces
$\s_{F^0(o)}(\check{E})$ and $\s'_{F^0(o)}(\check{E})$.

We have the following proposition.
\begin{prop}
Let $E = (I,F,V) \ic$ be a complete object.  Then the Fourier
transform $\F$ gives an isomorphism between representations of the
group $\widehat{\Aut\nolimits_{C_2^{ar}}(E)'}_{F(o)}$ on $\cc$-spaces
$\s_{F(o)}(E)$ and $\s_{F^0(o)} (\check{E})$ and on $\cc$-spaces
$\s'_{F(o)}(E)$ and $\s'_{F^0(o)} (\check{E})$.
\end{prop}
\proof. Using the definition of two-dimensional Fourier transform
$\F$ from section~\ref{2Ft}, and using
formulas~(\ref{2d})-(\ref{2d'}), we reduce the statements of this
proposition to corresponding statements about isomorphisms of
objects of category $C_1^{ar}$ and one-dimensional Fourier
transforms between them. These last statements follow from the
definition and properties of one-dimensional Fourier transform, or
can be again reduced to the case of objects from category $C_0^{ar}$
by formulas~(\ref{s})-(\ref{s'}) using the definition of
one-dimensional Fourier transforms from section~\ref{subsec4.2}. The
proposition is proved.

\begin{nt} \label{appl}{\em
We used in this section the subgroup $\Aut_{C_2^{ar}}(E)'$ instead
of the whole group $\Aut_{C_2^{ar}}(E)$ for $E \ic$.  In a lot of
important examples the groups which are considered are subgroups of
the group $\Aut_{C_2^{ar}}(E)'$ if we have the enough fine
filtration on $E$. For example, $E$ can be constructed from
two-dimensional local field $K$ with the finite last residue field,
or  from $K \simeq \dr((t))$, or from $E=K \simeq \dc((t))$, and the
group $K^* \subset \Aut_{C_2^{ar}}(E)'$ (see examples~\ref{exam1}
and~\ref{eexx}). In another example, $E$ can be the arithmetic
adelic ring $\da_X^{ar}$, where $X$ is an "arithmetic surface", and
the group of invertible elements (idelic group) $(\da_X^{ar})^*
\subset \Aut_{C_2^{ar}}(E)'$, see example~\ref{arsur}. }
\end{nt}

\section{Direct and inverse images} \label{secdi} In this section
we construct direct and inverse images of spaces of functions and
distributions on objects of category $C_2^{ar}$.
\subsection{} Let
$$
0\arrow{e} E_1\arrow{e,t}{\alpha} E_2\arrow{e,t}{\beta}E_3\arrow{e}0
$$
be an admissible triple from category $C_2^{ar}$, where
$E_i=(I_i,F_i,V_i)$, \linebreak $1 \le i \le 3$. By definition,
there are  order-preserving  functions:
\begin{itemize}
\item[(i)] $ \gamma \; \colon \; I_2\arrow{e} I_3\quad\mbox{such
that}\quad \beta(F_2(i))=F_3(\gamma(i))\quad\mbox{for any} ~i\in I_2
\mbox{,} $
\item[(ii)] $ \veps \; \colon \; I_2\arrow{e} I_1\quad\mbox{such that}\quad
F_2(i)\cap V_1=F_1(\veps(i))\quad\mbox{for any} ~i\in I_2. $
\end{itemize}

\begin{prop} \label{pr26}
We suppose that $E_1$ is a $c$-object. Let $o\in I_2$. Then there is
the direct image
$$
\beta_* \quad \colon \quad
\s_{F_{2}(o)}(E_2)\otimes_{\cc}\mu(F_1(\veps(o))\mid V_1)\arrow{e}
\s_{F_{3}(\gamma(o))}(E_3).
$$
\end{prop}
\proof. For any $i,j\in I_2$ we have canonically
\begin{equation} \label{*}
\mu(F_2(i)\mid F_2(j))=\mu(F_1(\veps(i))\mid
F_1(\veps(j)))\otimes_{\cc} \mu(F_3(\gamma(i))\mid F_3(\gamma(j)))
\mbox{.}
\end{equation}

Let $k\in I_2$ such that $F_1(\veps(k))=V_1$. Let $i\ge j\in I_2$ be
any such that $i\ge k\ge j$. Then we have an admissible triple from
category $C_1$:
$$
0\arrow{e} V_1/F_1(\veps(j))\arrow{e}
F_2(i)/F_2(j)\arrow{e,t}{\beta_{ij}}
F_3(\gamma(i))/F_3(\gamma(j))\arrow{e} 0.
$$

We have a well-defined map:
$$
(\beta_{ij})_*   \quad \colon  \quad
\s(F_2(i)/F_2(j))\otimes_{\cc}\mu(F_1(\veps(j))\mid V_1)\arrow{e}
\s(F_3(\gamma(i))/F_3(\gamma(j))) \mbox{,}
$$
where we use that $\mu(F_1(\veps(j))\mid V_1)=\mu(V_1 /
F_1(\veps(j)))$.

Therefore we have a well-defined map
$$
(\beta_{ij})_* \otimes \Id\nolimits_{\mu(F_3(\gamma(j))\mid
F_3(\gamma(o)))} \quad \colon
$$
$$
 \s(F_2(i)/F_2(j)) \otimes_{\cc} \mu(F_1(\veps(j)) \mid
V_1)\otimes_{\cc} \mu(F_3(\gamma(j))\mid F_3(\gamma(o)))\arrow{e}
\qquad \qquad \qquad \qquad
$$
$$
\qquad \qquad \qquad \qquad \qquad \qquad \arrow{e}
\s(F_3(\gamma(i))\mid
F_3(\gamma(j)))\otimes_{\cc}\mu(F_3(\gamma(j))\mid F_3(\gamma(o)))
\mbox{.}
$$
From (\ref{*}) we have that
$$
\mu(F_3(\gamma(j))\mid F_3(\gamma(o)) = \mu(F_2(j) \mid
 F_2(o)) \otimes_{\cc}  \mu(F_1(\veps(o))\mid F_1(\veps(j)))
 \mbox{.}
$$
Therefore
$$
\mu(F_1(\veps(j))\mid V_1)\otimes_{\cc}\mu(F_3(\gamma(j))\mid
F_3(\gamma(o)))=\mu(F_2(j) \mid
F_2(o))\otimes_{\cc}\mu(F_1(\veps(o))\mid V_1) \mbox{.}
$$

Thus we have a well-defined map
$$
(\beta_{ij})_*\otimes \Id\nolimits_{\mu(F_3(\gamma(j))\mid
F_3(\gamma(o)))} \quad \colon
$$
$$
\s(F_2(i)/F_2(j))\otimes_{\cc}\mu(F_2(j) \mid
F_2(o))\otimes_{\cc}\mu(F_1(\veps(o))\mid V_1)\arrow{e}  \qquad
\qquad \qquad \qquad \qquad
$$
$$
\qquad \qquad \qquad \qquad \qquad \qquad \arrow{e}
\s(F_3(\gamma(i)) /
F_3(\gamma(j)))\otimes_{\cc}\mu(F_3(\gamma(j))\mid F_3(\gamma(o)))
\mbox{.}
$$
Now we take the projective limits which are used to construct spaces
$\s_{F_2(o)}(E_2)$ and $\s_{F_3(\gamma(o))}(E_3)$. Then we obtain a
well-defined map $\beta_*$ from the maps $(\beta_{ij})_*\otimes
\Id\nolimits_{\mu(F_3(\gamma(j))\mid F_3(\gamma(o)))}$. The
proposition is proved.

\begin{nt} \em{
Dually to the maps above, we obtain a well-defined map
$$
\beta^{*} \quad \colon \quad
\s'_{F_3(\gamma(o))}(E_3)\otimes_{\cc}\mu(F_3(\veps(o))\mid
V_1)\arrow{e} \s'_{F_2(o)}(E_2)
$$
such that the maps $\beta_*$ and $\beta^*$ are conjugate maps with
respect to the pairings \linebreak
$<\cdot,\cdot>_{\s_{F_2(o)}(E_2)}$ and
$<\cdot,\cdot>_{\s_{F_3(\gamma(o))}(E_3)}$.}
\end{nt}

\subsection{} Let
$$
0\arrow{e} E_1\arrow{e,t}{\alpha} E_2\arrow{e,t}{\beta}E_3\arrow{e}0
$$
be an admissible triple from category $C_2^{ar}$, where
$E_i=(I_i,F_i,V_i)$, \linebreak$1 \le i \le 3$. By definition, there
are order-preserving  functions:
\begin{itemize}
\item[(i)]
$ \gamma \; \colon \; I_2 \arrow{e} I_3 \quad\mbox{such that}\quad
\beta(F_2(i))=F_3(\gamma(i))\quad\mbox{for any} \quad i\in I_2
\mbox{,} $
\item[(ii)]
$ \veps \; \colon \; I_2 \arrow{e} I_1 \quad\mbox{such that}\quad
F_2(i)\cap V_1=F_1(\veps(i))\quad\mbox{for any}\quad i\in I_2
\mbox{.} $
\end{itemize}

\begin{prop} \label{pr27}
We suppose that $E_3$ is a $d$-object. Let $o \in I_2$. Then there
is the inverse image
$$
\alpha^* \quad \colon \quad
\s_{F_{2}(o)}(E_2)\otimes_{\cc}\mu(F_3(\gamma(o))\mid
\{0\})\arrow{e} \s_{F_{1}(\veps(o))}(E_1) \mbox{,}
$$
where $\{0\}$ is  the zero subgroup of $V_3$.
\end{prop}
\proof. For any $i,j\in I_2$ we have canonically
\begin{equation} \label{**}
\mu(F_2(i)\mid F_2(j))=\mu(F_1(\veps(i))\mid
F_1(\veps(j)))\otimes_{\cc} \mu(F_3(\gamma(i))\mid F_3(\gamma(j)))
\mbox{.}
\end{equation}

Since $\beta$ is an admissible epimorphism, and $E_3$ is a
$d$-object, there is $k\in I_2$ such that $F_3(\gamma(k))=\{0\}$.
Let $i\ge j\in I_2$ be any such that $i\ge k\ge j$. Then we have an
admissible triple from category $C_1^{ar}$:
\begin{equation} \label{***}
0\arrow{e} F_1(\veps(i))/F_1(\veps(j))
\arrow{e,t}{\alpha_{ij}}F_2(i)/F_2(j)\arrow{e}
F_3(\gamma(i))\arrow{e}0 \mbox{.}
\end{equation}

We have a well-defined map
$$
\alpha_{ij}^* \quad \colon \quad \s(F_2(i)/F_2(j))\arrow{e}
\s(F_1(\veps(i))/F_1(\veps(j))) \mbox{.}
$$

Therefore we have a well-defined map
$$
\alpha_{ij}^* \otimes \Id\nolimits_{\mu(F_1(\veps(j))\mid
F_1(\veps(o)))} \quad \colon
$$
$$
\s(F_2(i)/F_2(j))\otimes_{\cc}\mu(F_1(\veps(j))\mid
F_1(\veps(o)))\arrow{e} \qquad \qquad \qquad \qquad \qquad \qquad $$
$$
\qquad \qquad \qquad \qquad \qquad \qquad \arrow{e}
 \s(F_1(\veps(i))\mid
F_1(\veps(j)))\otimes_{\cc}\mu(F_1(\veps(j))\mid
F_1(\veps(o)))\mbox{.}
$$
From (\ref{**}) and (\ref{***}) we have that
$$
\mu(F_2(o) \mid F_2(j))=\mu(F_1(\veps(o))\mid
F_1(\veps(j)))\otimes_{\cc}\mu(F_3(\gamma(o))\mid \{0\}) \mbox{.}
$$
Hence
$$
\mu(F_1(\veps(j))\mid F_1(\veps(o)))=\mu(F_2(j) \mid
F_2(o))\otimes_{\cc}\mu(F_3(\gamma(o))\mid \{0\}) \mbox{.}
$$
Thus we have a well-defined map
$$
\alpha_{ij}^*\otimes \Id\nolimits_{\mu(F_1(\veps(j))\mid
F_1(\veps(o)))} \quad \colon
$$
$$
\s(F_2(i)/F_2(j))\otimes_{\cc}\mu(F_2(j)/F_2(o))\otimes_{\cc}\mu(F_3(\gamma(o))\mid
\{0\})\arrow{e} \qquad \qquad \qquad \qquad \qquad \qquad
$$
$$
\qquad \qquad \qquad \qquad \qquad \qquad
 \arrow{e} \s(F_1(\veps(i))\mid
F_1(\veps(j)))\otimes_{\cc}\mu(F_1(\veps(j))\mid
F_1(\veps(o)))\mbox{.}
$$
Now we take the projective limits (with respect to $i\ge k\ge j$)
which we used to construct spaces $\s_{F_2(o)}(E_2)$ and
$\s_{F_1(\veps(o))}(E_1)$. Then we obtain a well-defined map
$\alpha^*$ from the maps $\alpha_{ij}^*\otimes
\Id\nolimits_{\mu(F_1(\veps(j))\mid F_1(\veps(o)))}$. The
proposition is proved.

\begin{nt} {\em
Dually to the maps above, we obtain a well-defined map
$$
\alpha_{*} \quad \colon \quad
\s'_{F_1(\veps(o))}(E_1)\otimes_{\cc}\mu(F_3(\gamma(o))\mid
\{0\})\arrow{e} \s'_{F_2(o)}(E_2)
$$
such that the maps $\alpha^*$ and $\alpha_*$ are conjugate maps with
respect to the pairings \linebreak
$<\cdot,\cdot>_{\s_{F_2(o)}(E_2)}$ and
$<\cdot,\cdot>_{\s_{F_1(\veps(o))}(E_1)}$. }
\end{nt}

\subsection{} Let
$$
0\arrow{e} E_1\arrow{e,t}{\alpha} E_2\arrow{e,t}{\beta}E_3\arrow{e}0
$$
be an admissible triple from category $C_2^{ar}$, where
$E_i=(I_i,F_i,V_i)$,  $1 \le i \le 3$. By definition, there are
order-preserving  functions:
\begin{itemize}
\item[(i)]
$ \gamma \; \colon \; I_2\arrow{e} I_3\quad\mbox{such that}\quad
\beta(F_2(i))=F_3(\gamma(i))\quad\mbox{for any}\quad i\in I_2
\mbox{,} $
\item[(ii)]
$ \veps \; \colon \; I_2\arrow{e} I_1\quad\mbox{such that}\quad
F_2(i)\cap V_1=F_1(\veps(i))\quad\mbox{for any}\quad i\in I_2
\mbox{.} $
\end{itemize}

\begin{prop} \label{pr28}
We suppose that $E_1$ is a $cf$-object. Let $o \in I_2$. Then there
is the inverse image
$$
\beta^* \quad \colon \quad \s_{F_{3}(\gamma(o))}(E_3)\arrow{e}
\s_{F_{2}(o)}(E_2) \mbox{.}
$$
\end{prop}
\proof. For any $i,j\in I_2$ we have canonically
\begin{equation} \label{sharp}
\mu(F_2(i)\mid F_2(j))=\mu(F_1(\veps(i))\mid
F_1(\veps(j)))\otimes_{\cc} \mu(F_3(\gamma(i))\mid F_3(\gamma(j)))
\mbox{.}
\end{equation}

Since $E_1$ is a $cf$-object, for any $i\ge j\in I_2$ we have a
canonical element $1_{ij}\in \mu(F_1(\veps(i)) / F_1(\veps(j)))$
constructed in remark~\ref{remark9}.

Therefore we have a canonical element $1_{ij}\in \mu(F_1(\veps(i))
\mid F_1(\veps(j)))$ for any \linebreak $i,j\in I_2$ such that
$1_{ij}\otimes 1_{jk}=1_{ik}$ for any $i,j,k\in I_2$. Thus,  we have
canonically that $\mu(F_1(\veps(i))\mid F_1(\veps(j)))=\cc$ for any
$i,j\in I_2$. (See also remark~\ref{remark9}).

Thus from formula (\ref{sharp}) we have for any $i,j\in I_2$
\begin{equation} \label{sharpsharp}
\mu(F_2(i)\mid F_2(j))=\mu(F_3(\gamma(i))\mid F_3(\gamma(j)))
\mbox{.}
\end{equation}

For any $i\ge j\in I_2$ we have an admissible triple from category
$C_1^{ar}$:
$$
0\arrow{e}
\frac{F_1(\veps(i))}{F_1(\veps(j))}\arrow{e}\frac{F_2(i)}{F_2(j)}
\arrow{e,t}{\beta_{ij}}
\frac{F_3(\gamma(i))}{F_3(\gamma(j))}\arrow{e}0 \mbox{.}
$$

We have a well-defined map:
$$
\beta_{ij}^* \quad \colon \quad \s(F_3(\gamma(i))/F_3(\gamma(j)))
\arrow{e} \s(F_2(i)/F_2(j)) \mbox{,}
$$
because $\frac{F_1(\veps(i))}{F_1(\veps(j))}$ is a compact object
from category $C_1^{ar}$.

From formula (\ref{sharpsharp}) we have that
$$
\mu(F_2(j) \mid F_2(o)) = \mu(F_3(\gamma(j))\mid
F_3(\gamma(o)))\mbox{.}
$$

Therefore we have a map:
$$
\beta_{ij}^* \otimes \Id\nolimits_{\mu(F_2(j)\mid F_2(o))} \quad
\colon
$$
$$
\s(F_3(\gamma(i)) /
F_3(\gamma(j)))\otimes_{\cc}\mu(F_3(\gamma(j))\mid F_3(\gamma(o)))
\arrow{e} \qquad \qquad \qquad \qquad \qquad \qquad \qquad \qquad
\qquad
$$
$$
\qquad \qquad \qquad \qquad \qquad \qquad \qquad \qquad \qquad
 \arrow{e} \s(F_2(i) / F_2(j))\otimes_{\cc}\mu(F_2(j)\mid
F_2(o)) \mbox{.}
$$
These maps (for various $i\ge j\in I_2$) are  compatible when we
take the projective limits according to formulas which  define
$\s_{F_3(\gamma(o))}(E_3)$ and $\s_{F_2(o)}(E_2)$. Therefore we
obtain the map $\beta^*$ from the maps $\beta_{ij}^* \otimes
\Id_{\mu(F_2(j)\mid F_2(o))}$. The proposition is proved.

\begin{nt} {\em
Dually to the maps, which we considered above, we have a
well-defined map
$$
\beta_{*} \quad \colon \quad \s'_{F_2(o)}(E_2)\arrow{e}
\s'_{F_3(\gamma(o))}(E_3)
$$
such that the maps $\beta^*$ and $\beta_*$ are conjugate maps with
respect to the pairings \linebreak$<\cdot,\cdot>_{\s_{F_2(o)}(E_2)}$
and $<\cdot,\cdot>_{\s_{F_3(\gamma(o))}(E_3)}$. }
\end{nt}

\subsection{}
 Let
$$
0\arrow{e} E_1\arrow{e,t}{\alpha} E_2\arrow{e,t}{\beta}E_3\arrow{e}0
$$
be an admissible triple from category $C_2^{ar}$, where
$E_i=(I_i,F_i,V_i)$, \linebreak $1 \le i \le 3$. By definition,
there are order-preserving functions:
\begin{itemize}
\item[(i)]
$ \gamma \; \colon \; I_2\arrow{e} I_3\quad\mbox{such that}\quad
\beta(F_2(i))=F_3(\gamma(i))\quad\mbox{for any}\quad i\in I_2
\mbox{,} $
\item[(ii)]
$ \veps \; \colon \; I_2\arrow{e} I_1\quad\mbox{such that}\quad
F_2(i)\cap V_1=F_1(\veps(i))\quad\mbox{for any}\quad i\in I_2
\mbox{.} $
\end{itemize}

\begin{prop} \label{pr29}
We suppose that $E_3$ is a $df$-object. Let $o\in I_2$. Then there
is the direct image
$$
\alpha_* \quad \colon \quad \s_{F_{1}(\veps(o))}(E_1)\arrow{e}
\s_{F_{2}(o)}(E_2) \mbox{.}
$$
\end{prop}
\proof. For any $i,j\in I_2$ we have canonically
\begin{equation} \label{natural}
\mu(F_2(i)\mid F_2(j))=\mu(F_1(\veps(i))\mid
F_1(\veps(j)))\otimes_{\cc} \mu(F_3(\gamma(i))\mid F_3(\gamma(j)))
\mbox{.}
\end{equation}
Since $E_3$ is a $df$-object, for any $i \ge j \in I_2$ we have a
canonical element  $\delta_{ij} \in \mu(F_3(\gamma(i)) /
F_3(\gamma(j)))$ constructed in remark~\ref{remark9}.  Therefore we
have a canonical ele\-ment $\delta_{ij}\in \mu(F_3(\gamma(i))\mid
F_3(\gamma(j)))$ for any $i,j\in I_2$ such that $\delta_{ij}\otimes
\delta_{jk}=\delta_{ik}$ for any $i,j,k \in I_2$. Therefore
 we have canonically $\mu(F_3(\gamma(i)) \mid
F_3(\gamma(j)))=\cc$ for any $i,j\in I_2$. (See also
remark~\ref{remark9}).

Thus from formula (\ref{natural}) we have for any $i,j\in I_2$
\begin{equation} \label{naturalnatural}
\mu(F_2(i)\mid F_2(j))=\mu(F_1(\veps(i))\mid F_1(\veps(j))) \mbox{.}
\end{equation}

For any $i\ge j\in I_2$ we have an admissible triple from
$C_1^{ar}$:
$$
0\arrow{e}
\frac{F_1(\veps(i))}{F_1(\veps(j))}\arrow{e,t}{\alpha_{ij}}\frac{F_2(i)}{F_2(j)}
\arrow{e} \frac{F_3(\gamma(i))}{F_3(\gamma(j))}\arrow{e}0 \mbox{.}
$$

We have a well-defined map:
$$
(\alpha_{ij})_* \quad \colon \quad
\s(F_1(\veps(i))/F_1(\veps(j)))\arrow{e} \s(F_2(i)/F_2(j)) \mbox{,}
$$
because $\frac{F_3(\gamma(i))}{F_3(\gamma(j))}$ is a discrete object
from category $C_1^{ar}$.

From formula (\ref{naturalnatural}) we have that
$$
\mu(F_2(j) \mid F_2(o))=\mu(F_1(\veps(j)) \mid F_1(\veps(o))).
$$

Therefore we have a map:
$$
(\alpha_{ij})_*\otimes \Id\nolimits_{\mu(F_2(j) \mid F_2(o))} \quad
\colon
$$
$$
\s(F_1(\veps(i))/F_1(\veps(j)))\otimes_{\cc}\mu(F_1(\veps(j))\mid
F_1(\veps(o))) \arrow{e} \qquad \qquad \qquad \qquad \qquad \qquad
\qquad \qquad
$$
$$
\qquad \qquad \qquad \qquad \qquad \qquad \qquad \qquad \arrow{e}
\s(F_2(i)/F_2(j)) \otimes_{\cc} \mu(F_2(j)\mid F_2(o)).
$$
These maps (for various $i\ge j\in I_2$) are  compatible when we
take the projective limits according to the formulas which  define
the spaces $\s_{F_1(\veps(o))}(E_1)$ and $\s_{F_2(o)}(E_2)$.
Therefore we obtain the map $\alpha_*$ from the maps
$(\alpha_{ij})_*\otimes \Id_{\mu(F_2(j) \mid F_2(o))}$. The
proposition is proved.

\begin{nt} {\em
Dually to the maps, which we considered above, we have a
well-defined map
$$
\alpha^{*} \quad \colon \quad \s'_{F_2(o)}(E_2) \arrow{e}
\s'_{F_1(\veps(o))}(E_1)
$$
such that the maps $\alpha_*$ and $\alpha^*$ are conjugate maps with
respect to the pairings \linebreak
$<\cdot,\cdot>_{\s_{F_2(o)}(E_2)}$ and
$<\cdot,\cdot>_{\s_{F_1(\veps(o))}(E_1)}$. }
\end{nt}

\section{Composition of maps and base change rules.} \label{bcr2}
In this section we consider  base change rules
(subsection~\ref{ss5.7.1}) and rules of composition of maps
(subsection~\ref{ss5.7.2}).

\subsection{} \label{ss5.7.1} Let
$$
0 \arrow{e} E_1  \arrow{e,t}{\alpha}  E_2 \arrow{e,t}{\beta}   E_3
\arrow{e} 0
$$
be an admissible triple from category $C_2^{ar}$, where $E_i = (I_i,
F_i, V_i)$, $ 1 \le i \le 3$.

Let
$$
0 \arrow{e} D  \arrow{e,t}{\gamma}  E_3 \arrow{e,t}{\delta}   B
\arrow{e} 0
$$
be an admissible triple from category  $C_2^{ar}$, where $D = (R, S,
Y)$, $B = (T, U, W)$. Then we have the following commutative diagram
of morphisms from category $C_2^{ar}$.

\begin{equation} \label{zvezda}
\begin{diagram}
\node[3]{0} \arrow{s} \node{0} \arrow{s} \\
 \node{0} \arrow{e} \node{E_1}
\arrow{e,t}{\gamma_{\alpha}} \arrow{s,=} \node{E_2
\mathop{\times}_{E_3} D} \arrow{e,t}{\gamma_{\beta}}
\arrow{s,l}{\beta_{\gamma}} \node{D} \arrow{e} \arrow{s,r}{\gamma} \node{0} \\
\node{0} \arrow{e} \node{E_1}  \arrow{e,t}{\alpha}  \node{E_2}
\arrow{e,t}{\beta}  \arrow{s,l}{\beta_{\delta}} \node{E_3} \arrow{e}
\arrow{s,r}{\delta} \node{0} \\
\node[3]{B} \arrow{r,=} \arrow{s} \node{B} \arrow{s} \\
\node[3]{0} \node{0}
\end{diagram}
\end{equation}

Let $X' = E_2 \mathop{\times}\limits_{E_3} D = (N, Q, X)$ as object
of category $C_2^{ar}$. In this diagram two horizontal triples and
two vertical triples are admissible triples from category
$C_2^{ar}$.

We choose any  $i \le j \in I_2$. Then diagram~(\ref{zvezda})
induces the following commutative diagram of morphisms  from
category $C_1^{ar}$ (see the construction of fibered product from
lemma~\ref{lbc2}).

\begin{equation} \label{2zvezda}
\begin{diagram}
\node[3]{0} \arrow{s} \node{0} \arrow{s} \\
 \node{0} \arrow{e}
 \node{\frac{F_2(j) \cap V_1}{F_2(i) \cap V_1}}
\arrow{e,t}{(\gamma_{\alpha})_{ji}} \arrow{s,=} \node{\frac{F_2(j)
\cap X}{F_2(i) \cap X}} \arrow{e,t}{(\gamma_{\beta})_{ji}}
\arrow{s,l}{(\beta_{\gamma})_{ji}} \node{\frac{\beta(F_2(j)) \cap
Y}{\beta(F_2(i)) \cap Y}}
\arrow{e} \arrow{s,r}{\gamma_{ji}} \node{0} \\
\node{0} \arrow{e} \node{\frac{F_2(j) \cap V_1}{F_2(i) \cap V_1}}
\arrow{e,t}{\alpha_{ji}}  \node{\frac{F_2(j)}{F_2(i)}}
\arrow{e,t}{\beta_{ji}} \arrow{s,l}{(\beta_{\delta})_{ji}}
\node{\frac{\beta(F_2(j))}{\beta(F_2(i))}} \arrow{e}
\arrow{s,r}{\delta_{ji}} \node{0} \\
\node[3]{\frac{\delta \beta (F_2(j))}{\delta \beta (F_2(i))}}
\arrow{r,=} \arrow{s}
\node{\frac{\delta \beta (F_2(j))}{\delta \beta (F_2(i))}} \arrow{s} \\
\node[3]{0} \node{0}
\end{diagram}
\end{equation}

Here two horizontal triples and two vertical triples are admissible
triples from category
 $C_1^{ar}$. At the same time,
$$
\frac{F_2(j) \cap X}{F_2(i) \cap X} = \frac{F_2(j)}{F_2(i)}
\mathop{\times}_{\frac{\beta(F_2(j))}{\beta(F_2(i))}}
\frac{\beta(F_2(j)) \cap Y}{\beta(F_2(i)) \cap Y}
$$
as objects from category $C_1^{ar}$.

\begin{prop} \label{przam1}
Using notations of diagram~(\ref{zvezda}),  let $E_1$ be a
$c$-object, $B$ be a $d$-object. Let $o \in I_2$, $\mu \in
\mu(F_2(o) \cap V_1 \mid V_1)$, $ \nu \in \mu(\delta \beta (F_2(o))
\mid \{ 0\})$. Then the following formulas are satisfied.
\begin{enumerate}
\item
For any $f \in  \s_{F_2(o)} (E_2)$
\begin{equation} \label{ut1}
\gamma^*( \beta_*(f \otimes \mu) \otimes \nu )= (\gamma_{\beta})_*
(\beta_{\gamma}^* (f \otimes \nu) \otimes \mu) \mbox{.}
\end{equation}
\item
For any  $G \in \s'_{\beta(F_2(o))\cap Y} (D)$
\begin{equation} \label{ut2}
\beta^* ( \gamma_* (G \otimes \nu) \otimes \mu) =
(\beta_{\gamma})_*(\gamma_{\beta}^* (G \otimes \mu) \otimes \nu)
\mbox{.}
\end{equation}
\end{enumerate}
\end{prop}
\begin{nt} {\em
The formulation of this proposition is well-defined  (with respect
to definitions of direct and inverse images), since $\alpha =
\beta_{\gamma} \circ \gamma_{\alpha}$, $\beta_{\delta} = \delta
\circ \beta $.}
\end{nt}
\proof  of proposition~\ref{przam1}. We prove at first
formula~(\ref{ut1}). Let $j_0 \in I_2$ be such that $F_2(j_0) \cap
V_1 = V_1$. (Such a $j_0$ it is possible to choose, since $E_1$ is a
$c$-object.) Let $i_0 \in I_2$ be such that  $\beta \delta
(F_2(i_0)) = \{ 0\}$. (Such an $i_0$ it is possible to choose, since
 $B$ is a $d$-object.)
We consider arbitrary  $j \ge i \in I_2$ such that
  $j \ge j_0$, $i
\le i_0$. Then  the base change formula on objects from category
$C_1^{ar}$ is satisfied for the maps $\beta_{ji}$, $\gamma_{ji}$,
$(\gamma_{\beta})_{ji}$, $(\beta_{\gamma})_{ji}$ from cartesian
square of diagram~(\ref{2zvezda}). (We have to use the analog of
formula~(\ref{e''1}) from proposition~\ref{pr18} for objects of
$C_1^{ar}$, see remark~\ref{remm8}.) We multiply this formula by
corresponding measure spaces, and then we take the projective limit
with respect to all such $j \ge i \in I_2$ ($j \ge j_0$, $i \le
i_0$). Using the constructions of maps (direct and inverse images)
$\beta_*$ and $\gamma^*$ from propositions~\ref{pr26}
and~\ref{pr27},  we obtain the base change formula~(\ref{ut1}).

Formula~(\ref{ut2}) is a dual formula to formula~(\ref{ut1}), and it
can be obtained by analogous reasonings as in the proof of
formula~(\ref{ut1}), but one has to take the inductive limit with
respect to $j, i \in I_2$, and one has to use the analog of
formula~(\ref{e''2}) instead of the analog of  formula~(\ref{e''1})
(see remark~\ref{remm8}). The proposition is proved. \vspace{0.3cm}

\begin{prop} \label{przam2}
Using notations of diagram~(\ref{zvezda}), let   $E_1$ be a
$cf$-object, $B$ be a $df$-object. Let $o \in I_2$. Then the
following formulas are satisfied.
\begin{enumerate}
\item For any $f \in \s_{\beta(F_2(o)) \cap Y} (D)$
\begin{equation} \label{ut3}
\beta^* \gamma_* (f) = (\beta_{\gamma})_* \gamma_{\beta}^* (f)
\mbox{.}
\end{equation}
\item
For any $G \in \s'_{F_2(o)}(E_2)$
\begin{equation} \label{ut4}
\gamma^* \beta_* (G) = (\gamma_{\beta})_* \beta_{\gamma}^* (G)
\mbox{.}
\end{equation}
\end{enumerate}
\end{prop}
\proof. We prove formula~(\ref{ut3}). We consider arbitrary  $j \ge
i \in I_2$. Then  the base change formula on objects from category
$C_1^{ar}$ is satisfied for the maps $\beta_{ji}$, $\gamma_{ji}$,
$(\gamma_{\beta})_{ji}$, $(\beta_{\gamma})_{ji}$ from cartesian
square of diagram~(\ref{2zvezda}). (We have to use the analog of
formula~(\ref{pr18f5}) from proposition~\ref{pr18} for objects of
$C_1^{ar}$, see remark~\ref{remm8}). We take the projective limit
with respect to such $j \ge i \in I_2$. Using the constructions of
maps (inverse and direct images) $\beta^*$  and
 $\gamma_*$ from propositions~\ref{pr28}
and~\ref{pr29},  we obtain the base change formula~(\ref{ut3}).

Formula~(\ref{ut4}) can be obtained by analogous reasonings, but one
has to take the inductive limit with respect to $ i, j \in I_2$, and
one has to use the analog of formula~(\ref{pr18f6}) instead of the
analog of formula~(\ref{pr18f5}) (see remark~\ref{remm8}). The
proposition is proved. \vspace{0.3cm}

\begin{prop}
Using notations of diagram~(\ref{zvezda}), let   $E_1$ be a
$c$-object, $B$ be a $df$-object. Let $o \in I_2$, $\mu \in
\mu(F_2(o) \cap V_1 \mid V_1)$. Then the following formulas are
satisfied.
\begin{enumerate}
\item For any $f \in \s_{F_2(o) \cap X}(X')$, $\mu \in \mu(E_1)$
\begin{equation} \label{ut5}
\beta_* ((\beta_{\gamma})_* (f) \otimes \mu) = \gamma_*
(\gamma_{\beta})_* (f \otimes \mu) \mbox{.}
\end{equation}
\item For any $G \in \s_{\beta(F_2(o))}'(E_3)$,
$\mu \in \mu(E_1)$
\begin{equation} \label{ut6}
\gamma_{\beta}^* (\gamma^* (G) \otimes \mu ) = \beta_{\gamma}^*
\beta^* (G \otimes \mu) \mbox{.}
\end{equation}
\end{enumerate}
\end{prop}
\proof  is analogous to the proofs of proposition~\ref{przam1}
and~\ref{przam2}, and it is reduced to the corresponding analog of
formula~(\ref{pr18f4c})
 for objects of category $C_1^{ar}$ (to prove formula~(\ref{ut5}))
and to the corresponding analog of  formula~(\ref{pr18f4d}) for
objects of category $C_1^{ar}$ (to prove formula~(\ref{ut6})), see
remark~\ref{remm8}. The proposition is proved. \vspace{0.3cm}

\begin{prop}
Using notations of diagram~(\ref{zvezda}), let   $E_1$ be a
$cf$-object, $B$ be a $d$-object. Let $o \in I_2$, $ \nu \in
\mu(\delta \beta (F_2(o)) \mid \{ 0\})$. Then the following formulas
are satisfied.
\begin{enumerate}
\item
For any $f \in \s_{\beta(F_2(o))}(E_3)$
\begin{equation} \label{ut7}
\gamma_{\beta}^* \gamma^* (f \otimes \nu) = \beta_{\gamma}^*
(\beta^* (f) \otimes \nu) \mbox{.}
\end{equation}
\item For any $G \in \s_{F_2(o) \cap X}'(X')$
\begin{equation} \label{ut8}
\beta_* (\beta_{\gamma})_* (G \otimes \nu) = \gamma_*
((\gamma_{\beta})_* (G) \otimes \nu ) \mbox{.}
\end{equation}
\end{enumerate}
\end{prop}
\proof  is analogous to the proofs of proposition~\ref{przam1}
and~\ref{przam2}, and it is reduced to  the corresponding analog of
formula~(\ref{e''3})
 for objects of category $C_1^{ar}$ (to prove formula~(\ref{ut7}))
and to the corresponding analog of  formula~(\ref{e''4}) for objects
of category $C_1^{ar}$ (to prove formula~(\ref{ut8})), see
remark~\ref{remm8}. The proposition is proved. \vspace{0.3cm}

\subsection{} \label{ss5.7.2} We consider again
diagram~(\ref{zvezda}). If $E_1$ and $D$ are $c$-objects, then $E_2
\mathop{\times}\limits_{E_3} D$ is a $c$-object. Indeed, it follows
from the admissible triple from category $C_2^{ar}$:
\begin{equation} \label{diez}
0 \arrow{e} E_1  \arrow{e,t}{\gamma_{\alpha}} E_2
\mathop{\times}\limits_{E_3} D \arrow{e,t}{\gamma_{\beta}} D
\arrow{e} 0 \mbox{.}
\end{equation}

Let $o \in I_2$. Then from~(\ref{diez}) we obtain the canonical
isomorphism
\begin{equation} \label{2diez}
\mu(F_2(o) \cap X \mid X) = \mu (F_2(o) \cap V_1 \mid V_1)
\otimes_{\cc} \mu (\beta (F_2(o)) \cap Y \mid Y) \mbox{.}
\end{equation}
(The subspaces, which  appear in formula~(\ref{2diez}), are elements
of the filtration of corresponding objects from category $C_2^{ar}$.
Therefore the spaces of virtual measures are well-defined, see
section~\ref{vm}.)

We have from diagram~(\ref{zvezda}) the following admissible triple
from  $C_2^{ar}$:
$$
0 \arrow{e} E_2 \mathop{\times}\limits_{E_3} D
\arrow{e,t}{\beta_{\gamma}} E_2 \arrow{e,t}{\delta  \beta} B
\arrow{e} 0 \mbox{.}
$$

\begin{prop} \label{prop34}
Using notations of diagram~(\ref{zvezda}), let   $E_1$ and $D$ be
$c$-objects. Let $o \in I_2$, $\mu \in \mu(F_2(o) \cap V_1 \mid
V_1)$, $\nu \in \mu (\beta (F_2(o)) \cap Y \mid Y)$. Then the
following formulas are satisfied.
\begin{enumerate}
\item
For any $f \in \s_{F_2(o)}(E_2)$
\begin{equation} \label{ut9}
(\delta \beta)_* (f \otimes (\mu \otimes \nu)) = \delta_* ( \beta_*
(f \otimes \mu) \otimes \nu) \mbox{.}
\end{equation}
\item
For any $G \in \s'_{\delta \beta (F_2(o))}(B)$
\begin{equation} \label{ut10}
(\delta \beta)^* (G \otimes (\mu \otimes \nu))= \beta^* (\delta^* (G
\otimes \nu) \otimes \mu) \mbox{.}
\end{equation}
\end{enumerate}
\end{prop}
\begin{nt} {\em
Due to~(\ref{2diez}) we have that $\mu \otimes \nu \in \mu(F_2(o)
\cap X \mid X)$. }
\end{nt}
\proof of  proposition~\ref{prop34}. Let $j_0 \in I_2$ be such that
$F_2(j_0) \cap V_1 = V_1$ and $\quad \beta (F_2(j_0)) \cap Y =Y$. We
consider arbitrary  $j \ge i \in I_2$ such that $j \ge j_0$. Then
the analog of formula~(\ref{e1})  from proposition~\ref{pr16} for
objects of category $C_1^{ar}$ (see remark~\ref{remm8}) is satisfied
for the maps $\beta_{ji}$, $\delta_{ji}$ and $(\beta_{\delta})_{ji}
= \delta_{ji} \beta_{ji}$ (see diagram~(\ref{2zvezda})). We multiply
this formula by corresponding measure spaces, and then we take the
limit with respect to  $j \ge i \in I_2$ ($j \ge j_0$). Using
explicit construction of direct image from proposition~\ref{pr26},
we obtain formula~(\ref{ut9}).

Formula~(\ref{ut10}) can be obtained by analogous reasonings, but
one has to use the analog of  formula~(\ref{e2}) from
proposition~\ref{pr16} for objects of category $C_1^{ar}$ (see
remark~\ref{remm8}). The proposition is proved. \vspace{0.3cm}

Using notations of diagram~(\ref{zvezda}), we suppose that $E_1$ and
$D$ are $cf$-objects. Then from admissible triple from category
$C_2^{ar}$~(\ref{diez}) it follows that $E_2
\mathop{\times}\limits_{E_3} D$ is also a $cf$-object.

\begin{prop} \label{predut}
Using notations of diagram~(\ref{zvezda}), let   $E_1$ and $D$ be
$cf$-objects. Let $o \in I_2$. Then the following formulas are
satisfied.
\begin{enumerate}
\item
For any $f \in \s_{\delta \beta (F_2(o))}(B)$
\begin{equation} \label{ut11}
(\delta \beta)^* (f )= \beta^* \delta^* (f) \mbox{.}
\end{equation}
\item
For any $G \in \s'_{F_2(o)}(E_2)$
\begin{equation} \label{ut12}
(\delta \beta)_* (G ) = \delta_*  \beta_* (G) \mbox{.}
\end{equation}
\end{enumerate}
\end{prop}
\proof. We consider arbitrary elements $j \ge i \in I_2$. Then the
analog of formula~(\ref{pr16f5}) from proposition~\ref{pr16} for
objects of category $C_1^{ar}$ (see remark~\ref{remm8}) is satisfied
for the maps $\beta_{ji}$, $\delta_{ji}$ and $(\beta_{\delta})_{ji}
= \delta_{ji} \beta_{ji}$ (see diagram~(\ref{2zvezda})). We multiply
this formula by corresponding measure spaces, and then we take the
limit with respect to  $j \ge i \in I_2$. Using explicit
construction of inverse image from proposition~\ref{pr28}, we obtain
formula~(\ref{ut11}).

Formula~(\ref{ut12}) can be obtained by analogous reasonings, but
one has to use the analog of formula~(\ref{pr16f6}) from
proposition~\ref{pr16} for objects of category~$C_1^{ar}$ (see
remark~\ref{remm8}). The proposition is proved. \vspace{0.3cm}

A little bit rewriting  (and redenoting) diagram~(\ref{zvezda}), we
obtain the following commutative diagram  of morphisms between
objects from category $C_2^{ar}$:
\begin{equation} \label{3zvezda}
\begin{diagram}
\node[3]{0} \arrow{s} \node{0} \arrow{s} \\
 \node{0} \arrow{e} \node{E_1}
\arrow{e,t}{\alpha} \arrow{s,=} \node{E_2} \arrow{e,t}{\beta}
\arrow{s,l}{\alpha'} \node{E_3} \arrow{e} \arrow{s,r}{\beta_{\alpha'}} \node{0} \\
\node{0} \arrow{e} \node{E_1}  \arrow{e,t}{\alpha' \alpha} \node{H'}
\arrow{e,t}{\alpha'_{\beta}}  \arrow{s,l}{\theta} \node{E_3
\mathop{\amalg}\limits_{E_2} H'} \arrow{e}
\arrow{s,r}{\theta'} \node{0} \\
\node[3]{L'} \arrow{r,=} \arrow{s} \node{L'} \arrow{s} \\
\node[3]{0} \node{0}
\end{diagram}
\end{equation}
In this diagram two horizontal and two vertical triples are
admissible triples from category  $C_2^{ar}$.

Let $E_i = (I_i, F_i, V_i)$ ($1 \le i \le 3$), $H' = (J', T', W')$
as objects from category $C_2^{ar}$. From admissible triple from
category $C_2^{ar}$
\begin{equation} \label{zv}
0 \arrow{e} E_3  \arrow{e,t}{\beta_{\alpha'}} E_3
\mathop{\amalg}\limits_{E_2} H' \arrow{e,t}{\theta'} L' \arrow{e} 0
\end{equation}
we obtain that if  $E_3$ and $L'$ are $d$-objects, then $ E_3
\mathop{\amalg}\limits_{E_2} H'$ is also a $d$-object.

Let $ o \in J'$ be any, then from~(\ref{zv}) we canonically have
\begin{equation} \label{2zv}
\mu (\alpha_{\beta}' (T'(o)) \mid \{ 0 \} ) = \mu(\beta(T'(o) \cap
V_2) \mid \{ 0 \}) \otimes_{\cc} \mu(\theta(T'(o)) \mid \{ 0\})
\mbox{.}
\end{equation}
(The subspaces, which  appear in formula~(\ref{2zv}), are elements
of the filtration of corresponding objects from category $C_2^{ar}$.
Therefore the spaces of virtual measures are well-defined, see
section~\ref{vm}.)

From diagram~(\ref{3zvezda}) we have an admissible triple from
category $C_2^{ar}$
$$
0 \arrow{e} E_1  \arrow{e,t}{\alpha' \alpha} H'
\arrow{e,t}{\alpha_{\beta}'}  E_3 \mathop{\amalg}\limits_{E_2} H'
\arrow{e} 0 \mbox{.}
$$

\begin{prop}
Using notations of diagram~(\ref{3zvezda}), let   $E_3$ and $L'$ be
$d$-objects. Let $o \in J'$, $\mu \in \mu(\beta(T'(o) \cap V_2) \mid
\{ 0 \}) $, $\nu \in \mu(\theta(T'(o)) \mid \{ 0\})$. Then the
following formulas are satisfied.
\begin{enumerate}
\item
For any $f \in \s_{T'(o)} (H')$
\begin{equation} \label{ut13}
(\alpha' \alpha)^* (f \otimes (\mu \otimes \nu))= \alpha^*
((\alpha')^*(f \otimes \nu) \otimes \mu ) \mbox{.}
\end{equation}
\item
For any $G \in \s'_{T'(o) \cap V_1} (E_1)$
\begin{equation} \label{ut14}
(\alpha' \alpha)_* (G \otimes (\mu \otimes \nu)) = (\alpha')_*
(\alpha_*(G \otimes \mu) \otimes \nu) \mbox{.}
\end{equation}
\end{enumerate}
\end{prop}
\begin{nt} {\em
Due to~(\ref{2zv}) we have that $\mu \otimes \nu \in
\mu(\alpha_{\beta}' (T'(o)) \mid \{ 0 \})$. }
\end{nt}
\proof of proposition. Formula~(\ref{ut13}) can be reduced to the
analog of  formula~(\ref{e'1}) from proposition~\ref{pr17} for
objects of category $C_1^{ar}$ (see remark~\ref{remm8}) after the
taking the projective limit with respect to $j \ge i \in J'$, where
$i \le i_0$, and $i_0$ is chosen such that $\theta(T'(i_0)) = \{
0\}$ and $\beta (T'(i_0) \cap V_2) = \{ 0\}$. At that we use the
explicit construction of inverse image from proposition~\ref{pr27}.
(Compare
 with the proof of proposition~\ref{predut}.)

Formula~(\ref{ut14}) can be proved analogously, but one has to use
the analog of  formula~(\ref{e'2}) from proposition~\ref{pr17} for
objects of category~$C_1^{ar}$ (see remark~\ref{remm8}) and take the
inductive limit. The proposition is proved. \vspace{0.3cm}

Using notations of diagram~(\ref{3zvezda}) we suppose that  $E_3$
and $L'$ are $df$-objects. Then from admissible triple from category
$C_2^{ar}$~(\ref{zv}) it follows that $E_3
\mathop{\amalg}\limits_{E_2} H'$ is also a  $df$-object.

\begin{prop} Using notations of diagram~(\ref{3zvezda}), let  $E_3$ and $L'$
are $df$-objects. Let $o \in J'$. Then the following formulas are
satisfied.
\begin{enumerate}
\item
For any $f \in \s_{T'(o) \cap V_1} (E_1)$
\begin{equation} \label{ut15}
(\alpha' \alpha)_* (f ) = (\alpha')_* \alpha_*(f ) \mbox{.}
\end{equation}
\item
For any $G \in \s'_{T'(o)} (H')$
\begin{equation} \label{ut16}
(\alpha' \alpha)^* (G )= \alpha^* ((\alpha')^*(G ) \mbox{.}
\end{equation}
\end{enumerate}
\end{prop}
\proof of formula~(\ref{ut15}) can be reduced to the analog of
formula~(\ref{pr17f5}) from proposition~\ref{pr17} for objects of
category~$C_1^{ar}$ (see remark~\ref{remm8}) after the choice of
 $j \ge i \in J'$ and taking the projective limit with respect to this  $j \ge
i \in J'$.  At that we use the explicit construction of direct image
from proposition~\ref{pr29}.

The proof of formula~(\ref{ut16}) is analogous, but one has to take
the inductive limit with respect to   $j \ge i \in J'$ and use the
analog of formula~(\ref{pr17f6}) instead of the analog of
formula~(\ref{pr17f5}) (see remark~\ref{remm8}). The proposition is
proved. \vspace{0.3cm}

\section{Fourier transform and direct and inverse images}
\label{ftdi} Let
$$
0\arrow{e} E_1\arrow{e,t}{\alpha} E_2\arrow{e,t}{\beta}E_3\arrow{e}0
$$
be an admissible triple from category $C_2^{ar}$, where
$E_i=(I_i,F_i,V_i)$, \linebreak $1 \le i \le 3$. We suppose that
$E_i$ ($1 \le i\le 3$) is a {\it complete} object from category
$C_2^{ar}$. Then $\check{\check{E_i}}=E_i$, $1 \le i \le 3$.

By definition, there are order-preserving  functions:
\begin{itemize}
\item[(i)]
$ \gamma \; \colon \; I_2\arrow{e} I_3\quad\mbox{such that}\quad
\beta(F_2(i))=F_3(\gamma(i))\quad\mbox{for any}\quad i\in I_2
\mbox{;} $
\item[(ii)]
$ \veps \; \colon \; I_2\arrow{e} I_1\quad\mbox{such that}\quad
F_2(i)\cap V_1=F_1(\veps(i))\quad\mbox{for any}\quad i\in
I_2\mbox{.} $
\end{itemize}
We consider the following admissible triple from category
$C_2^{ar}$:
$$
0\arrow{e} \check{E_3} \arrow{e,t}{\check\beta} \check{ E_2}
\arrow{e,t}{\check\alpha} \check{E_1} \arrow{e}0 \mbox{.}
$$
We have $\check{E_i}=(I^0_i,F^0_i,\check{V_i})$, $1 \le i \le 3$. We
recall that the partially ordered set $I^{0}_i$ ($1 \le i \le 3$) is
equal to $I_i$ as a set, but has inverse partial order comparing
with $I_i$ ($1 \le i \le 3$). We have that for any $i,j \in I_2$:
\begin{equation} \label{bigstar}
\begin{array}{c}
\check\alpha (F^{0}_2(i))=F^{0}_1(\veps(i)) \\[0.2cm]
F_2^0(i)\cap \check{V_3}=F^{0}_3(\gamma(i)) \mbox{.}
\end{array}
\end{equation}

For any $i,j\in I_2$ we have from formula~(\ref{fff}) that
\begin{equation} \label{bigstarbigstar}
\begin{array}{c}
\mu(F_1(\veps(i))\mid F_1(\veps(j)))=\mu(F^{0}_1(\veps(i))\mid
F^{0}_1(\veps(j))) \\[0.2cm]
\mu(F_3(\gamma(i))\mid F_3(\gamma(j)))=\mu(F^{0}_3(\gamma(i))\mid
F^{0}_3(\gamma(j))) \mbox{.}
\end{array}
\end{equation}

We have the following statements.
\begin{prop} \label{pppp}
Let $o \in I_2$. The following diagrams are commutative.
\begin{enumerate}
\item If $E_1$ is a $c$-object, then
$$
\begin{diagram}
\node{\s_{F_2(o)}(E_2)\otimes_{\cc}\mu(F_1(\veps(o))\mid V_1)}
\arrow{e,t}{\beta_*} \arrow{s,l}{\F\otimes
\Id_{\mu(F_1(\veps(o))\mid V_1)}}
\node{\s_{F_3(\gamma(o))}(E_3)}\arrow{s,r}{\F}\\
\node{\s_{F^{0}_2(o)}(\check{
E_2})\otimes_{\cc}\mu(F^{0}_1(\veps(o))\mid \{0\})}
\arrow{e,t}{(\check\beta)^*} \node{\D_{F^{0}_3(\gamma(o))}(\check
E_3) \mbox{.}}
\end{diagram}
$$
\item If $E_3$ is a $d$-object, then
$$
\begin{diagram}
\node{\s_{F_2(o)}(E_2)\otimes_{\cc}\mu(F_3(\gamma(o))\mid\{0\})}
\arrow{e,t}{\alpha^*} \arrow{s,l}{\F\otimes
\Id_{\mu(F_3(\gamma(o))\mid\{0\})}}
\node{\s_{F_1(\veps(o))}(E_1)}\arrow{s,r}{\F}\\
\node{\s_{F^{0}_2(o)}(\check{
E_2})\otimes_{\cc}\mu(F^{0}_3(\gamma(o))\mid \check{ V_3})}
\arrow{e,t}{(\check\alpha)_*}
\node{\s_{F^{0}_1(\veps(o))}(\check{E_1}) \mbox{.}}
\end{diagram}
$$
\item If $E_1$ is a $c$-object, then
$$
\begin{diagram}
\node{\s'_{F_3(\gamma(o))}(E_3)\otimes_{\cc}\mu(F_1(\veps(o))\mid
V_1)} \arrow{e,t}{\beta^*}\arrow{s,l}{\F\otimes
\Id_{\mu(F_1(\veps(o))\mid V_1)}}
\node{\s'_{F_2(o)}(E_2)}\arrow{s,r}{\F}\\
\node{\s'_{F^{0}_3(\gamma(o))}(\check{
E_3})\otimes_{\cc}\mu(F^{0}_1(\veps(o))\mid \{0\})}
\arrow{e,t}{(\check\beta)_*} \node{\s'_{F^{0}_2(o)}(\check{E_2})
\mbox{.}}
\end{diagram}
$$
\item If $E_3$ is a $d$-object, then
$$
\begin{diagram}
\node{\s'_{F_1(\veps(o))}(E_1)\otimes_{\cc}\mu(F_3(\gamma(o))\mid
\{0\})} \arrow{e,t}{\alpha_*}\arrow{s,l}{\F\otimes
\Id_{\mu(F_3(\gamma(o))\mid \{0\})}}
\node{\s'_{F_2(o)}(E_2)}\arrow{s,r}{\F}\\
\node{\s'_{F^{0}_1(\veps(o))}(\check{
E_1})\otimes_{\cc}\mu(F^{0}_3(\gamma(o))\mid \check{V_3})}
\arrow{e,t}{(\check\alpha)^*} \node{\s'_{F^{0}_2(o)}(\check{ E_2})
\mbox{.}}
\end{diagram}
$$
\item If $E_1$ is a $cf$-object, then
$$
\begin{diagram}
\node{\s_{F_3(\gamma(o))}(E_3)} \arrow{e,t}{\beta^*}\arrow{s,l}{\F}
\node{\s_{F_2(o)}(E_2)}\arrow{s,r}{\F}\\
\node{\s_{F^{0}_3(\gamma(o))}(\check{E_3})}
\arrow{e,t}{(\check\beta)_*} \node{\s_{F^{0}_2(o)}(\check{E_2})
\mbox{.}}
\end{diagram}
$$
\item If $E_3$ is a $df$-object, then
$$
\begin{diagram}
\node{\s_{F_1(\veps(o))}(E_1)}\arrow{e,t}{\alpha_*}\arrow{s,l}{\F}
\node{\s_{F_2(o)}(E_2)}\arrow{s,r}{\F}\\
\node{\s_{F^{0}_1(\veps(o))}(\check{E_1})}
\arrow{e,t}{(\check\alpha)^*} \node{\s_{F^{0}_2(o)}(\check{E_2})
\mbox{.} }
\end{diagram}
$$
\item If $E_1$ is a $cf$-object, then
$$
\begin{diagram}
\node{\s'_{F_2(o)}(E_2)} \arrow{e,t}{\beta_*}\arrow{s,l}{\F}
\node{\s'_{F_3(\gamma(o))}(E_3)}\arrow{s,r}{\F}\\
\node{\s'_{F^{0}_2(o)}(\check{E_2})} \arrow{e,t}{(\check\beta)^*}
\node{\s'_{F^{0}_3(\gamma(o))}(\check{E_3}) \mbox{.}}
\end{diagram}
$$
\item If $E_3$ is a $df$-object, then
$$
\begin{diagram}
\node{\s'_{F_2(o)}(E_2)} \arrow{e,t}{\alpha^*}\arrow{s,l}{\F}
\node{\s'_{F_1(\veps(o))}(E_1)}\arrow{s,r}{\F}\\
\node{\s'_{F^{0}_2(o)}(\check{E_2})} \arrow{e,t}{(\check\alpha)_*}
\node{\s'_{F^{0}_1(\veps(o))}(\check{E_1}) \mbox{.}}
\end{diagram}
$$
\end{enumerate}
\end{prop}
\proof  follows from formulas (\ref{bigstar}) and
(\ref{bigstarbigstar}), constructions of two-dimensional direct and
inverse images from section~\ref{secdi}, construction of
two-dimensional Fourier transform from section~\ref{ft}, and
corresponding analogs of propositions~\ref{prpr} and~\ref{comdia}
for objects of category $C_1^{ar}$ (see remark~\ref{remm8}), which
connect direct and inverse images with Fourier transform on objects
of category $C_1^{ar}$. The proposition is proved.

\section{Two-dimensional Poisson formulas} \label{2PF} In this
section we obtain two-dimensional Poisson formulas: the Poisson
formula I (subsection~\ref{s5.9.1}) and the Poisson formula II
(subsection~\ref{s5.9.2}).
\subsection{Poisson formula I} \label{s5.9.1} Let
$$
0\arrow{e} E_1\arrow{e,t}{\alpha} E_2\arrow{e,t}{\beta}E_3\arrow{e}0
$$
be an admissible triple from category  $C_2^{ar}$, where \linebreak
$E_i=(I_i,F_i,V_i)$, $1 \le i \le 3$. We suppose that $E_i$ ($1 \le
i\le 3$) is a {\it complete} object from category $C_2^{ar}$.
 By
definition, there are order-preserving  functions:
\begin{itemize}
\item[(i)]
$ \gamma \; \colon \; I_2\arrow{e} I_3\quad\mbox{such that}\quad
\beta(F_2(i))=F_3(\gamma(i))\quad\mbox{for any}\quad i\in I_2
\mbox{;} $
\item[(ii)]
$ \veps \; \colon \; I_2\arrow{e} I_1\quad\mbox{such that}\quad
F_2(i)\cap V_1=F_1(\veps(i))\quad\mbox{for any}\quad i\in I_2
\mbox{.} $
\end{itemize}

Let $o\in I_2$. If we suppose that $E_1$ is a $c$-object, then from
any \linebreak $\mu\in \mu(F_1(\veps(o))\mid V_1)$ we canonically
construct $\znakneponyaten_{\mu} \in \s'_{F_1(\veps(o))}(E_1)$ in
the following way. By definition, we have
\begin{equation} \label{nabla}
\s'_{F_1(\veps(o))}(E_1)=\mathop{\Lim_{\lto}}_{k\in I_1}
\s'(V_1/F_1(k))\otimes_{\cc}\mu(F_1(\veps(o))\mid k) \mbox{.}
\end{equation}
By definition, $\mu(F_1(\veps(o))\mid V_1)\subset
\s'(V_1/F_1(\veps(o)))$. Therefore $\mu\in \s'(V_1/F_1(\veps(o)))$.
In formula (\ref{nabla}) we take $k=\veps(o)$. Then
$\mu(F_1(\veps(o))\mid k)=\cc$ and we put \linebreak $1 \in
\mu(F_1(\veps(o))\mid \veps(o))$. Therefore $\mu$ gives
$\znakneponyaten_\mu\in \s'_{F_1(\veps(o))}(E_1)$ by formula
(\ref{nabla}).

Let $o \in I_2$. If we suppose that $E_3$ is a $d$-object, then from
any \linebreak $\nu \in \mu(F_3(\gamma(o))\mid \{0\})$ we
canonically construct $\delta_\nu \in \s'_{F_3(\gamma(o))}(E_3)$ in
the following way. By definition, we have
\begin{equation} \label{nablanabla}
\s'_{F_3(\gamma(o))}(E_3)=\mathop{\Lim_{\lto}}_{k\in I_2}
\s'(F_3(k))\otimes_{\cc}\mu(F_3(\gamma(o))\mid \{0\}) \mbox{.}
\end{equation}
In formula (\ref{nablanabla}) we take  $k=\gamma(o)$. Then we put
$\delta \in \s'(F_3(\gamma(o)))$ in formula~(\ref{nablanabla}),
where  $<\delta,f> \eqdef f(0)$, $f\in \s(F_3(\gamma(o)))$ with
respect to the pairing~(\ref{par}). (The element $\delta$ can be
also defined as the projective limit of Dirac delta functions
according to formula~(\ref{s'}).) We put also $\nu \in
\mu(F_3(\gamma(o))\mid \{0\})$ in formula~(\ref{nablanabla}). Thus
we have constructed $\delta_\nu\in \s'_{F_3(\gamma(o))}(E_3)$.

We note that
\begin{equation} \label{triangle}
\F(\znakneponyaten_\mu)=\delta_\mu\quad\mbox{and}\quad
\F(\delta_\nu)=\znakneponyaten_\nu
\end{equation}
under the maps
$$\F \quad \colon \quad \s'_{F_1(\veps(o))}(E_1)\arrow{e}
\s'_{F^{0}_1(\veps(o))}(\check{E_1}) \quad  \mbox{and}
$$
$$
\F \quad \colon \quad \s'_{F_3(\gamma(o))}(E_3)\arrow{e}
\s'_{F^{0}_3(\gamma(o))}(\check{ E_3}) \mbox{.}
$$
(We used that
$$
\mu(F_1(\veps(o))\mid V_1)=\mu(F^{0}_1(\veps(o))\mid
\{0\})\quad\mbox{and}\quad \mu(F_3(\gamma(o))\mid
\{0\})=\mu(F^{0}_3(\gamma(o))\mid \check{V_3}) \mbox{.})
$$

Now we suppose that simultaneously $E_1$ is a $c$-object, and $E_3$
is a $d$-object. Let $o \in I_2$. Then from any $\mu\in
\mu(F_1(\veps(o))\mid V_1)$, $\nu\in \mu(F_3(\gamma(o))\mid \{0\})$
it is well-defined {\em the characteristic function}
\begin{equation}
\delta_{E_1,\mu\otimes \nu}\eqdef\alpha_*(\znakneponyaten_\mu\otimes
\nu)\in \s'_{F_2(o)}(E_2).
\end{equation}

We have the following statement.
\begin{lemma} \label{lem5}
The following equality is satisfied
$$
\delta_{E_1,\mu\otimes \nu}=\beta^*(\delta_\nu\otimes \mu).
$$
\end{lemma}
\proof follows from constructions and corresponding statements for
the following admissible triple from category $C_1^{ar}$:
$$
0\arrow{e}
V_1/F_1(\veps(k_2))\arrow{e}\frac{F_2(k_1)}{F_2(k_2)}\arrow{e}
F_3(\gamma(k_1))\arrow{e}0 \mbox{,}
$$
where $k_1\ge o \ge k_2\in I_2$ and $F_1(\veps(k_1))=V_1$,
$F_3(\gamma(k_2))=\{0\}$. The lemma is proved. \vspace{0.3cm}

Now we have the following first two-dimensional Poisson formula.
\begin{Th}[\it Poisson formula I]. \label{tpf1}
Let $E_1$ be a $c$-object. Let $E_3$ be a $d$-object. Let $o \in
I_2$. Then for any $\mu\in \mu(F_1(\veps(o))\mid V_1)$, $\nu\in
\mu(F_3(\gamma(o))\mid \{0\})$ we have
$$
\F(\delta_{E_1,\mu\otimes \nu})=\delta_{\check{E_3},\nu\otimes \mu}
\mbox{.}
$$
\end{Th}
\proof follows from lemma~\ref{lem5}, from
formulas~(\ref{triangle}), and from proposition~\ref{pppp} about
connection of direct and inverse images with two-dimensional Fourier
transform. The theorem is proved.
 \vspace{0.3cm}

\begin{defin} \label{deflast}
Let $E=(I,F,V) \ic$. We will say that an element $g \in
\Aut_{C_2^{ar}}(E)$ satisfies condition~$(*)$ iff the following two
conditions are satisfied.
\begin{enumerate}
\item $g \in \Aut_{C_2^{ar}}(E)'$ (see definition~\ref{grdef}).
\item Let $i \ge j \in I$ be any. Let $g(F(i))=F(p)$, $g(F(j))=F(q)$
for some $p \ge q \in I$. Let $E_{j,i}=(I_{j,i}, F_{j,i},
F(i)/F(j))$, and $E_{q,p}=(I_{q,p}, F_{q,p}, F(p)/F(q))$ be the
corresponding objects of category $C_1^{ar}$ constructed from $E$.
Then for any $l \in I_{i,j}$ there is $r_l \in I_{q,p}$ such that
$gF_{j,i}(l)=F_{q,p}(r_l)$, and for $l_1 \le l_2 \in I_{j,i}$ we
have that $r_{l_1} \le r_{l_2} \in I_{q,p}$.
\end{enumerate}
\end{defin}

\begin{nt}{\em
It is not difficult to see that subgroups of $\Aut_{C_2^{ar}}(E)$
for $E \ic$ which was considered in remark~\ref{appl} consist of
elements satisfying condition~$(*)$.}
\end{nt}

We consider any admissible triple of complete objects from category
$C_2^{ar}$:
$$
0\arrow{e} E_1\arrow{e,t}{\alpha} E_2\arrow{e,t}{\beta}E_3\arrow{e}0
\mbox{,}
$$
where $E_i=(I_i,F_i,V_i)$, $1 \le i \le 3$, and $E_1$ is a
$c$-object, $E_3$ is a $d$-object. Let
  $g \in \Aut_{C_2}(E_2)$ satisfies condition~$(*)$ (see definition~\ref{deflast}). Then
   we have the following
admissible triple from category $C_2^{ar}$:
\begin{equation} \label{kvadrat}
0\arrow{e} g{E_1} \arrow{e,t}{g(\alpha)} g{ E_2}
\arrow{e,t}{g(\beta)} g {E_3} \arrow{e} 0 \mbox{,}
\end{equation}
where $g  E_2=E_2$, $g {E_3}= { E_2}/g {E_1}$.

From~(\ref{kvadrat}) we have canonically
$$
\mu(F_2(o)\mid gF_2(o))=\mu(F_2(o) \cap gV_1 \mid
gF_1(\veps(o)))\otimes_{\cc}\mu(g(\beta)(F_2(o)) \mid
gF_3(\gamma(o))) \mbox{.}
$$
Then for any $a\in\mu(F_2(0)\mid gF_2(0))$ we consider $a=b\otimes
c$, where
$$
b\in\mu(F_2(o) \cap gV_1 \mid gF_1(\veps(o))) \mbox{,} \quad
c\in\mu(g(\beta)(F_2(o)) \mid gF_3(\gamma(o))).
$$
Under conditions and notations of theorem~\ref{tpf1}  we have
$$
g\mu\in\mu(gF_1(\veps(o))\mid gV_1) \mbox{,} \qquad
g\nu\in\mu(gF_3(\gamma(o))\mid \{0\}) \mbox{.}
$$
We have $\widetilde
g=(g,a)\in\widehat{\Aut_{C_2^{ar}}(E_2)'}_{F_2(o)}$. Then we obtain
$$
R'_{\widetilde g}(\delta_{E_1,\mu\otimes\nu})=\delta_{gE_1,
(g\mu\otimes b)\otimes(g\nu\otimes c)}=\delta_{gE_1, g\mu\otimes
g\nu\otimes a}.
$$

\begin{cons}
Let $(g,a)\in\widehat{\Aut_{C_2^{ar}}(E_2)'}_{F_2(o)}$ such that $g \in
\Aut_{C_2^{ar}}(E_2)$ satisfies condition~(*) from
definition~\ref{deflast}. Then we have
$$
\F(\delta_{gE_1, g\mu\otimes g\nu\otimes a})=\delta_{\check
g^{-1}(\check{E_3}),\check g^{-1}(\nu)\otimes\check
g^{-1}(\mu)\otimes a}.
$$
\end{cons}
\proof follows from theorem~\ref{tpf1} (Poisson formula I)  and
section~\ref{sec5.5.5}. The corollary is proved. \vspace{0.3cm}

\subsection{Poisson formula II} \label{s5.9.2}
 Let
$$
0\arrow{e} E_1\arrow{e,t}{\alpha} E_2\arrow{e,t}{\beta}E_3\arrow{e}0
$$
be an admissible triple from category $C_2^{ar}$, where
$E_i=(I_i,F_i,V_i)$, $1 \le i \le 3$. We suppose that $E_i$ ($1 \le
i\le 3$) is a {\it complete} object from category $C_2^{ar}$.
 By
definition, there are order-preserving  functions:
\begin{itemize}
\item[(i)]
$ \gamma \; \colon \; I_2\arrow{e} I_3\quad\mbox{such that}\quad
\beta(F_2(i))=F_3(\gamma(i))\quad\mbox{for any}\quad i\in I_2
\mbox{;} $
\item[(ii)]
$ \veps \; \colon \; I_2\arrow{e} I_1\quad\mbox{such that}\quad
F_2(i)\cap V_1=F_1(\veps(i))\quad\mbox{for any}\quad i\in I_2
\mbox{.} $
\end{itemize}

If $E_1$ is a $cf$-object, then using remark~\ref{remark9} we have
$$
\s(E_1)=\mathop{\Lim_{\longleftarrow}}_{i\ge j\in I_1}
\s(F_1(i)/F_1(j)) \mbox{.}
$$
We have that $F_1(i)/F_1(j)$ is a compact object from category
$C_1^{ar}$ (for any $i\ge j\in I_1$). Therefore  $1\in
\s(F_1(i)/F_1(j))$. Taking the projective limits, we obtain that in
this case the following element is well-defined:
\begin{equation}
\znakneponyaten \in \s(E_1) \mbox{.}
\end{equation}

If $E_3$ is a $df$-object, then using remark~\ref{remark9} we have
$$
\s(E_3)=\mathop{\Lim_{\longleftarrow}}_{i\ge j\in I_3}
\s(F_3(i)/F_3(j)) \mbox{.}
$$
We have that $F_3(i)/F_3(j)$ is a discrete object from category
$C_1^{ar}$ (for any $i\ge j\in I_3$). Therefore  $\delta_0 \in
\s(F_3(i)/F_3(j))$, where $\delta_0(0)=1$ and $\delta_0(x)=0$ for $x
\ne 0$ ($x$ is from underlying abelian group of $C_1^{ar}$-object
$F_3(i)/F_3(j)$). Taking the projective limits, we obtain in this
case a well-defined element:
\begin{equation}
\delta_0 \in \s(E_3) \mbox{.}
\end{equation}

We note that
\begin{equation} \label{triangletriangle}
\F(\znakneponyaten)=\delta_0\quad\mbox{and}\quad
\F(\delta_0)=\znakneponyaten
\end{equation}
under the maps
$$\F \quad \colon \quad \s(E_1)\arrow{e} \s(\check E_1) \quad \mbox{and} $$
$$\F \quad \colon \quad \s(E_3)\arrow{e}
\s(\check E_3) \mbox{. }$$

Now we suppose that simultaneously $E_1$ is a $cf$-object, and $E_3$
is a $df$-object. Let $o \in I_2$. Then an element
\begin{equation}
\delta_{E_1} \eqdef \alpha_*(\znakneponyaten) \in \s_{F_2(o)}(E_2)
\end{equation}
is well-defined.

We have the following statement.
\begin{lemma} \label{lem6}
The following equality is satisfied
$$
\delta_{E_1}=\beta^*(\delta_0) \mbox{.}
$$
\end{lemma}
\proof follows from constructions and corresponding statements for
the following admissible triples from category $C_1^{ar}$ (for any
$i\ge j\in I_2$):
$$
0\arrow{e}\frac{F_1(\veps(i))}{F_1(\veps(j))}\arrow{e}
\frac{F_2(i)}{F_2(j)}\arrow{e}
\frac{F_3(\gamma(i))}{F_3(\gamma(j))}\arrow{e}0 \mbox{,}
$$
where $\frac{F_1(\veps(i))}{F_1(\veps(j))}$ is a compact object from
category $C_1^{ar}$, and $\frac{F_3(\gamma(i))}{F_3(\gamma(j))}$ is
a discrete object from category $C_1^{ar}$. The lemma is proved.
\vspace{0.3cm}

Now we have the following second two-dimensional Poisson formula.
\begin{Th}[\it Poisson formula II]. \label{th3}
Let $E_1$ be a $cf$-object. Let $E_3$ be a $df$-object. Let $o \in
I_2$. Then
$$
\F(\delta_{E_1})=\delta_{\check{E_3}} \mbox{.}
$$
\end{Th}
\proof follows from lemma~\ref{lem6}, from
formulas~(\ref{triangletriangle}), and from proposition~\ref{pppp}
about direct and inverse images and two-dimensional Fourier
transform. The theorem is proved. \vspace{0.3cm}

We consider any admissible triple of complete objects from category
$C_2^{ar}$:
$$
0\arrow{e} E_1\arrow{e,t}{\alpha} E_2\arrow{e,t}{\beta}E_3\arrow{e}0
\mbox{,}
$$
where $E_i=(I_i,F_i,V_i)$, $1 \le i \le 3$, and $E_1$ is a
$cf$-object, $E_3$ is a $df$-object. Let
  $g\in \Aut_{C_2^{ar}}(E_2)$ satisfies condition~$(*)$ from definition~\ref{deflast}. Then
   we have the following
admissible triple from category $C_2^{ar}$:
\begin{equation}
0\arrow{e} g{E_1} \arrow{e,t}{g(\alpha)} g { E_2}
\arrow{e,t}{g(\beta)} g {E_3} \arrow{e} 0 \mbox{,}
\end{equation}
where  $g { E_2}={E_2}$, $g {E_3}= { E_2}/g {E_1}$.

\begin{cons}
Let $g \in \Aut_{C_2^{ar}}(E_2)$ satisfies condition~$(*)$ from
definition~\ref{deflast}. Then we have
$$
\F(\delta_{gE_1})=\delta_{\check g^{-1}(\check{E_3})}.
$$
\end{cons}
\proof follows from theorem~\ref{th3} (Poisson formula II)  and
section~\ref{sec5.5.5}. The corollary is proved. \vspace{0.3cm}

\section{An example} \label{lastexample}
In this section we will compute some quotient groups of adelic
groups $\da_X$ of a two-dimensional scheme $X$.
\subsection{Some quotient groups of  adelic groups of algebraic
surface} \label{sqr}
 We recall that in example~\ref{exad} we
computed the group $\da_K/K$ for a number field $K$.

Let $D$ be an irreducible  projective curve over a field $k$. Let $p
\in D$ be a smooth $k$-rational point. Then by a similar reasoning
as in example~\ref{exad} we obtain the following exact sequence:
\begin{equation} \label{adseq}
0 \lto \prod_{q \in D, q \ne p} \hat{\oo}_q \lto \da_D/k(D) \lto K_p/ A_p  \lto 0  \mbox{,}
\end{equation}
where $k(D)$ is the field of rational functions on the curve $D$,
$\hat{\oo}_q$ is the completion by maximal ideal $m_q$ of the local
ring $\oo_q$ of point $q \in D$, $K_p=k(D)_p$ is the fraction field
of the ring $\hat{\oo}_p$ if $p $ is a smooth point (in general, the
ring $K_p$ is the localization of the ring $\hat{\oo}_p$ with
respect to the multiplicative system $\oo_p \setminus 0$), $A_p$ is
the ring of regular functions on the affine curve $D \setminus p$.
We note that $K_p \simeq k((t))$. If the field $k$ is a finite
field, then all terms of sequence~(\ref{adseq}) have the structure
of compact objects from the category $C_1^{ar}$ such that this
sequence is an admissible triple from the category $C_1^{ar}$.

We note that the adelic ring $\da_D$ and the field of rational
functions $k(D)$ are the groups which appear in adelic complex of
the curve $D$. Indeed, for any quasicoherent sheaf $\ff$ on the
curve $D$ there is the following adelic complex (see~\cite{PF},
\cite[\S 3.1]{Osi}, \cite{B}, \cite{H}):
$$
\da_{D,0}(\ff) \oplus \da_{D,1}(\ff) \lto \da_{D,01}(\ff)
$$
whose cohomology groups coincide with the cohomology groups $H^*(D, \ff)$. For the sheaf $\ff = \oo_D$ we have
$$
\da_{D, 0} (\oo_D) = k(D) \mbox{,} \qquad \da_{D, 1} (\oo_D)=\prod_{q \in D} \hat{\oo}_q \mbox{,}
$$
$$
\da_{D, 01}(\oo_D) = \da_D = {\prod_{q \in D}}' K_q \mbox{,}
$$
where $\prod'$ denotes the restricted (adelic) product. Thus, we obtained that
\begin{equation} \label{analog}
\da_D / k(D) = \da_{D, 01}(\oo_D) / \da_{D, 0}(\oo_D) \mbox{.}
\end{equation}

\vspace{0.5cm}

Now let $X$ be a two-dimensional scheme. For any quasicoherent sheaf $\ff$ on $X$ there is
the following adelic complex (see~\cite{PF}, \cite[\S 3.3]{Osi}):
\begin{equation} \label{adcompl}
\da_{X,0}(\ff) \oplus \da_{X,1}(\ff) \oplus \da_{X,2}(\ff) \lto
\da_{X, 01}(\ff) \oplus \da_{X,02}(\ff) \oplus \da_{X,12}(\ff) \lto
\da_{X,012}(\ff)
\end{equation}
whose cohomology groups coincide with the cohomology groups $H^*(X, \ff)$.

Let $X$ be a smooth irreducible surface over a field $k$, and $C$ be a divisor on $X$.
For the sheaf $\ff = \oo_X(C)$ we define the components of  adelic complex~(\ref{adcompl}) as subgroups of the group
$$
\prod_{p \in D} K_{p,D} \mbox{,}
$$
where this product is taken over all pairs:  irreducible curves $D$ on $X$ and  points $x$ on $D$. For every such pair $x \in D $ the ring $K_{x,D}$ is a finite product of two-dimensional local fields such that if the point $x$ is a smooth point of the curve $D$, then the ring
$K_{x,D}= k(x)((u))((t))$ is a two-dimensional local field (see~\cite[\S 2.2, \S 3.3]{Osi}).
Now
$$\da_{X, 012} (\oo_X(C))= \da_X = {\prod_{x \in D}}' K_{x,D} \mbox{,}$$
where the adelic  product $\prod' f_{x,D}$ is given inside the usual product $\prod f_{x,D}$ by the following two conditions:
\begin{itemize}
\item for a fixed irreducible curve $D \subset X$ and  a local equation $t_D=0$ of the curve $D$ on some open affine subset $U \subset X$ we demand that
$$ \{ f_{x,D}\} \in \da_D((t_D)) \mbox{,}
$$
\item for  all except a finite number of irreducible curves $D \subset X$ we demand that
$$ \{ f_{x,D}\} \in \da_D[[t_D]] \mbox{.}
$$
\end{itemize}

Let $\hat{\oo}_D$ be the completion of local ring of an irreducible curve $D \subset X$, $K_D$ be its fraction field, $\hat{\oo}_x$ be the
completion of local ring of a point $x \in X$, $K_x =\hat{\oo}_x \cdot k(X) $ inside the fraction field of the ring $\hat{\oo}_x$,
$\hat{\oo}_{x,D}$ be the product of rings of discrete valuation of two-dimensional local fields that belong to the ring $K_{x,D}$.

We have the following groups, which are diagonally embedded into the group $\da_X$:
$$\da_{X, 0} (\oo_X(C)) = k(X) \mbox{,}  \qquad \da_{X, 1} (\oo_X(C)) = \prod_{ D \subset X} \hat{\oo}_D \otimes_{\oo_X} \oo_X(C)\mbox{,}
$$
$$
\da_{X, 2} (\oo_X(C)) =\prod_{ x \in  X} \hat{\oo}_x \otimes_{\oo_X} \oo_X(C)\mbox{,}
$$
$$
\da_{X, 01} (\oo_X(C)) ={\prod_{ D \subset  X}}' K_D :=  \left({\prod_{ D \subset  X}} K_D \right) \cap \da_X \mbox{,}
$$
$$
\da_{X, 12} (\oo_X(C)) ={\prod_{ x \in  D}}' \hat{\oo}_{x,D} \otimes_{\oo_X} \oo_X(C) :=
\left(\prod_{ x \in  D} \hat{\oo}_{x,D} \otimes_{\oo_X} \oo_X(C) \right) \cap \da_X \mbox{.}
$$

In the sequel we will omit indication on the sheaf in adelic notations when this sheaf is the structure sheaf $\oo_X$,
and  we will use simple  notations in this case:
$\da_X=\da_{X, 012}$, $\da_{X,0}$, $\da_{X,1}$, $\da_{X,2}$, $\da_{X,01}$, $\da_{X,02}$, $\da_{X,12}$.

The analog of formula~(\ref{analog}) for the case of projective algebraic surface $X$ is the following expression:
\begin{equation}  \label{expad}
\frac{\da_X}{\da_{X, 01} + \da_{X , 02}} = \frac{\da_{X,012}(\oo_X)}{\da_{X, 01}(\oo_X) + \da_{X,02}(\oo_X)} \mbox{.}
\end{equation}
We want to obtain for expression~(\ref{expad}) an exact sequence
which will be analogous to exact sequence~(\ref{adseq}), and also to
the sequences from example~\ref{exad}  for the case of arithmetical
surfaces.

\begin{Th} \label{theorem3}
Let $X$ be a projective smooth irreducible algebraic surface over a field $k$. Let $C_i$, $1 \le i \le k$ be irreducible curves on $X$.
We suppose that $C = \bigcup\limits_{1 \le i \le k} C_i$ is an ample divisor on $X$.
Then there is the following exact sequence:
\begin{equation}  \label{t3seq}
0 \to
\frac{\mathop{{\prod}'}\limits_{x \in D , D \nsubset C} \hat{\oo}_{x,D}}
{\prod\limits_{D \subset X, D \nsubset C} \hat{\oo}_D + \prod\limits_{x \in X, x \notin C} \hat{\oo}_x }
\to
\frac{\da_X}{\da_{X, 01} + \da_{X, 02}} \to
\frac{\prod\limits_{1 \le i \le k} \mathop{{\prod}'}\limits_{x \in C_i} K_{x, C_i}}{\prod\limits_{1 \le i \le k} K_{C_i}
+ \mathop{\prod'}\limits_{x \in C} B_x }
\to 0 \mbox{,}
\end{equation}
where $B_x$ is a subring of the ring $K_x$ defined as
$$
B_x \eqdef \bigcap_{D \subset X, D \nsubset C} \left(K_x \cap \hat{\oo}_{x,D}    \right) \mbox{.}
$$
(The intersection  $K_x \cap \hat{\oo}_{x,D}$ is taken in the ring $K_{x,D}$.)
\end{Th}
\begin{nt} {\em This theorem is am elaboration of the constructions
introduced in~\cite{Par2, Par3}.}
\end{nt}
\begin{nt} \label{adcomp} {\em
We can rewrite the first term in sequence~(\ref{t3seq}) in the following way:
\begin{equation} \label{comob}
\frac{\mathop{{\prod}'}\limits_{x \in D , D \nsubset C} \hat{\oo}_{x,D}}
{\prod\limits_{D \subset X, D \nsubset C} \hat{\oo}_D + \prod\limits_{x \in X, x \notin C} \hat{\oo}_x }=
\frac{ \prod\limits_{D \subset X, D \nsubset C} \left(\frac{\mathop{\prod'}\limits_{x \in D} \hat{\oo}_{x,D}}{\hat{\oo}_D}  \right) }
{\prod\limits_{x \in X, x \notin C } \hat{\oo}_x    }
\end{equation}
For a fixed irreducible curve $D \subset X$ such that $D \not{\hspace{-3pt} \subset} C$ we have that
$$\mathop{\prod\nolimits'}\limits_{x \in D} \hat{\oo}_{x,D} = \da_D[[t_D]] \mbox{,} \qquad
\hat{\oo}_D = k(D)[[t_D]] \mbox{.}
$$
Hence,
\begin{equation} \label{comob1}
\frac{\mathop{\prod'}\limits_{x \in D} \hat{\oo}_{x,D}}{\hat{\oo}_D}= ( \da_D/k(D)  )[[t]]
\end{equation}

We suppose that the field $k$ is a finite field. Then by exact sequence~(\ref{adseq}) we have that the group $\da_D/k(D)$ has a structure of compact object from category $C_1^{ar}$. More exactly, it is an object from category $C_1^{fin}$. Therefore the group $( \da_D/k(D)  )[[t]]$
has also a structure of compact object from category $C_1^{fin}$. (We have to use the structure of object of category $C_1^{fin}$ which is similar to the Tikhonov topology on the infinite product of  compact topological spaces.) For a point $x \in X$  the group $\hat{\oo}_x$ is a profinite group,
therefore it also has a structure of  compact object from category $C_1^{fin}$. Hence we obtain that the group given by expression~(\ref{comob})  has also a structure of  compact object from category $C_1^{fin}$.
}
\end{nt}
\proof (of theorem~\ref{theorem3}). We define an open subset $U = X \setminus C$.
By condition of theorem, the divisor $C$ on the surface $X$ is an ample divisor. Therefore for the quasicoherent sheaf
$\ff \eqdef \mathop{\Lim\limits_{\lto}}\limits_{n \in \dn} \oo_X(nC)$ we have that
$$
H^1(X, \ff) = 0 \quad \mbox{and} \quad H^2(X, \ff)=0 \mbox{,}
$$
because $H^i(X, \ff) = \mathop{\Lim\limits_{\lto}}\limits_{n \in
\dn} H^i(X, \oo_X(nC))$, and $H^i(X, \oo_X(nC)) = 0$ when $i >0$ and
 $n \in \dn$ is sufficiently large. By
definition (see~\cite{B}, \cite{H}, \cite{PF}, \cite{Osi}), the
adelic groups of quasicoherent sheaves commute with inductive
limits. Hence, we have that inside of the group $\da_{X, 012}(\ff)
=\da_{X,012}=\da_X$ the following properties are satisfied:
$$
\da_{X,1}(\ff) = \da_{X, 12}(\ff) \cap \da_{X, 01}(\ff)
\quad \mbox{and} \quad
\da_{X,2}(\ff) = \da_{X, 12}(\ff) \cap \da_{X, 02}(\ff) \mbox{,}
$$
because, by construction,  these properties are satisfied for the
adelic groups of sheaves $\oo_X(nC)$, $n \in\dn$. Therefore from
$H^1(X, \ff) =0$ and adelic complex~(\ref{adcompl}) we obtain that
inside of the group $\da_X$:
\begin{equation} \label{eqin}
\da_{X, 12}(\ff) \cap (\da_{X, 01} + \da_{X, 02}) = (\da_{X,12}(\ff) \cap \da_{X, 01})+(\da_{X, 12}(\ff) \cap \da_{X, 02}) \mbox{,}
\end{equation}
because $\da_{X, 01} (\ff)=\da_{X, 01}$ and $\da_{X, 02}(\ff) =\da_{X, 02}$. By definition,
\begin{equation} \label{ea1}
\da_{X,12}(\ff) = \mathop{{\prod}'}\limits_{x \in D, D \nsubset C} \hat{\oo}_{x,D} \oplus
\prod\limits_{1 \le i \le k} \mathop{{\prod}'}\limits_{x \in C_i} K_{x, C_i}  \mbox{,}
\end{equation}
\begin{equation}  \label{ea2}
\da_{X, 01} = \mathop{{\prod}'}_{D \subset X, D \nsubset C} K_D \oplus \prod_{1 \le i \le k} K_{C_i} \mbox{,}
\end{equation}
\begin{equation}   \label{ea3}
\da_{X, 02} = \mathop{{\prod}'}_{x \in U} K_x \oplus  \mathop{{\prod}'}_{x \in C} K_x \mbox{.}
\end{equation}
Hence, and using formula~(\ref{eqin}), we obtain that
\begin{equation} \label{eqad1}
\da_{X, 12} (\ff) \cap (\da_{X, 01} + \da_{X, 02}) = \left( \prod_{D \subset X,D \nsubset C} \hat{\oo}_D  \oplus \prod_{1 \le i \le k} K_{C_i} \right)
+ \left( \prod_{x \in U} \hat{\oo}_x  \oplus \mathop{{\prod}'}_{x \in C} B_x \right) \mbox{.}
\end{equation}

We consider now a natural map
\begin{equation}  \label{adel1}
\phi \; :  \; \mathop{{\prod}'}_{x \in D, D \nsubset C}  \hat{\oo}_{x,D} \lto  \frac{\da_X}{\da_{X, 01} + \da_{X,02}}  \mbox{.}
\end{equation}
By definition, we have that
$$
\Ker \phi \: = (\mathop{{\prod}'}_{x \in D, D \nsubset C}  \hat{\oo}_{x,D}) \cap (\da_{X, 01} + \da_{X, 02}) \: \supset \:
\prod_{D \subset X,D \nsubset C} \hat{\oo}_D + \prod_{x \in U} \hat{\oo}_x \mbox{.}
$$
From formula~(\ref{ea1}) we have that $\da_{X,12}(\ff) \supset \mathop{{\prod}'}\limits_{x \in D, D \nsubset C}  \hat{\oo}_{x,D}$.
Therefore
$$\da_{X, 12}(\ff) \cap (\da_{X, 01} + \da_{X, 02}) \: \supset \: \Ker \phi  \mbox{.}$$
Hence, by formula~(\ref{eqad1}) we have that
\begin{equation}  \label{adel2}
\Ker \phi = \prod_{D \subset X, D \nsubset C} \hat{\oo}_D + \prod_{x \in U} \hat{\oo}_x \mbox{,}
\end{equation}
because $\mathop{{\prod}'}\limits_{x \in D, D \nsubset C}  \hat{\oo}_{x,D} \supset \Ker \phi $, and  elements of the group
$\mathop{{\prod}'}\limits_{x \in D, D \nsubset C}  \hat{\oo}_{x,D} $ have only the zero projections to any subgroups $K_{x,C_i}$
($x \in C_i$, $1 \le i \le k$)
of adelic group $\da_X$, but any nonzero element from the group  $\prod\limits_{1 \le i \le k} K_{C_i}$  or the group
$\mathop{{\prod}'}\limits_{x \in C} B_x$
has the nonzero projection to some subgroup $K_{x,C_i}$ ($x \in C_i$, $1 \le i \le k$)
of adelic group $\da_X$.

Now from formulas~(\ref{adel1}) and~(\ref{adel2}) it follows the exactness of the beginning of sequence~(\ref{t3seq}).

From $H^2(X, \ff)=0$ and adelic complex~(\ref{adcompl}) we have that
\begin{equation}  \label{adform}
\da_X = \da_{X,12}(\ff) + \da_{X, 01} + \da_{X, 02} \mbox{.}
\end{equation}
We denote the subgroup $\da^U_{X, 12}(\ff) \eqdef \mathop{{\prod}'}\limits_{x \in D, D \nsubset C}  \hat{\oo}_{x,D} $
of the group $\da_X$. Now formula~(\ref{ea1}) can be rewritten as the following formula:
\begin{equation} \label{reformu}
\da_{X, 12}(\ff) = \da_{X, 12}^U (\ff) \oplus \prod_{1 \le i \le k} \mathop{{\prod}'}_{x \in C_i} K_{x,C_i} \mbox{.}
\end{equation}
From formula~(\ref{adform}) we have the following evident exact sequence:
\begin{equation} \label{eeqq1}
 \frac{\da_{X,12}(\ff) \cap (\da_{X, 12}^U (\ff) + \da_{X, 01} + \da_{X, 02} )}{\da_{X,12}^U (\ff)}
\hookrightarrow \frac{\da_{X,12}(\ff)}{\da_{X,12}^U(\ff)}
\twoheadrightarrow
\frac{\da_X}{\da_{X,12}^U(\ff) + \da_{X,01} + \da_{X,02}}  \mbox{.}
\end{equation}
We denote the group
$$G \: \eqdef \: \frac{\da_{X,12}(\ff) \cap (\da_{X, 01} + \da_{X, 02} )}{ \da_{X, 12}^U(\ff) \cap ( \da_{X, 01} + \da_{X, 02} ) }  \mbox{,}$$
and we
 consider a map:
 \begin{equation} \label{eeqq2}
 \psi \: : \:
G \lto \frac{\da_{X,12}(\ff) \cap (\da_{X, 12}^U (\ff) + \da_{X, 01} + \da_{X, 02} )}{\da_{X,12}^U (\ff)} \mbox{.}
\end{equation}
It is clear that the map $\psi$ is an isomorphism.

We want to prove the following isomorphism of groups:
\begin{equation} \label{formadel}
G \simeq \prod_{1 \le i \le k} K_{C_i} + \mathop{{\prod}'}_{x \in C} B_x \mbox{.}
\end{equation}
By definition,
\begin{equation} \label{fad}
\da_{X, 12}^U (\ff)  \cap (\da_{X, 01} + \da_{X, 02})  \:  \supset  \:
\prod_{D \subset X , D \nsubset C} \hat{\oo}_D + \prod_{x \in U} \hat{\oo}_x  \mbox{.}
\end{equation}
From formula~(\ref{eqad1}) we have a natural map
$$
\xi \: : \: \prod_{1 \le i \le k} K_{C_i} + \mathop{{\prod}'}_{x \in C_i} B_x  \lto G \mbox{.}
$$
From formulas~(\ref{eqad1}) and~(\ref{fad}) we obtain that the map $\xi$ is a surjective map. The map $\xi$ is an injective map,
because the image of map $\xi$ does not contain the subgroup $\da_{X,12}^U (\ff)$. Thus, the map $\xi$ is an isomorphism, and we have proved
formula~(\ref{formadel}).

From formula~(\ref{reformu}) we have that
$$
\frac{\da_{X,12}(\ff)}{\da_{X,12}^U(\ff)} \simeq \prod_{1 \le i \le k} \mathop{{\prod}'}_{x \in C_i} K_{x,C_i} \mbox{.}
$$
Hence and from formulas~(\ref{eeqq1})--(\ref{formadel}), we obtain
the exactness of the middle and the end of sequence~(\ref{t3seq}).
The theorem is proved.

\subsection{More refined computation} \label{mrc}
In this subsection we will compute the last nonzero term of
sequence~(\ref{t3seq}) more explicitly.

\begin{Th} \label{te4}
Let $X$ be a projective smooth irreducible algebraic surface over a field $k$. Let $C_i$, $1 \le i \le k$ be irreducible curves on $X$.
Let $C = \bigcup\limits_{1 \le i \le k} C_i$. Let the point $p = \bigcap\limits_{1 \le i \le k} C_i$. There is the following isomorphism:
$$
\frac{\prod\limits_{1 \le i \le k} \mathop{{\prod}'}\limits_{x \in C_i} K_{x, C_i}}{\prod\limits_{1 \le i \le k} K_{C_i}
+ \mathop{\prod'}\limits_{x \in C} B_x } \: \simeq \:
\frac{\prod\limits_{1 \le i \le k} K_{p,C_i}}{\prod\limits_{1 \le i \le k} B_{C_i} + B_p}  \mbox{,}
$$
where the ring $B_{C_i}$ ($1 \le i \le k$) is a subring of the field $K_{C_i}$ defined as
$$
B_{C_i}  \eqdef \bigcap_{x \in C_i \setminus p } (K_{C_i} \cap B_x) \mbox{.}
$$
(The intersection $K_{C_i} \cap B_x$ is taken in the ring $K_{x,C_i}$.)
\end{Th}
Before the proof of this theorem we will prove the following lemma.
\begin{lemma} \label{adlemm}
Under conditions of theorem~\ref{te4}, for every $1 \le i \le k$ there is the following isomorphism:
$$
K_{C_i}/B_{C_i}  \simeq  \mathop{{\prod}'}_{x \in C_i \setminus p} K_{x,C_i}/B_x \mbox{.}
$$
\end{lemma}
\proof of lemma~\ref{adlemm}. We fix $1 \le i \le k$. We have the diagonal embedding: $K_{C_i} \to   \mathop{{\prod}'}\limits_{x \in C_i \setminus p} K_{x,C_i}$,
which induces the map:
$$
\eta \: : \: K_{C_i}/B_{C_i}  \lto  \mathop{{\prod}'}_{x \in C_i \setminus p} K_{x,C_i}/B_x \mbox{.}
$$
Now the injectivity of the map $\eta$ follows from the definition of the ring $B_C$.

We will prove that the map $\eta$ is a surjective map.

Let $J_{C_i}$ be the ideal sheaf of the curve $C_i$ on the surface $X$. For every $n \in \dn$ we consider a scheme
$$Y_n =(C_i \setminus p \:,\: \oo_X/ J_{C_i}^n \mid_{ C_i \setminus p}) \mbox{,}$$
which is some infinitesimal  neighbourhood of the affine curve $C_i \setminus p$ in $X$. For every $l, m \in \dz$ such that $l < m$  we consider a coherent sheaf
$J_{C_i}^l /J_{C_i}^m \mid_{C_i \setminus p}$  on the scheme $Y_{m-l}$. We have that $\dim Y_{m-l} =1$, and $H^1(Y_{m-l},J_{C_i}^l /J_{C_i}^m
\mid_{ C_i \setminus p})=0$ for any $l < m \in \dz$, because $Y_{m -l}$ is an affine scheme. Therefore the adelic complex (see~\cite{B},~\cite{Osi}) gives us the following exact sequence $\kk_{l,m}$ for every
$l < m \in \dz$:
$$
0 \lto H^0(Y_{m-l},J_{C_i}^l /J_{C_i}^m \mid_ {C_i \setminus p})  \lto \da_{Y_{m-l}, 0}(J_{C_i}^l /J_{C_i}^m \mid _{C_i \setminus p}) \oplus \da_{Y_{m-l}, 1}(J_{C_i}^l /J_{C_i}^m \mid _{C_i \setminus p})
\qquad \qquad \qquad \qquad \qquad
$$
$$
  \quad \qquad \qquad \qquad \qquad \qquad \qquad \qquad \qquad \qquad \qquad
\lto \da_{Y_{m-l}, 01}(J_{C_i}^l /J_{C_i}^m \mid _{C_i \setminus p}) \lto 0 \mbox{.}
$$
For every fixed $l \in \dz$ the projective system of abelian groups
$$(H^0(Y_{m-l},J_{C_i}^l /J_{C_i}^m \mid_ {C_i \setminus p}) \: , \: m>l)$$ satisfies the
ML-condition (the Mittag-Leffler condition), because for any $m > l$
the map
$$H^0(Y_{m-l+1},J_{C_i}^l /J_{C_i}^{m+1} \mid _ {C_i \setminus p})  \lto H^0(Y_{m-l},J_{C_i}^l /J_{C_i}^m \mid _{C_i \setminus p}) $$
is a surjective map, as it follows
from the following  exact sequence of sheaves on the scheme $Y_{m-l+1}$:
$$
0 \lto J_{C_i}^m /J_{C_i}^{m+1} \mid_{C_i \setminus p} \lto J_{C_i}^l /J_{C_i}^{m+1} \mid _ {C_i \setminus p}
\lto J_{C_i}^l /J_{C_i}^m \mid _{C_i \setminus p} \lto 0 \mbox{,}
$$
and $H^1(Y_{m-l+1},J_{C_i}^m /J_{C_i}^{m+1} \mid_{C_i \setminus p})=0$, because  $J_{C_i}^m /J_{C_i}^{m+1}$ is a coherent sheaf on the affine scheme $Y_{m-l+1}$.

Now we consider the sequence $\kk :=
\mathop{\Lim\limits_{\lto}}\limits_{l \in \dz}
\mathop{\Lim\limits_{\longleftarrow}}\limits_{m > l}  \kk_{l,m}$,
which is an exact sequence, because the inductive limit preserves
the exactness of sequences, and the projective limit preserves the
exactness due to ML-condition which we have explained above. The
sequence $\kk$ is of the following type:
\begin{equation} \label{dad}
0 \lto M \lto \da_0 \oplus \da_1 \lto \da_{01} \lto 0 \mbox{,}
\end{equation}
where as subgroups of the group $\da_X$ we have
$$
\da_0= K_{C_i}  \mbox{,} \qquad
\da_1= \mathop{{\prod}'}_{x \in C_i \setminus p} B_x \mbox{. } \qquad
\da_{10}= \mathop{{\prod}'}_{x \in C_i \setminus p} K_{x,C_i} \mbox{,}
$$

  Now from exactness of sequence~(\ref{dad}) we obtain that the map $\eta$ is a surjective map. The lemma is proved.

\vspace{0.4cm}
\noindent \proof (of theorem~\ref{te4}).
Let $k=1$, i.e, $C=C_1$. We have
$$
\frac{ \mathop{{\prod}'}\limits_{x \in C} K_{x, C}}{ K_{C}
+ \mathop{\prod'}\limits_{x \in C} B_x }
=
\frac{ \mathop{{\prod}'}\limits_{x \in C\setminus p} K_{x, C} \oplus K_{p,C}}{ K_{C}
+ \mathop{\prod'}\limits_{x \in C\setminus p} B_x \oplus B_p }
=
\frac{ \frac{\mathop{{\prod}'}\limits_{x \in C\setminus p} K_{x, C}}{ \mathop{\prod'}\limits_{x \in C\setminus p} B_x }
\oplus \frac{K_{p,C}}{B_p}}{ K_{C}
}
$$
$$
\stackrel{\mbox{by lemma} \; \ref{adlemm}}{=}  \:
\frac{ \frac{K_C}{ B_C }
\oplus \frac{K_{p,C}}{B_p}}{ K_{C}
}=
\frac{ \frac{K_C}{ B_C }
\oplus \frac{K_{p,C}}{B_p+B_C}}{ \frac{K_{C}}{B_C}
}
=
\frac{ K_{p,C}}{ B_{C} + B_p} \mbox{.}
$$

If $k >1$, then we have to do the analogous computation applying lemma~{\ref{adlemm}} $k$ times: for the  curves $C_i$, $1 \le i \le k$.
For example, we do it when $k=2$. In this case we have
$$
\frac{\mathop{{\prod}'}\limits_{x \in C_1} K_{x,C_1} \oplus \mathop{{\prod}'}\limits_{x \in C_2} K_{x,C_2}}{K_{C_1} \oplus K_{C_2}
+ \mathop{{\prod}'}\limits_{x \in C} B_x} \:
= \:
\frac{
(\mathop{{\prod}'}\limits_{x \in C_1 \setminus p} K_{x,C_1}) \oplus K_{p,C_1} \oplus
(\mathop{{\prod}'}\limits_{x \in C_2 \setminus p} K_{x,C_2}) \oplus K_{p,C_2}
}{
K_{C_1} \oplus K_{C_2} + (\mathop{{\prod}'}\limits_{x \in C_1 \setminus p} B_x) \oplus B_p \oplus
(\mathop{{\prod}'}\limits_{x \in C_2 \setminus p} B_x)
}
$$
$$
\stackrel{\mbox{by lemma} \; \ref{adlemm}}{=}  \:
\frac{K_{p,C_1} \oplus K_{p,C_2}}{B_{C_1} \oplus B_{C_2} + B_p} \mbox{.}
$$
Theorem~\ref{te4} is proved.

\subsection{Case of an arithmetic surface}
Let $X = C_1 \times_k C_2$, and $C = C_1  \cup C_2$, $p = C_{1} \cap C_{2}$. We have
 the transversal intersection of the curves $C_1$ and $C_2$ on $X$.
Assume that $C_1$ and $C_2$ are projective lines over the field $k$ with
coordinates  $t$ and $u$ such that the point $p$ is given by local equations $t=0$, $u=0$.

By definition, we have the following rings and fields:
$$
k(C_1) = k(t)
\mbox{,} \qquad
k(C_2) = k(u) \mbox{,}
$$
$$
K_{p,C_1} = k((t))((u))
\mbox{,}
\qquad
K_{p,C_2} = k((u))((t))  \mbox{,}
$$
$$
B_{C_1} = k[t^{-1}]((u))
\mbox{,}
\qquad
B_{C_2} = k[u^{-1}]((t)) \mbox{,}
$$
$$
B_p  = \mathop{\Lim_{\lto}}_{m \in \dn, n \in \dn} u^{-m}t^{-n} k[[u, t]] \mbox{.}
$$

We note that the quotient group (see theorem~\ref{te4})
$$
F = \frac{K_{p,C_1} \oplus K_{p,C_2}}{ B_{C_1} \oplus  B_{C_2} +     B_p} \mbox{,}
$$
where the group $B_p$ is diagonally embedded to the group
$K_{p,C_1} \oplus K_{p,C_2}$, can be computed by reduction to a single two-dimensional
local field, say, the field  $K_{p,C_2}$. Namely, we have
that $ 0 \oplus K_{p,C_2} + B_{C_1} \oplus  B_{C_2} +     B_p=K_{p,C_1} \oplus K_{p,C_2}$. Therefore
\begin{equation} \label{factsec}
0 \lto K_{p,C_2}  \cap (B_{C_1} \oplus  B_{C_2} +     B_p)  \lto  K_{p,C_2} \rightarrow  F \lto 0
\mbox{,}
\end{equation}
where the intersection is taken inside the group $K_{p,C_1} \oplus
K_{p,C_2}$, and $K_{p,C_2} = 0 \oplus K_{p,C_2}$ is considered
inside of this group. The first nonzero group in
sequence~(\ref{factsec}) is equal to the group $B_{C_2} + K_{p,C_2}
\cap (B_{C_1} +  B_p)$, since the group $B_{C_2}$  belongs only to
the group $K_{p,C_2}$ as a subgroup. The group $K_{p,C_2}  \cap
(B_{C_1} +     B_p)$ is equal to the subgroup $B_{C_1} \cap  B_p$ of
the group $K_{p,C_1} \oplus 0$, which is embedded  to the group $0
\oplus K_{p,C_2}$ through the inclusion $B_p \subset  K_{p,C_2}$.

At last, we obtain
$$
F =  \frac{K_{p,C_2}}{B_{C_2} +  (B_{C_1} \cap B_{p})} =
\frac{k((u))((t))}{ k[u^{-1}] ((t)) + k((u))[t^{-1}] } \simeq k[[u, t]]ut \mbox{.}
$$

We have the following analogy of this computation for an arithmetic
surface. We consider  the surface $X$ as a fibration with projection
onto the curve $C_2$. We want to compare the surface $X$ fibered
over the curve $C_2$ with the simplest arithmetic surface $P^1$ over
$\Spec \mathbb{Z}$. Point $p \in C_2$ will correspond to the place
$\infty$ added to $\Spec \mathbb{Z}$. The curve  $C_1$ corresponds
to the non-existing closed fiber over the place $\infty$. There are
the following analogies based on the classical analogy between
algebraic surfaces and arithmetic surfaces (see, for
example,~\cite{Par}):
$$
k(C_2)_p  = k((u))   \sim    \mathbb{R}
$$
$$
K_{p, C_2}    \sim    \mathbb{R}((t))
$$
$$
B_{C_2}  =  \{\mbox{regular functions on} \; \, C_2 \setminus p\} ((t))  \sim \mathbb{Z}((t))
$$
$$
B_{C_1} \cap B_p  = k((u))[t^{-1}]  \sim \mathbb{R}[t^{-1}]
$$

Using this dictionary one obtains the following analog of the group $F$:
$$
F \sim \frac{\mathbb{R} ((t))}{ \mathbb{Z} ((t))  +  \mathbb{R} [t^{-1}]}      =
(\mathbb{R} / \mathbb{Z} ) [[t]]t \simeq \dt[[t]] t \mbox{.}
$$

Obviously, on the group $\dt[[t]] t$ we have a structure of an object from the category $C_2^{ar}$  such that it
is a $c$-object and a $cf$-object.

\begin{nt}{\em
It is also interesting to consider an analog of the first nonzero
term in sequence~(\ref{t3seq})  for the case of an arithmetical
surface. We recall the situation of example~{\ref{arsur}}. We
consider a regular two-dimensional scheme $X$ with a projective
surjective morphism $X \lto \Spec E$, where $E$ is the ring of
integers in a number field $K$, $[K : \dq] =n$.  It means that $X$
is an arithmetical surface. Let $X_K$ be the generic fibre of this
morphism. Let $C \subset X$ be a "horizontal" irreducible arithmetic
curve, i.e.  an integral one-dimensional subscheme $C$ which is
surjectively  mapped onto $\Spec E$. The subscheme $C$ corresponds
to some closed point $p_C \in X_K$. We recall that $p_1, \ldots,
p_l$ are all Archimedean places of the field $K$, and $K_{p_i}$, $1
\le i \le l$ is the corresponding completion field. Then according
to the definition of adelic ring $\da_X^{ar}$ (see
example~{\ref{arsur}}) we have that the following group $\Psi$ is an
analog of the first nonzero term in sequence~(\ref{t3seq}):
\begin{equation} \label{hex}
\Psi \: \eqdef \: \frac{(\mathop{{\prod}'}\limits_{x \in D , D
\nsubset C} \hat{\oo}_{x,D}) \; \oplus \prod\limits_{q \in X_K, \, q
\ne p_C} \hat{\oo}_q \hat{\otimes}_K (\prod\limits_{1 \le i \le l}
K_{p_i}) } {\prod\limits_{D \subset X, D \nsubset C} \hat{\oo}_D +
\prod\limits_{x \in X, x \notin C} \hat{\oo}_x } \mbox{,}
\end{equation}
where $x \in D$ runs over the pairs: integral one-dimensional
subschemes $D$ of $X$ and closed points $x$ on $D$. A point $q$ is a
closed point of the curve $X_K$, the ring $\hat{\oo}_q$ is the
completion of the local ring of the point $q$ on $X_K$ with respect
to the maximal ideal of the point $q$. By definition,
$$
\hat{\oo}_q \hat{\otimes}_K (\prod\limits_{1 \le i \le l} K_{p_i})
\;=\; \mathop{\mathop{\Lim}\limits_{\longleftarrow}}\limits_{n \in
\dn} \left(\frac{\hat{\oo}_{q}}{\hat{m}_q^n}  \otimes_K (
\prod\limits_{1 \le i \le l} K_{p_i} ) \right) \mbox{,}
$$
where $\hat{m}_q$ is the maximal ideal of the local ring
$\hat{\oo}_q$.

Now, we clarify how subgroups from (\ref{hex}) are
embedded, one inside another. We note that for a fixed
"horizontal" arithmetic irreducible curve $D \subset X$ there is  the
diagonal embedding:
$$
\hat{\oo}_D \; \hookrightarrow \; (\mathop{{\prod}'} \limits_{x \in
D} \hat{\oo}_{x,D}) \; \oplus \; \hat{\oo}_{p_D} \hat{\otimes}_K (
\prod\limits_{1 \le i \le l} K_{p_i} ) \mbox{,}
$$
where the point $p_D \in X_K$ corresponds to the irreducible
arithmetic curve $D \subset X$, and, by construction, $ \hat{\oo}_D
= \hat{\oo}_{p_D}$.

For a fixed "vertical" irreducible curve $D \subset X$ (i.e. $D$ is
defined over some finite field) we have only the following diagonal
embedding:
$$
\hat{\oo}_D \; \hookrightarrow \; \mathop{{\prod}'} \limits_{x \in
D} \hat{\oo}_{x,D} \mbox{.}
$$

Also, for a fixed closed point $x \in X$ we have only the following
diagonal embedding:
$$
\hat{\oo}_x \hookrightarrow \prod\limits_{D \ni x} \hat{\oo}_{x,D}
\mbox{.}
$$

Now, similar to~ (\ref{comob}) we can rewrite the group $\Psi$
in the following way:
\begin{equation} \label{finfin}
\Psi = \frac{\Psi_1 \oplus \Psi_2}{\prod\limits_{x \in X, x \notin
C} \hat{\oo}_x } \mbox{,}
\end{equation}
where
$$
\Psi_1 = \mathop{\prod\limits_{D \subset X, D \nsubset C}}\limits_{D
"horizontal"} \left( \frac{\mathop{\prod'}\limits_{x \in D}
\hat{\oo}_{x,D} \; \oplus \; \hat{\oo}_{p_D} \hat{\otimes}_K (
\prod\limits_{1 \le i \le l} K_{p_i} ) }{\hat{\oo}_D} \right)
\mbox{,}
$$
and
$$
\Psi_2 = \mathop{\prod\limits_{D \subset X, D \nsubset C}}\limits_{D
"vertical"} \left( \frac{\mathop{\prod'}\limits_{x \in D}
\hat{\oo}_{x,D}  }{\hat{\oo}_D} \right) \mbox{.}
$$

For a fixed irreducible "horizontal" arithmetic curve $D \subset X$
let $k(D)$ be the field of rational functions on $D$, and let $t_D=0$ be
the local equation of $D$ on some open subset of $X$.  We have
$$
\frac{\mathop{\prod'}\limits_{x \in D} \hat{\oo}_{x,D} \; \oplus \;
\hat{\oo}_{p_D} \hat{\otimes}_K ( \prod\limits_{1 \le i \le l}
K_{p_i} )}{\hat{\oo}_D} =
\left(\da_{k(D)}/k(D)\right)[[t_D]] \mbox{.}
$$
From example~\ref{exad} we have that $\da_{k(D)}/k(D)$ is a compact
object in category $C_1^{ar}$. Therefore,
$\left(\da_{k(D)}/k(D)\right)[[t_D]]$ and, consequently, $\Psi_1$
has structure of $c$-object and $cf$-object from category $C_2^{ar}$
(compare with remark~\ref{adcomp}).

For a fixed irreducible "vertical"  curve $D \subset X$   we have
$$
\frac{\mathop{\prod'}\limits_{x \in D} \hat{\oo}_{x,D}
}{\hat{\oo}_D} =
\mathop{\mathop{\Lim}\limits_{\longleftarrow}}\limits_{n \in \dn}
\da_{D_n} / \da_{D_n, 0} \mbox{,}
$$
where the closed subscheme $D_n \subset X$ is given as
$$
D_n = (D, \oo_X/J_D^n) \mbox{,}
$$
and $J_D$ is the ideal sheaf of the curve $D$ on $X$. We note that
$D_1=D$. By induction on $n$ and by reasoning from the beginning of
section~\ref{sqr} we obtain that $\da_{D_n} / \da_{D_n, 0}$ has
structure of a compact object from category $C_1^{ar}$. Hence
$\Psi_2$ has a structure  of $c$-object and $cf$-object from
category $C_2^{ar}$.

Like in remark~\ref{adcomp}, the group $\prod\limits_{x \in X, x
\notin C} \hat{\oo}_x $ has a structure of $c$-object and
$cf$-object from category $C_2^{ar}$.

Hence and from~(\ref{finfin}), we obtain that the group
$\Psi$ has a structure of $c$-object and $cf$-object from category
$C_2^{ar}$ (compare with remark~\ref{adcomp}).

}
\end{nt}

\vspace{0.3cm}

\noindent D.V. Osipov\\
Steklov Mathematical Institute RAS \\
{\it E-mail:}  ${d}_{-} osipov@mi.ras.ru$ \\

\vspace{0.3cm}

\noindent A.N. Parshin\\
Steklov Mathematical Institute RAS \\
{\it E-mail:}  $parshin@mi.ras.ru$
\end{document}